\documentclass[11pt, a4paper]{article}

\RequirePackage{amsthm,amsmath,amsfonts,amssymb}
\RequirePackage[authoryear]{natbib} 
\RequirePackage[colorlinks,citecolor=blue,urlcolor=blue,breaklinks]{hyperref}

\usepackage{bm}
\usepackage{mathrsfs}  
\usepackage{appendix}
\usepackage{commath}
\usepackage{dsfont}
\usepackage{subfig}
\usepackage{xcolor}
\usepackage{mathtools}
\usepackage{enumitem}
\usepackage{fancyhdr} 
\usepackage[a4paper, margin=1in, headheight=14pt]{geometry}
\usepackage{caption}
\usepackage{booktabs}
\usepackage{multirow}
\usepackage{lscape}

\usepackage{setspace}
\newcommand{\N}{\mathbb{N}}
\newcommand{\R}{\mathbb{R}}
\newcommand{\eps}{\varepsilon}
\renewcommand{\P}{\mathbb{P}}
\newcommand\E{\mathbb{E}}
\newcommand{\mc}{\mathcal}

\theoremstyle{plain}
\newtheorem{theorem}{Theorem}[section]
\newtheorem{lemma}[theorem]{Lemma}

\newtheorem{proposition}[theorem]{Proposition}

\theoremstyle{remark}
\newtheorem{remark}{Remark}[section]

\title{A sharp dimension barrier for uniform Gaussian approximations}

\author{
\begin{tabular}{c}
Anders Bredahl Kock\footnote{Kock's research was supported by the European Research Council (ERC) grant number 	 -- HIDI (UKRI EP/Z002222/1).  He is also a member of, and grateful for support from, i) the Aarhus Center for Econometrics (ACE), funded by the Danish National Research Foundation grant number DNRF186,  and ii) the Center for Research in Energy: Economics and Markets (CoRe).} \\ 
\footnotesize	University of Oxford \\
\footnotesize	10 Manor Rd, Oxford OX1 3UQ
\\
\footnotesize	{\footnotesize	\href{mailto:anders.kock@economics.ox.ac.uk}{anders.kock@economics.ox.ac.uk}} 
\end{tabular}
}

\date{\today}

\begin{document}
\maketitle

\begin{abstract}
We establish a universal dimension barrier for Gaussian approximations based on~$n$ i.i.d.~random vectors~$X_1,\hdots,X_n$ in~$\R^d$. Fix~$m\in(2,\infty)$ and consider mean-zero distributions with coordinate variances bounded from below by a fixed positive constant and coordinatewise~$m$th absolute moments bounded from above by a fixed constant. We show that
\begin{equation*}
	\log (d)=o\del[1]{n^{(m-2)/(3m-2)}}
\end{equation*}
is necessary for \emph{any} sequence of statistics to admit a Gaussian approximation to its distribution uniformly over this distribution class and all hyperrectangles. Here, the centered Gaussian target is required to have the covariance or correlation matrix of the observations. The obstruction includes randomized procedures and already holds for coordinate maxima. Conversely, for~$m\in(2,4]$, we establish that this growth condition on~$d$ is sufficient for centered and scaled quantile-winsorized means and their suitably studentized counterparts to admit uniform Gaussian approximations with covariance and correlation targets, respectively, without further restrictions on the covariance or knowledge of~$m$. Thus, the dimension barrier is sharp in this moment range. For~$m\in(4,\infty)$, both procedures attain the barrier up to the same power of~$\log (n)$.
\end{abstract}

\section{Introduction}\label{sec:Intro}
In this article we establish a universal dimension barrier for uniform Gaussian approximation and show that it can be attained by winsorized means, whereas standard statistics such as ordinary sums and self-normalized sums generally face stricter dimension restrictions.

Let~$X_1,\ldots,X_n$ be i.i.d.~random
vectors in~$\R^d$ with distribution~$P$, mean zero, and covariance
matrix~$\Sigma=\Sigma(P)$. Write~$\gamma_{A}:=\mathsf N_d(0_d,A)$ for the centered Gaussian
probability measure on~$\R^d$ with covariance matrix~$A$ (which need not
be nonsingular) and let $\mc H_d$ be the class of (generalized) hyperrectangles in~$\R^d$, that is all sets of the form
\begin{equation*}
	H=\cbr[1]{x\in\R^d:a_j\leq x_j\leq b_j\text{ for all }j=1,\hdots, d},
\end{equation*}
where~$-\infty\leq a_j\leq b_j\leq \infty$ for all~$j = 1, \hdots, d$. Throughout, $d=d_n$ may increase with~$n$. Writing~$P^{\otimes n}$ for the joint (product) distribution
of~$X_1,\hdots,X_n$, the Gaussian approximation error to the distribution of a statistic~$T_n:(\R^d)^n\to\R^d$ over~$\mc{H}_d$ is
\begin{equation}\label{eq:intro:sumerror}
 \rho_n(T_n,A):=\sup_{H\in\mc H_d}
 \left|P^{\otimes n}(T_n\in H)
       -\gamma_{A}(H)\right|.
\end{equation}
Since the seminal paper of~\cite{chernozhukov2013gaussian}, increasingly refined upper bounds on~$\rho_n(S_n,\Sigma)$ and related quantities have been established for~$S_n=n^{-1/2}\sum_{i=1}^nX_i$, cf., e.g.,
\cite{chernozhukov2017central,deng2020beyond,
kuchibhotla2020high,koike2021notes,kuchibhotla2021high,
lopes2022central,chernozhuokov2022improved,fang2023high,koike2024high, gao2026cubic, LopesMasini26, CCKK26}.\footnote{These papers do not always impose that the~$X_i$ are identically distributed, but the lower bounds in Theorem~\ref{thm:sharp} hold already in this setting.} We refer to
\cite{chernozhukov2023high} for a review and further references. Additional structure, such as variance decay or suitable covariance nondegeneracy, permits improvements; see~\cite{lopes2020bootstrapping,das2021central,fang2021high,chernozhukov2023nearly}. For bounded summands with uniformly nondegenerate covariance, see~\cite{aronow2026optimal}.

Despite the progress on bounds on~$\rho_n(S_n,\Sigma)$ when the~$X_{i,j}$ have relatively light tails, it follows from Remark~2 in~\cite{zhangwu2017} and Theorem 2.1 in~\cite{kock2024remark} that for every~$m\in(2,\infty)$ there exist~$X_i$ with i.i.d.~coordinates having mean zero, variance one and bounded~$m$th absolute moment, yet~$\limsup_{n\to\infty}\rho_n(S_n,\Sigma)=1$ if~$\limsup_{n\to\infty}\frac{d}{n^{m/2-1+\xi}}>0$ for some~$\xi\in(0,\infty)$. In particular,~$\rho_n(S_n,\Sigma) \to 0$ does not hold uniformly over all distributions specified by such moment bounds when~$d$ exceeds this \emph{polynomial} growth barrier. \footnote{Further moment-dependent necessary and sufficient conditions for~$\rho_n(S_n,\Sigma)\to 0$ under i.i.d.~coordinates have been studied in~\cite{das2025necessary}.} 

The polynomial growth barrier faced by~$S_n$ can be overcome by changing the statistic~$S_n$, while keeping~$\gamma_{\Sigma}$ as the approximating Gaussian. \cite{resende2024robust} established Gaussian approximations for maxima of trimmed means under~$m\in(2,\infty)$ moments and~\cite{kock2025high} obtained approximations over all of~$\mc{H}_d$ for winsorized and trimmed means. If the variances are bounded away from zero and the absolute centered~$m$th moments are bounded from above, the Gaussian approximation errors in these papers converge to zero if~$\log(d)=o\del[1]{n^{(m-2)/(5m-2)}}$ without any further restrictions on the covariance matrix~$\Sigma$. This exponential rate remains valid also in the presence of a sufficiently quickly vanishing fraction of adversarially contaminated~$X_i$. Changing the Gaussian target to~$\gamma_{\Sigma_0}$, where~$\Sigma_0=\Sigma_0(P)$ denotes the correlation matrix of the~$X_i$, studentized winsorized means likewise admit approximations under the same conditions, cf.~Theorem~4.1 of~\cite{kock2025high}. Robust max-statistics under additional variance-decay, moment, and dependence conditions have been studied by~\cite{liu2024robust}.

Another approach to accommodating larger dimensions is to self-normalize the entries of~$S_n$. Let~$U_n=(U_{n,1},\ldots,U_{n,d})'$ with
\begin{equation}\label{eq:selfnorm}
U_{n,j}
=
\frac{\sum_{i=1}^nX_{i,j}}
{\bigl(\sum_{i=1}^nX_{i,j}^2\bigr)^{1/2}},
\qquad j=1,\ldots,d.
\end{equation}
The relevant Gaussian target is now~$\gamma_{\Sigma_0}$. If the~$X_{i,j}$ are i.i.d.~for each~$i$ (so that~$\Sigma_0=I_d$) with uniformly bounded~$m$th moments and variances bounded away from zero, \cite{das2024central} showed that~$\rho_n(U_n,I_d)\to 0$ if~$\log(d)=o(n^{(m-2)/m})$ for~$2<m\leq 3$, while \cite{qiu2025self} established relative Gaussian tail approximations under suitable restrictions on~$\Sigma_0$.\footnote{See also \cite{chang2025berry} for Berry--Esseen bounds for absolute maxima of self-normalized sums.}

 Under \emph{unrestricted} $\Sigma_0$, \citet[Theorem~3.1]{wu2026counterexample} shows that Gaussian approximation for maxima of self-normalized sums can fail already at polynomial dimension~$d=\lceil n^m\rceil $, under unit coordinate variances and uniformly bounded coordinatewise~$m$th absolute moments for any fixed~$m>2$. Recently, \citet[Theorem 2.1]{su2026sharp} showed that, under unrestricted~$\Sigma_0$,~$\rho_n(U_n,\Sigma_0)\to 0$ uniformly over this class of distributions if and only if~$d=o\del[0]{n^{m/2-1}}$ for~$m\geq 3$.\footnote{\cite{su2026sharp} appeared
after the first version of the present article.} The exponent~$m/2-1$ on~$n$ is the same as in the polynomial obstruction for ordinary sums~$S_n$ described above. Thus, both sums and self-normalized sums can fail in dimension regimes where winsorized means admit uniform Gaussian approximations without additional covariance restrictions.

These results show that a dimension barrier for one statistic need not
be a barrier for another statistic. It is therefore natural to ask
\begin{quote}
\emph{Is there a universal dimension barrier that cannot be overcome by changing the statistic, while retaining $\gamma_{\Sigma}$ or $\gamma_{\Sigma_0}$ as the Gaussian target?}
\end{quote}
Existing lower bounds for particular statistics, such as sums or self-normalized sums, do not answer this question. Our main result, Theorem~\ref{thm:sharp}, shows that, for every fixed~$m\in(2,\infty)$ and~$\Gamma\in\cbr[0]{\Sigma,\Sigma_0}$, the condition
\begin{equation}\label{eq:barrierintro}
	 \log(d)=o\del[1]{n^{(m-2)/(3m-2)}}
\end{equation}
is necessary for \emph{any} sequence of measurable statistics~$T_n:(\R^d)^n\to\R^d$, not depending on the unknown~$P$, to satisfy~$\rho_n(T_n,\Gamma)\to 0$ uniformly over the mean-zero moment class with coordinate variances bounded away from zero and~$m$th absolute moments bounded above. This impossibility remains valid for randomized statistics represented by Markov kernels and the obstruction already occurs for approximation of the distribution of~$\max_{1\leq j\leq d}T_{n,j}$. 

Conversely, Theorems~\ref{thm:SufficientWins} and~\ref{thm:CorrTargetSharp} show that, for every fixed~$m\in(2,4]$, the quantile-winsorized means~$S_{n,W}$ and their suitably studentized counterparts~$S_{n,W,S}$ studied in~\cite{kock2025high} satisfy~$\rho_n(S_{n,W},\Sigma)\to 0$ and~$\rho_n(S_{n,W,S},\Sigma_0)\to 0$, respectively, uniformly over this moment class if~\eqref{eq:barrierintro} is satisfied. No further assumptions on~$\Sigma$ or~$\Sigma_0$ are imposed and the winsorization rules do not require knowledge of~$m$. Thus, for~$m\in(2,4]$ and either Gaussian target, condition~\eqref{eq:barrierintro} characterizes exactly the dimension growth for which a uniform Gaussian approximation is possible.

For~$m\in(4,\infty)$, our sufficient condition on~$d$ for~$\rho_n(S_{n,W},\Sigma)\to 0$ and $\rho_n(S_{n,W,S},\Sigma_0)\to 0$ differs from the necessary condition~\eqref{eq:barrierintro} only by a power of~$\log n$. Thus, the necessary and sufficient conditions on~$\log(d)$ have the same exponent~$(m-2)/(3m-2)$ on~$n$ for every fixed~$m\in(2,\infty)$, allowing for a larger growth of~$d$ than the requirement~$\log(d)=o\del[1]{n^{(m-2)/(5m-2)}}$ in~\cite{kock2025high}.

Theorems~\ref{thm:HDGaussSharp} and~\ref{thm:HDGauss_studentized} also establish fixed-$n$ approximations to the distribution of~$S_{n,W}$ and~$S_{n,W,S}$, respectively, in the presence of some coordinatewise adversarial contamination.
 
The proof of the lower bound constructs a baseline distribution of the~$n$ observations that is close in total variation to a mixture of alternative distributions, while the corresponding Gaussian maximum probabilities remain separated. Rare perturbations arranged in dependent blocks produce this separation and contraction of total variation under Markov kernels then yields an approximation obstruction for every procedure. The exact exponent~$(m-2)/(3m-2)$ emerges from the interaction between moment restrictions, total-variation indistinguishability of the constructed distributions and the fluctuation scale of the corresponding Gaussian maxima.

\section{The setting}
Fix~$m>2$ and $0<b_1<b_2<\infty$, and let~$\mc P_d(b_1,b_2,m)$ be the class of Borel probability measures~$P$ on~$\R^d$ satisfying
\begin{align*}
 &\min_{1\leq j\leq d}\int (x_j-\mu_j(P))^2P(dx)\geq b_1^2
 \qquad\text{and}\qquad
 \max_{1\leq j\leq d}\int |x_j-\mu_j(P)|^mP(dx)\leq b_2^m,\\
 &\text{where}\qquad\mu(P):=\int xP(dx)\in\R^d.  
\end{align*}
For~$P\in \mc P_d(b_1,b_2,m)$ we let~$\Sigma(P)$ and~$\Sigma_0(P)$ denote its covariance and correlation matrix, respectively. We stress that~$\mc P_d(b_1,b_2,m)$ does not restrict these matrices beyond what is implied by the previous display. The sample space is $(\R^d)^n$ with its product Borel sigma-field and distributions~$P^{\otimes n}$ with~$P\in \mc P_d(b_1,b_2,m)$. Let~$X_i$ denote the $i$th coordinate projection to~$\R^d$ so that $X_1,\ldots,X_n$ are i.i.d.~with common distribution $P$ under $P^{\otimes n}$. Furthermore, let
\begin{equation}\label{eq:dist_class0}
	\mc P_d^0(b_1,b_2,m):=\cbr[1]{P\in \mc P_d(b_1,b_2,m): \mu(P)=0_d}
\end{equation}
denote the subset of~$\mc P_d(b_1,b_2,m)$ with mean zero; here~$0_d=(0,\hdots,0)'$.

Let~$\mathfrak T_{n,d}$ denote the collection of Borel measurable maps~$T_n:(\R^d)^n\to\R^d$. A map in this class may depend on~$n,d,b_1,b_2$, and~$m$, but not on the unknown~$P$. The discussed suitably scaled sums, self-normalized sums, and trimmed/winsorized means all belong to~$\mathfrak T_{n,d}$. Furthermore, for all of these examples the relevant covariance target for the uniform Gaussian approximation error (recall that~$\gamma_{A}=\mathsf{N}_d(0_d,A)$ for~$A$ a~$d\times d$ real covariance matrix)
\begin{equation} \label{eq:error}
 \Delta_{n,d}(T_n,\Gamma)
 :=
 \sup_{P\in\mc P_d^0(b_1,b_2,m)}\sup_{H\in\mc H_d}
 \big|P^{\otimes n}(T_n\in H)
 -\mathsf \gamma_{\Gamma(P)}(H)\big|,
\end{equation}
is either~$\Gamma=\Sigma$ (the covariance matrix of the~$X_i$) or~$\Gamma=\Sigma_0$ (the correlation matrix of the~$X_i$). We shall henceforth focus on necessary and sufficient conditions on the growth rate of~$d$ for~$\Delta_{n,d}(T_n,\Gamma)\to 0$ when~$\Gamma\in\cbr[0]{\Sigma,\Sigma_0}$.\footnote{Some restriction on the covariance functional~$\Gamma$ is essential. If it can be chosen freely, Proposition~\ref{prop:arbitrarytarget} in Appendix~\ref{sec:structure_needed} provides a trivial example of a~$T_n^*\in\mathfrak{T}_{n,d}$ and a covariance functional~$\Gamma_*$ for which~$\Delta_{n,d}(T_n^*,\Gamma_*)\to 0$ \emph{irrespective} of the growth rate of~$d=d_n$.}

\section{Necessity of~$\log(d)/n^{(m-2)/(3m-2)}\to 0$}
The dimension barrier in Theorem~\ref{thm:sharp} below remains valid also for randomized statistics. A randomized statistic can be represented by a Markov kernel
$\mathsf K_n$ from~$(\R^d)^n$ to~$\R^d$. Under the sample
distribution~$P^{\otimes n}$, its induced distribution on
$(\R^d,\mc B(\R^d))$ is the probability measure
$\mathsf K_nP^{\otimes n}$ defined by
	\begin{equation*}
		(\mathsf{K}_n P^{\otimes n})(A):=\int \mathsf{K}_n(A,x)P^{\otimes n}(dx),\qquad A\in \mc{B}(\R^d)  
	\end{equation*}  
	for~$P\in\mc{P}_d(b_1,b_2,m)$. Denote the class of such Markov kernels by~$\mathfrak{K}_{n,d}$, the elements of which can depend on~$n,d,b_1,b_2$, and~$m$, but not on the unknown~$P$. As always, any non-randomized statistic~$T_n\in\mathfrak{T}_{n,d}$ can be represented via the Markov kernel
	\begin{equation*}
		\mathsf{K}_{T_n}(A,x)
		=
		\delta_{\cbr[0]{T_n(x)}}(A)
		=
		\mathds{1}_{A}(T_n(x))
		\qquad \text{so that}\qquad (\mathsf{K}_{T_n}P^{\otimes n})(A)
		=
		P^{\otimes n}(T_n\in A).
	\end{equation*}
	In this sense~$\mathfrak{T}_{n,d}\subseteq \mathfrak{K}_{n,d}$.\footnote{Statistics~$T_n(X_1,\ldots,X_n,U)$ that include randomization via a random variable~$U$ can also be represented via a Markov kernel: if the kernel $\mathsf L_n(\cdot,x)$ is the conditional distribution of $U$ given the sample $x=(x_1,\hdots,x_n)\in(\R^d)^n$, then the corresponding kernel is
\begin{equation*}
 \mathsf K_n(A,x)
 =
 \int \mathds{1}_{A}(T_n(x,u))\mathsf L_n(du,x).	
\end{equation*}
Here $T_n$ is Borel measurable, and both $T_n$ and
$\mathsf L_n$ are common to all~$P$.} We now state the main dimension barrier of this paper, which holds already over the mean-zero family~$\mc{P}_d^0(b_1,b_2,m)$ and the equal-threshold lower orthants~$\mc R_d:=\cbr[0]{(-\infty,t]^d:t\in\R}\subseteq\mc H_d$.

\begin{theorem}\label{thm:sharp}
Let~$d=d_n\geq2$,~$m\in(2,\infty)$ and suppose that
\begin{equation}\label{sharp:growth}
 \limsup_{n\to\infty}\frac{\log(d)}{n^{(m-2)/(3m-2)}}>0.
\end{equation}
Then, for each~$\Gamma\in\{\Sigma,\Sigma_0\}$,
\begin{equation}\label{sharp:conclusion}
 \limsup_{n\to\infty}
 \inf_{\mathsf{K}_n\in\mathfrak K_{n,d}}
 \sup_{P\in\mc P_d^0(b_1,b_2,m)}\sup_{H\in\mc R_d}
 \envert[1]{(\mathsf{K}_nP^{\otimes n})(H)-\gamma_{\Gamma(P)}(H)}>0.
\end{equation}
If even~$\limsup_{n\to\infty}\log(d)/
n^{(m-2)/(3m-2)}=\infty$, then the
left-hand side of~\eqref{sharp:conclusion} is at least~$1/2$.  If~$\log(d)/n^{(m-2)/(3m-2)}$ converges to a limit in~$(0,\infty]$, then the corresponding conclusion holds with a limit inferior in place of the limit superior in~\eqref{sharp:conclusion}.
\end{theorem}
Theorem~\ref{thm:sharp} shows that~$\log(d)=o\del[0]{n^{(m-2)/(3m-2)}}$ is necessary for the distribution of any statistic, including randomized ones, to admit a Gaussian approximation uniformly over the equal-threshold lower orthants~$\mc{R}_d$ as well as~$\mc{P}^0_d(b_1,b_2,m)$ with covariance functional~$\Gamma\in\{\Sigma,\Sigma_0\}$. In particular, if~\eqref{sharp:growth} holds, then 
\begin{equation*}
 \limsup_{n\to\infty}
 \inf_{T_n\in\mathfrak T_{n,d}}
 \sup_{P\in\mc P_d^0(b_1,b_2,m)}\sup_{t\in\R}
 \envert[3]{P^{\otimes n}\del[2]{\max_{1\leq j\leq d}T_{n,j}\leq t}-\gamma_{\Gamma(P)}\del[1]{(-\infty,t]^d}}>0.	
\end{equation*}
Thus, even the distribution
of the coordinate maximum of an arbitrary 
statistic cannot be approximated uniformly over
$\mc P_d^0(b_1,b_2,m)$ by the corresponding Gaussian maximum. Under~\eqref{sharp:growth}, Theorem~\ref{thm:sharp} rules out neither i) valid approximations for particular sequences  of distributions nor ii) uniform approximations over more structured subsets of~$\mc P_d^0(b_1,b_2,m)$ (imposing, e.g., additional restrictions on~$\Sigma$ or~$\Sigma_0$). 
 
Appendix~\ref{sec:proofSharp} proves a refinement, Theorem~\ref{thm:sharpStronger}, in which the offending orthants in~$\mc{R}_d$ can be chosen independently of the procedure and their upper endpoints are controlled. The endpoint bounds are used in Remark~\ref{rem:conjectureselfnormalized}. 

\subsection{Sketch of proof of Theorem~\ref{thm:sharp}.}
Let~$\mathsf{TV}(Q,Q'):=\sup_E|Q(E)-Q'(E)|$ be the total variation distance between~$Q$ and~$Q'$, where the supremum is over all measurable sets in the common domain of $Q$ and $Q'$.

We outline the argument for~$\Gamma=\Sigma$, following the
five steps of the proof in Appendix~\ref{sec:proofSharp}. Along a suitable subsequence, we construct distributions~$P_{b,n}\in\mc P_d^0(b_1,b_2,m)$,~$b\in\{0,\ldots,k\}$ (Step~1) such that
\begin{equation}\label{eq:TVclose}
 \mathsf{TV}\del[3]{
 P_{0,n}^{\otimes n},
 \frac1k\sum_{b=1}^k P_{b,n}^{\otimes n}}
 \to0
 \qquad\text{as }k=k_n\to\infty\qquad \text{(established in Step~3)}.
\end{equation}
At the same time, the covariance matrices~$\Sigma(P_{b,n})$ of the~$P_{b,n}$ are chosen  (Step 2) so as to ensure that the purported approximating Gaussian probabilities~$\gamma_{\Sigma(P_{0,n})}(A_n)$ and~$\gamma_{\Sigma(P_{b,n})}(A_n)$,~$b\in\cbr[0]{1,\hdots,k}$ are separated uniformly in~$b$ when~\eqref{sharp:growth} holds for a suitably chosen sequence of equal-threshold lower orthant sets~$A_n=(-\infty,b_1t_n]^d\in \mc{R}_d$ (Step 4). In particular,
\begin{equation}\label{eq:Gausssep}
 \limsup_{n\to\infty}
 \min_{1\leq b\leq k}
 \left\{
 \gamma_{\Sigma(P_{0,n})}(A_n)
 -
 \gamma_{\Sigma(P_{b,n})}(A_n)
 \right\}>0.
\end{equation}
By the structure of~$A_n$,~\eqref{eq:Gausssep} amounts to separating the probabilities of coordinate maxima of Gaussian vectors~$Z_{0,n}\sim\gamma_{\Sigma(P_{0,n})}$ and~$Z_{b,n}\sim\gamma_{\Sigma(P_{b,n})}$ in~$\R^d$,~$b\in\cbr[0]{1,\hdots,k}$, respectively, falling below~$b_1t_n$. This step uses the specific block structure of the~$\Sigma(P_{b,n})$ as coordinate maxima of the~$Z_{b,n}$,~$b\in\cbr[0]{0,\hdots,k}$ decompose into the sum of two independent maxima. Classical extreme-value theory applied to suitable comparison maxima yields centering sequences~$t_{0,n}$ (for~$b=0$) and~$t_{1,n}$ (for~$b\in\cbr[0]{1,\hdots,k}$) for these decompositions. Choosing~$t_n:=(t_{0,n}+t_{1,n})/2$ halfway between these centers, we obtain%
\begin{equation}\label{eq:Gaussprobssep}
	\gamma_{\Sigma(P_{0,n})}(A_n)
	=
	\Lambda(D_n)+o(1)\qquad \text{and}\qquad \max_{1\leq b\leq k}\gamma_{\Sigma(P_{b,n})}(A_n)\leq 
	\Lambda(-D_n)+o(1),
\end{equation}
where~$\Lambda$ is the cdf of the standard Gumbel distribution and~$D_n=\sqrt{2\log(k)}\cdot(t_{1,n}-t_{0,n})/2$ is the suitably scaled distance of~$t_n$ to each of the centers~$t_{0,n}$ and~$t_{1,n}$. Since~$\Lambda$ is strictly increasing,~$\liminf_{n\to\infty} D_n>0$ then suffices for~\eqref{eq:Gausssep}. This condition explains the specific dimension exponent~$(m-2)/(3m-2)$ because
\begin{equation*}
	D_n
	\asymp
\del[2]{\frac{\log(k)}{n}}^{(m-2)/(2m)}\log(k)
 =
 \left(
 \frac{\log(k)}{n^{(m-2)/(3m-2)}}
 \right)^{(3m-2)/(2m)}
\end{equation*}
and the choice of~$k$ in Step~1 ensures~$\liminf_{n\to\infty} D_n>0$  whenever~\eqref{sharp:growth} holds. In Step 5 we combine, by standard arguments, the total variation closeness of~$P_{0,n}^{\otimes n}$ and~$\frac{1}{k}\sum_{b=1}^{k}P_{b,n}^{\otimes n}$ in~\eqref{eq:TVclose} with the Gaussian separation in~\eqref{eq:Gausssep} to achieve
\begin{equation*}
\limsup_{n\to\infty}\inf_{\mathsf{K}_{n}\in\mathfrak{K}_{n,d}}\max_{0\leq b\leq k}\envert[1]{(\mathsf{K}_nP_{b,n}^{\otimes n})(A_n)-\gamma_{\Sigma(P_{b,n})}(A_n)}	> 0\qquad\text{if}\qquad \limsup_{n\to\infty}\frac{\log(d)}{n^{\frac{m-2}{3m-2}}}>0,
\end{equation*}      	 
which implies~\eqref{sharp:conclusion} for~$\Gamma=\Sigma$. The proof also establishes the result for the case of~$\Gamma=\Sigma_0$.

\section{Sufficiency of~$\log(d)/n^{(m-2)/(3m-2)}\to 0$}
For~$m\in(2,4]$ and~$\Gamma\in\cbr[0]{\Sigma,\Sigma_0}$, the necessary condition in Theorem~\ref{thm:sharp} is also
sufficient and for~$m\in(4,\infty)$ the condition is sufficient up to a power of~$\log(n)$. The convergence is attained by the (studentized) quantile-winsorized means studied in~\cite{kock2025high}, whose construction we now recall.

For real numbers~$x_1,\hdots,x_n$, denote by~$x_1^*\leq \hdots\leq x_n^*$ their non-decreasing rearrangement. Let~$-\infty<\alpha\leq\beta<\infty$ and define
\begin{align*}
	\phi_{\alpha,\beta}(x)
	=
	\begin{cases}
		\alpha\qquad \text{if }x<\alpha\\
		x\qquad \text{if }x\in[\alpha,\beta]\\
		\beta\qquad \text{if }x>\beta.
	\end{cases}
\end{align*} 
Abbreviate~$\mu=\mu(P)$ and let~$S_{n,W}(P)=S_{n,W}\in\R^d$ be the vector of centered and scaled winsorized means defined via
\begin{equation}\label{eq:WinsMain}
	S_{n,W,j}
	=
	n^{-1/2}\sum_{i=1}^n\del[1]{\phi_{\hat\alpha_j,\hat\beta_j}(X_{i,j})-\mu_j}
	=
	\sqrt{n}(\hat\mu_{n,j}-\mu_j),\qquad j=1,\hdots,d;
\end{equation}
here~$\hat\mu_{n,j}=n^{-1}\sum_{i=1}^n\phi_{\hat\alpha_j,\hat\beta_j}(X_{i,j})$ with~$\hat{\alpha}_j= X_{\lceil \eps_n n \rceil,j}^*$  and $\hat{\beta}_j= X_{\lfloor(1-\eps_n )n\rfloor+1,j}^*$ for~$\eps_n\in(0,1/2)$ of the form 
\begin{equation}\label{eq:epsfamMain}
	\eps_n=\lambda\cdot \frac{\log(dn)}{n}\qquad \text{with }\lambda\in (0,\infty)\text{ fixed}.
\end{equation} 

The choice of~$\eps_n$ determines the fraction of observations that are winsorized. For each coordinate~$j\in\cbr[0]{1,\hdots,d}$, the winsorization points~$\hat{\alpha}_j$ and~$\hat{\beta}_j$ are the lower- and upper-$\eps_n$ order statistics of the~$X_{1,j},\hdots,X_{n,j}$. The winsorization rule does not require knowledge of~$m$ or any other unknown population quantities. We refer to~\cite{kock2025high} for further discussion of the construction of~$S_{n,W}$, whose distribution obeys the following Gaussian approximation.

\subsection{The case of covariance target~$\gamma_\Sigma$}
We first consider the case where the approximating Gaussian has covariance matrix~$\Sigma$.
\begin{theorem}\label{thm:SufficientWins}
	Let~$d=d_n\geq 2$ and~$\eps_n\in(0,1/2)$ with~$\eps_n$ as in~\eqref{eq:epsfamMain}. If~$m\in(2,4]$ and
	\begin{equation}\label{eq:SuffWins}
	 \lim_{n\to\infty}\frac{\log(d)}{n^{(m-2)/(3m-2)}}
	 =
	 0,
	\end{equation}
	then
	\begin{equation}\label{eq:SuffconvWins}
		 \lim_{n\to\infty}
 \sup_{P\in\mc P_d(b_1,b_2,m)}\sup_{H\in\mc H_d}
 \envert[1]{P^{\otimes n}(S_{n,W}\in H)-\gamma_{\Sigma(P)}(H)}
 =
 0.
	\end{equation}
If~$m\in(4,\infty)$ and
\begin{equation}\label{eq:m>4cond}
	\frac{\log(d)}{n^{\frac{m-2}{3m-2}}/\log^{\frac{m}{3m-2}}(n)}\to 0,
\end{equation}
then~\eqref{eq:SuffconvWins} holds.
\end{theorem} 
For~$m\in(2,4]$ and~$\Gamma=\Sigma$, Theorems~\ref{thm:sharp} and~\ref{thm:SufficientWins} together show that
\begin{equation*}
\limsup_{n\to\infty}
 \inf_{\mathsf{K}_n\in\mathfrak K_{n,d}}
 \sup_{P\in\mc P_d^0(b_1,b_2,m)}\sup_{H\in\mc H_d}
 \envert[1]{(\mathsf{K}_nP^{\otimes n})(H)-\gamma_{\Sigma(P)}(H)}=0
\end{equation*}		
\emph{if and only if}
\begin{equation}\label{eq:NecSuff}
	\lim_{n\to\infty}\frac{\log(d)}{n^{(m-2)/(3m-2)}}
	 =
	 0,
\end{equation}
and thus characterize the exact dimension barrier for a Gaussian approximation to hold uniformly over~$\mc P_d^0(b_1,b_2,m)$ and~$\mc{H}_d$. Put differently, no statistic --- whether randomized or not --- can have a distribution that obeys a Gaussian approximation for larger~$d$ than the winsorized means~$S_{n,W}$. We also stress that the Gaussian approximation to the distribution of~$S_{n,W}$ in Theorem~\ref{thm:SufficientWins} is valid over~$\mc P_d(b_1,b_2,m)$ and does not impose~$\mu(P)=0_d$. For example, with bounded fourth moments, the optimal dimension growth condition is~$\log(d)=o(n^{1/5})$ absent further restrictions on~$\Sigma$. In contrast, ordinary sums can fail
to admit a Gaussian approximation already for~$d=\lceil n^{1+\xi}\rceil$, for any fixed~$\xi>0$; cf.~Theorem~2.1 in~\cite{kock2024remark}.

In case~$m\in(4,\infty)$, the necessary condition and our sufficient condition on~$\log(d)$ in~\eqref{eq:m>4cond} differ by the factor~$\log^{\frac{m}{3m-2}}(n)$, but the optimal exponent on~$n$ remains~$(m-2)/(3m-2)$ for all~$m\in(2,\infty)$.  We conjecture that the factor~$\log^{\frac{m}{3m-2}}(n)$ is also redundant when~$m\in(4,\infty)$ so that the necessary and sufficient condition is~\eqref{eq:NecSuff} for all~$m\in(2,\infty)$. This conjecture is true if Conjecture 1 in~\cite{CCKK26} is true.\footnote{This follows from inspection of the penultimate display of the proof of Lemma~\ref{lem:HDGaussBD}, which shows that for~$m\in(4,\infty)$ we inherit a~$\log^{1/2}(n)$-factor from Corollary~1 in~\cite{CCKK26}, which can be removed under their Conjecture 1. For~$m\in(2,4]$, we do not rely on that result and hence do not inherit this~$\log^{1/2}(n)$-factor.}

Theorem \ref{thm:SufficientWins} is a special case of Theorem~\ref{thm:HDGaussSharp} in Appendix~\ref{sec:WinsProof}, which also allows for~$X_1,\hdots,X_n$ to be coordinatewise adversarially contaminated and provides a fixed-$n$ bound on the approximation error. Both theorems allow~$\log(d)$ to grow at a strictly faster rate for all~$m\in(2,\infty)$ than Theorem~2.1 in \cite{kock2025high}, which required the suboptimal~$\log(d)=o\del[1]{n^{\frac{m-2}{5m-2}}}$.

 The global proof strategy of Theorem~\ref{thm:HDGaussSharp} is the same as that of Theorem~2.1 in \cite{kock2025high}. The improved growth rate of~$d$ stems from using two separate Gaussian approximations tailored to sums of centered independent bounded random vectors in~$\R^d$ absent further assumptions on~$\Sigma$ in Lemma~\ref{lem:HDGaussBD} instead of the more general-purpose Gaussian approximation for sums of sub-exponential random vectors employed in Lemma B.2 of \cite{kock2025high}. The first of these two alternative approximations, supplemental Lemma A.1 in~\cite{chernozhuokov2022improved}, yields our growth condition on~$d$ for~$m\in(2,4]$ whereas the other one, Corollary 1 in~\cite{CCKK26}, yields our condition for~$m\in(4,\infty)$. 
 
\subsection{The case of correlation target~$\gamma_{\Sigma_0}$}
We now show that also for Gaussian approximations with the correlation target~$\Sigma_0$ and~$m\in(2,4]$, the necessary condition in~\eqref{eq:NecSuff} is sufficient.~\eqref{eq:m>4cond} remains sufficient for~$m\in(4,\infty)$. This is witnessed by the suitably studentized (self-normalized) winsorized means in~\cite{kock2025high}, which we briefly recall. Set
\begin{equation}\label{eq:epsprimeMain}
	\eps_n'= \lambda'\cdot\frac{\log(d^2n)}{n}\qquad \text{with }\lambda'\in (0,\infty)\text{ fixed}.
\end{equation}
Writing~$\hat{a}_j=X_{\lceil \eps_n'n\rceil,j}^*$,~$\hat{b}_j=X_{\lfloor (1-\eps_n')n\rfloor+1,j}^*$, and~$\tilde{\mu}_{n,j}=n^{-1}\sum_{i=1}^n\phi_{\hat{a}_j,\hat{b}_j}(X_{i,j})$ for~$j=1,\hdots,d$, define~$\tilde{\Sigma}_n$ as the~$d\times d$ matrix with entries 
\begin{equation}\label{eq:tildeSigmaMain}
	\tilde{\Sigma}_{n,j,k}
	=
	\frac{1}{n}\sum_{i=1}^n\sbr[1]{\phi_{\hat{a}_j,\hat{b}_j}(X_{i,j})-\tilde{\mu}_{n,j}}\sbr[1]{\phi_{\hat{a}_k,\hat{b}_k}(X_{i,k})-\tilde{\mu}_{n,k}},\qquad 1\leq j,k\leq d,
\end{equation}
which is positive semi-definite by virtue of being a Gram matrix. Let~$\tilde{\sigma}_{n,j}=(\tilde{\Sigma}_{n,j,j})^{1/2}$ be the square root of the~$j$th diagonal element of~$\tilde{\Sigma}_n$ and let~$S_{n,W,S}=S_{n,W,S}(P)\in\R^d$ be the studentized winsorized means defined via
\begin{equation*}
S_{n,W,S,j}
	=
	\frac{1}{\sqrt{n}\tilde{\sigma}_{n,j}}\sum_{i=1}^n\del[1]{\phi_{\hat\alpha_j,\hat\beta_j}(X_{i,j})-\mu_j}
	=
	\frac{\sqrt{n}}{\tilde{\sigma}_{n,j}}(\hat\mu_{n,j}-\mu_j),\qquad j=1,\hdots,d	
\end{equation*}
with~$S_{n,W,S,j}=0$ in case~$\tilde{\sigma}_{n,j}=0$. The following theorem is the correlation-target analogue of Theorem~\ref{thm:SufficientWins}.
\begin{theorem}\label{thm:CorrTargetSharp}
Let~$d=d_n\geq 2$ and~$\eps_n,\eps_n'\in(0,1/2)$ with~$\eps_n$ as in~\eqref{eq:epsfamMain} and~$\eps_n'$ as in~\eqref{eq:epsprimeMain}. If~$m\in(2,4]$ and~\eqref{eq:SuffWins} holds, then
	\begin{equation}\label{eq:SuffconvWinsCorr}
		 \lim_{n\to\infty}
 \sup_{P\in\mc P_d(b_1,b_2,m)}\sup_{H\in\mc H_d}
 \envert[1]{P^{\otimes n}(S_{n,W,S}\in H)-\gamma_{\Sigma_0(P)}(H)}
 =
 0.
	\end{equation}
If~$m\in(4,\infty)$ and~\eqref{eq:m>4cond} holds, then~\eqref{eq:SuffconvWinsCorr} holds.
\end{theorem}
For~$m\in(2,4]$, Theorems~\ref{thm:sharp}
and~\ref{thm:CorrTargetSharp} establish that the dimension growth
condition~\eqref{eq:SuffWins} is optimal also for the correlation
target~$\gamma_{\Sigma_0}$. No further structure beyond that implied by~$P\in\mc P_d(b_1,b_2,m)$ is imposed
on~$\Sigma_0$, which may in particular be singular. For~$m\in(4,\infty)$,
the sufficient condition~\eqref{eq:m>4cond} has the same logarithmic loss as in Theorem~\ref{thm:SufficientWins}.

Note, in particular, that the studentized winsorized means~$S_{n,W,S}$ admit a Gaussian approximation by~$\gamma_{\Sigma_0}$ uniformly over~$\mc P_d(b_1,b_2,m)$ and~$\mc{H}_d$ in much higher (exponential) dimension than ``plain'' self-normalized sums, for which Gaussian approximation can fail already at polynomial dimension~$d=\lceil n^{m/2-1}\rceil$ absent further restrictions on~$\Sigma_0$, cf.~\cite[Theorem 2.1]{su2026sharp} for~$m\geq 3$.

Theorem~\ref{thm:CorrTargetSharp} is a special case of Theorem~\ref{thm:HDGauss_studentized}, which
also provides a finite-sample bound under coordinatewise adversarial contamination. As was the case for approximations by~$\gamma_{\Sigma}$ in Theorem~\ref{thm:SufficientWins}, the rate condition in Theorems~\ref{thm:CorrTargetSharp} and~\ref{thm:HDGauss_studentized} allows for a strictly faster growth rate of~$\log(d)$ than the suboptimal condition~$\log(d)=o\del[1]{n^{\frac{m-2}{5m-2}}}$ imposed in Theorem 4.1 of~\cite{kock2025high}. 

In light of the polynomial growth barrier that self-normalized sums are subject to, it is perhaps not surprising that we do not directly exploit the self-normalized nature of~$S_{n,W,S}$ in establishing Theorem~\ref{thm:HDGauss_studentized} and hence Theorem~\ref{thm:CorrTargetSharp}. Writing~$\sigma_{2,j}^2$ for the~$j$th diagonal element of~$\Sigma$, we instead use Theorem~\ref{thm:HDGaussSharp} to deduce that the distribution of the population-standardized winsorized means~$V_n:=(S_{n,W,1}/\sigma_{2,1},\hdots,S_{n,W,d}/\sigma_{2,d})$ can be approximated by~$\gamma_{\Sigma_0}$. The precision guarantee for~$\tilde{\Sigma}_{n}$  proved in Theorem~\ref{thm:covestimGRAM} then allows us to deduce that~$S_{n,W,S}=(S_{n,W,1}/\tilde{\sigma}_{n,1},\hdots,S_{n,W,d}/\tilde{\sigma}_{n,d})$ is sufficiently close to~$V_n$ to inherit the approximation by~$\gamma_{\Sigma_0}$ via Gaussian anti-concentration.

\begin{remark}\label{rem:conjectureselfnormalized}
In their Remark~2.4, equation~(2.16), \cite{qiu2025self} conjecture a relative Gaussian upper-tail approximation for the maximum of self-normalized sums, under their moment condition~(2.1), for arbitrary correlation matrices and thresholds~$0\leq x_n=o(n^{\tau/(4+2\tau)})$. Counterexamples to this extension to unrestricted~$\Sigma_0$ for self-normalized sums are given in~\cite{wu2026counterexample}. In Section~\ref{sec:remarkdetails} we explain how the failure of the conjecture to hold uniformly over~$\mc P_d^0(b_1,b_2,3)$ also follows from the general impossibility result in Theorem~\ref{thm:sharpStronger}.
\end{remark}

\section{Conclusion}
This paper established that for \emph{any} sequence of statistics to admit a Gaussian approximation by~$\gamma_{\Sigma}$ or~$\gamma_{\Sigma_0}$ uniformly over~$\mc{P}_d^0(b_1,b_2,m)$ and~$\mc{H}_d$ for any fixed~$m\in(2,\infty)$, it is necessary that~$\log(d)/n^{(m-2)/(3m-2)}\to 0$. Conversely, this condition is sufficient for~$m\in(2,4]$ as witnessed by winsorized means for the target~$\gamma_{\Sigma}$ and studentized winsorized means for the target~$\gamma_{\Sigma_0}$. For~$m\in(4,\infty)$ only a power of~$\log(n)$ separates the necessary and sufficient conditions.

\begin{appendix}
\numberwithin{equation}{section}

\section{Proof of Theorem~\ref{thm:sharp} and a refinement}\label{sec:proofSharp}
Theorem~\ref{thm:sharp} is an immediate consequence of Theorem~\ref{thm:sharpStronger} below, which for each~$\Gamma\in\cbr[0]{\Sigma,\Sigma_0}$ gives an upper bound on the endpoints of some sets in~$\mc{R}_d$ for which Gaussian approximations fail. Furthermore, Theorem~\ref{thm:sharpStronger} shows that the offending sets in~$\mc{R}_d$ can be chosen independently of the kernels~$\mathsf{K}_n$.
\begin{theorem}\label{thm:sharpStronger}
Let~$d=d_n\geq2$,~$m\in(2,\infty)$ and suppose that~\eqref{sharp:growth} holds. For each~$\Gamma\in\{\Sigma,\Sigma_0\}$ there exists a real sequence~$t_n=t_n(\Gamma)$ satisfying~$0\leq t_n\leq 2\sqrt{\log(d)}$ such that for
\begin{equation*}
	A_n^\Gamma:=
 \begin{cases}
  (-\infty,t_n]^d \quad&\text{if}\quad\Gamma=\Sigma_0,\\
  (-\infty,b_1t_n]^d \quad &\text{if}\quad\Gamma=\Sigma,
 \end{cases}
\end{equation*}
it holds that
\begin{equation}\label{sharp:conclusionStrong}
 \limsup_{n\to\infty}
 \inf_{\mathsf{K}_n\in\mathfrak K_{n,d}}
 \sup_{P\in\mc P_d^0(b_1,b_2,m)}
 \envert[1]{(\mathsf{K}_nP^{\otimes n})(A_n^\Gamma)-\gamma_{\Gamma(P)}(A_n^\Gamma)}>0.
\end{equation}
If even~$\limsup_{n\to\infty}\log(d)/n^{(m-2)/(3m-2)}=\infty$, then the
left-hand side of~\eqref{sharp:conclusionStrong} is at least~$1/2$. If~$\log(d)/n^{(m-2)/(3m-2)}$ converges to a limit in~$(0,\infty]$, then the corresponding conclusion holds with a limit inferior in place of the limit superior in~\eqref{sharp:conclusionStrong}.
\end{theorem}

\begin{proof}[Proof of Theorem~\ref{thm:sharp}]
	Because~$A_n^\Gamma\in\mc{R}_d$ in Theorem~\ref{thm:sharpStronger}, all statements of~Theorem~\ref{thm:sharp} follow from Theorem~\ref{thm:sharpStronger}.
\end{proof}

\begin{proof}[Proof of Theorem~\ref{thm:sharpStronger}]
The proof consists of five steps:
\medskip

\noindent\textbf{Step 1: Construction of distributions in~$\mc{P}_d^0(b_1,b_2,m)$.}

Abbreviate the exponent~$\alpha:=(m-2)/(3m-2)\in(0,1/3)$.
By~\eqref{sharp:growth} we can and do select a subsequence along which
\begin{equation*}
 \frac{\log(d_n)}{n^\alpha}\to\kappa
 \qquad\text{for some }\kappa\in(0,\infty].
\end{equation*}
If the limit superior in~\eqref{sharp:growth} is infinite, select
this subsequence with~$\kappa=\infty$. Unless indicated otherwise,
limits below are taken along the selected subsequence. Furthermore, set
\begin{equation*}
 k_n:=\left\lfloor
 \exp\left(\min\left\{\frac{\log(d_n)}{3},\sqrt n\right\}\right)
 \right\rfloor
 \qquad\text{and}\qquad q_n:=\log(k_n).
\end{equation*}
We choose~$k_n$ to ensure that~$k_n^3\leq d_n$, while~$q_n/n\to0 $ and~$q_n/n^\alpha$
has limit as in~\eqref{eq:kappafiniteconv} or~\eqref{eq:kappaINfiniteconv} depending on whether~$\kappa\in(0,\infty)$ or~$\kappa=\infty$. For~$\mathfrak e_n:=\exp(\min\{\log(d_n)/3,\sqrt n\})$, we have~$\mathfrak e_n\to\infty$ and, for all sufficiently large~$n$,
\begin{equation*}
 \log(\mathfrak e_n-1)\leq q_n\leq\log(\mathfrak e_n),
 \qquad\text{where}\qquad
 0\leq\log(\mathfrak e_n)-\log(\mathfrak e_n-1)
 \leq\frac1{\mathfrak e_n-1}\to0.
\end{equation*}
Thus,
\begin{equation*}
 q_n=\log(\mathfrak e_n)+o(1)
 =\min\left\{\frac{\log(d_n)}{3},\sqrt n\right\}+o(1).
\end{equation*}
If~$\kappa<\infty$, then~$\log(d_n)=O(n^\alpha)=o(\sqrt n)$
along the chosen subsequence, and hence
\begin{equation}\label{eq:kappafiniteconv}
 \frac{q_n}{n^\alpha}\to\frac{\kappa}{3}.
\end{equation}
If instead~$\kappa=\infty$, both~$\log(d_n)/(3n^\alpha)$
and~$n^{1/2-\alpha}$ tend to infinity, so that
\begin{equation}\label{eq:kappaINfiniteconv}
 \frac{q_n}{n^\alpha}\to\infty.
\end{equation}
In the following we let~$n$ be sufficiently large to ensure that~$k_n\geq2$,
and write~$k=k_n$ and~$d=d_n$ to simplify notation.

We first construct suitable distributions in~$\mc P_d^0(b_1,b_2,m)$.
To this end, let~$a:=(b_2-b_1)/2>0$,~$r:=(m-2)/m\in(0,1)$, and define
\begin{equation}\label{sharp:constants}
 p_n:=\frac{q_n}{4n},\qquad
 h_n:=ap_n^{-1/m}\qquad\text{and}\qquad
 s_n:=h_n\sqrt{p_n}
     =a\left(\frac{q_n}{4n}\right)^{r/2}.
\end{equation}
For all sufficiently large~$n$ one has~$0<p_n<1$. In addition,~$p_n\to0$ and~$s_n\to0$ because~$q_n/n\to0$. 

On an auxiliary probability space with probability measure~$\P$
and corresponding expectation~$\E$, let~$V=(V_1,\hdots,V_k)'$
have independent Rademacher coordinates. Independently of~$V$,
let~$R_1,\hdots,R_k$ be independent vectors in~$\{-1,1\}^k$,
each having independent Rademacher coordinates, and
write~$R=(R_1',\hdots,R_k')'\in\R^{k^2}$.
For each~$b\in\cbr[0]{0,\hdots,k}$, take independent Bernoulli
variables~$B_1^{(b)},\hdots,B_k^{(b)}$, independent of~$(V,R)$,
with
\begin{equation*}
 B_\ell^{(b)}\sim
 \begin{cases}
  \operatorname{Bernoulli}(p_n)\qquad &\text{if }\ell=b,\\
  \operatorname{Bernoulli}(p_n/3)\qquad &\text{if }\ell\neq b,
 \end{cases}
 \qquad \text{for }\ell\in[k]:=\{1,\hdots,k\}.
\end{equation*}
In particular, all~$B_\ell^{(0)}$ have success
probability~$p_n/3$. 
Write~$B^{(b)}=(B_1^{(b)},\hdots,B_k^{(b)})'$, and set
\begin{equation}\label{sharp:signblocks}
 Y_\ell^{(b)}:=h_nB_\ell^{(b)}R_\ell\in\R^k,
 \qquad \ell\in[k],\quad b\in\cbr[0]{0,\hdots,k}.
\end{equation}
Thus,~$ Y_\ell^{(b)}=0_k$ if~$B_\ell^{(b)}=0$ and~$ Y_\ell^{(b)}=h_nR_\ell$ if~$B_\ell^{(b)}=1$. Let~$Y^{(b)}=((Y_1^{(b)})',\hdots,(Y_k^{(b)})')'\in\R^{k^2}$ and define
\begin{equation}\label{sharp:observations}
 \widetilde X^{(b)}
 =
 \del[2]{(\widetilde X^{(b)}_1)',\hdots,(\widetilde X^{(b)}_k)'}'
 :=
 \begin{pmatrix}
  b_1V_1\mathbf1_{k^2}+Y^{(b)}\\
  \vdots\\
  b_1V_k\mathbf1_{k^2}+Y^{(b)}
 \end{pmatrix}
 \in\R^{k^3},
 \qquad b\in\cbr[0]{0,\hdots,k}.
\end{equation}
We use an additional subscript to denote coordinates within
a block; thus,~$R_{\ell,j}$ and~$Y_{\ell,j}^{(b)}$ denote
the~$j$th coordinates of~$R_\ell$ and~$Y_\ell^{(b)}$,
respectively. For every~$b\in\cbr[0]{0,\hdots,k}$ and~$\ell,j\in[k]$, independence of~$B_\ell^{(b)}$ and~$R_{\ell,j}$, the identity~$(B_\ell^{(b)})^2=B_\ell^{(b)}$, and~$|R_{\ell,j}|=1$ give
\begin{equation*}
 \E Y_{\ell,j}^{(b)}
 =
 h_n\E B_\ell^{(b)}\E R_{\ell,j}=0,\qquad
 \E(Y_{\ell,j}^{(b)})^2
 =
 h_n^2\E B_\ell^{(b)}\qquad\text{and}\qquad
 \E|Y_{\ell,j}^{(b)}|^m=h_n^m\E B_\ell^{(b)}.
\end{equation*}
It follows that~$\E\widetilde  X^{(b)}=0_{k^3}$. Furthermore, since~$(B^{(b)},R)$ is independent of~$V$, it follows that~$Y^{(b)}$ is independent of~$V$ for~$b\in\cbr[0]{0,\hdots,k}$. Therefore, for $u,\ell,j\in[k]$ the variance of
\begin{equation*}
\widetilde X_{u,\ell,j}^{(b)}
=
b_1V_u+Y_{\ell,j}^{(b)}
\end{equation*}
is 
\begin{equation*}
	\E \del[2]{\widetilde X^{(b)}_{u,\ell,j}}^2
	=
	\E\del[2]{b_1V_{u}+Y^{(b)}_{\ell,j}}^2
	=
	\begin{cases}
		b_1^2+h_n^2p_n/3=b_1^2+s_n^2/3\qquad &\text{if } \ell\neq b\\
		b_1^2+h_n^2p_n=b_1^2+s_n^2\qquad &\text{if } \ell= b. 
	\end{cases}
\end{equation*}
Thus the variance of~$\widetilde X^{(b)}_{u,\ell,j}$ is greater than~$b_1^2$ for~$b\in\cbr[0]{0,\hdots,k}$. For later use, define
\begin{equation}\label{sharp:variance}
 \sigma_{0,n}^2:=b_1^2+\frac{s_n^2}{3}
 \qquad \text{and}\qquad  \sigma_{1,n}^2:=b_1^2+s_n^2.
\end{equation}
Note also that~$b_1<\sigma_{0,n}\leq \sigma_{1,n}\leq b_1+s_n< b_1+a< b_2$.  Furthermore, Minkowski's inequality gives (irrespectively of whether~$\ell=b$ or not)
\begin{equation*}
 \|\widetilde  X^{(b)}_{u,\ell,j}\|_m
 \leq b_1+h_np_n^{1/m}
 =b_1+a<b_2.
\end{equation*}
Denote by~$\iota_n:\R^{k^3}\to \R^d$ the embedding into~$\R^d$ that is the identity map for the first~$k^3$ coordinates and a copy of the first coordinate for the remaining~$(d-k^3)$ ones in case~$k^3<d$. Writing~$P_{b,n}$ for the distribution of~$\iota_n(\widetilde X^{(b)})\in\R^d$, it follows by the preceding moment calculations that~$P_{b,n}\in\mc P_d^0(b_1,b_2,m)$ for every~$b\in\cbr[0]{0,\hdots,k}$.

\medskip
\noindent\textbf{Step 2: Covariance matrices of the~$P_{b,n}$.}
We now compute the covariance matrices of the~$\widetilde{X}^{(b)}$, which determine the Gaussian distributions compared
in Step 4. For~$\ell\in[k]$,~$b\in\cbr[0]{0,\hdots,k}$ and distinct coordinates~$j,j'\in[k]$,
independence of~$R_{\ell,j}$ and~$R_{\ell,j'}$ gives
\begin{equation*}
 \E(Y_{\ell,j}^{(b)}Y_{\ell,j'}^{(b)})
 =h_n^2\E B_\ell^{(b)}\,\E(R_{\ell,j}R_{\ell,j'})=0.
\end{equation*}
Furthermore, for every~$b\in\cbr[0]{0,\hdots,k}$ it holds that~$Y_{\ell}^{(b)}$ and~$Y_{\ell'}^{(b)}$ are mean zero and independent for distinct~$\ell$ and~$\ell'$ (each taking values in~$\R^k$). Consequently, writing~$D_{b,n}:=\E[Y^{(b)}(Y^{(b)})']$, we have (using also that~$\E(Y_{\ell,j}^{(b)})^2=s_n^2/3$ for~$\ell\neq b$ and~$\E(Y_{\ell,j}^{(b)})^2=s_n^2$ for~$\ell=b$)
\begin{equation*}
 D_{0,n}=\frac{s_n^2}{3}I_{k^2}\qquad\text{and}\qquad D_{b,n}=\text{diag}\del[2]{\frac{s_n^2}{3},\hdots,s_n^2,\hdots,\frac{s_n^2}{3}}\otimes I_k\quad \text{for }b\in[k],
\end{equation*}
the~$s_n^2$ being the~$b$th diagonal term of the~$k\times k$ diagonal matrix. By~\eqref{sharp:observations} it then follows that for~$u,u'\in[k]$
\begin{equation}
\E\del[1]{\widetilde X_{u}^{(b)}(\widetilde X_{u'}^{(b)})'}=
b_1^2\mathds{1}(u=u')\cdot\mathbf1_{k^2}\mathbf1_{k^2}'+D_{b,n}.
\end{equation}
To summarize, it follows that
\begin{align}
 \operatorname{Cov}\del[1]{\widetilde X^{(b)}}
 &=
 \begin{pmatrix}
  b_1^2\mathbf1_{k^2}\mathbf1_{k^2}'+D_{b,n}
   &D_{b,n}&\cdots&D_{b,n}\\
  D_{b,n}
   &b_1^2\mathbf1_{k^2}\mathbf1_{k^2}'+D_{b,n}
   &\cdots&D_{b,n}\\
  \vdots&\vdots&\ddots&\vdots\\
  D_{b,n}&D_{b,n}&\cdots
   &b_1^2\mathbf1_{k^2}\mathbf1_{k^2}'+D_{b,n}
 \end{pmatrix}\notag\\
 &=
I_k\otimes (b_1^2\mathbf1_{k^2}\mathbf1_{k^2}') +\mathbf1_{k}\mathbf1_{k}'\otimes D_{b,n}\label{sharp:covariance},
\end{align}
so that the covariance matrices only differ via~$\mathbf1_{k}\mathbf1_{k}'\otimes D_{b,n}$ across~$b\in\cbr[0]{0,\hdots,k}$. 

\medskip

\noindent\textbf{Step 3: Total variation bound between~$P_{0,n}^{\otimes n}$ and~$\frac1k\sum_{b=1}^kP_{b,n}^{\otimes n}$.}
We next show that the total variation distance between the probability measures 
\begin{equation*}
 Q_{0,n}:=P_{0,n}^{\otimes n},
 \qquad\text{and}\qquad 
 \overline Q_n:=\frac1k\sum_{b=1}^kP_{b,n}^{\otimes n}
\end{equation*}
converges to zero. For each~$b\in\cbr[0]{0,\hdots,k}$, let~$\del[0]{V^{(i)},R^{(i)},B^{(b,i)}},
 i\in[n]:=\cbr[0]{1,\hdots,n}$ be independent copies of~$(V,R,B^{(b)})$ defined prior to~\eqref{sharp:signblocks}, and denote their
joint distribution by~$\widetilde Q_{b,n}$. Applying the construction in~\eqref{sharp:observations}
and then the embedding~$\iota_n$ to each independent copy gives a
common (across~$b$) Borel measurable map
$F_n:\R^{n(2k+k^2)}\to(\R^d)^n$ such that~$ P_{b,n}^{\otimes n}
 =\widetilde Q_{b,n}\circ F_n^{-1}$, $b\in\cbr[0]{0,\hdots,k}.$
Thus, also
\begin{equation*}
\overline Q_n
=
\frac1k\sum_{b=1}^k\del[1]{\widetilde Q_{b,n}\circ F_n^{-1}}
=	
\del[3]{\frac1k\sum_{b=1}^k\widetilde Q_{b,n}}\circ F_n^{-1}
\end{equation*} 
and so it suffices to bound the total variation distance between~$\widetilde Q_{0,n}$ and~$\frac1k\sum_{b=1}^k\widetilde Q_{b,n}$. To this end, we study the likelihood ratio of these. For~$b\in\cbr[0]{0,\hdots,k}$ and~$\ell\in[k]$, write (recall that~$0<p_n<1$ for~$n$ sufficiently large)
\begin{equation*}
 \pi_{\ell,n}^{(b)}:=
 \begin{cases}
  p_n,&\ell=b,\\
  p_n/3,&\ell\neq b.
 \end{cases}
\end{equation*}
The joint distribution~$\widetilde Q_{b,n}$ of~$(V^{(i)},R^{(i)},B^{(b,i)})_{i=1}^n$,  is supported on 
\begin{equation*}
 \mathcal W_n:=
 \del[1]{\{-1,1\}^k
          \times\{-1,1\}^{k^2}
          \times\{0,1\}^k}^{n}\qquad \text{for all}\qquad b\in\cbr[0]{0,\hdots,k}.
\end{equation*}
Write a generic element of~$ \mathcal W_n$ as~$w=((v_i,r_i,z_i))_{i=1}^n$,
where~$z_i=(z_{i,1},\hdots,z_{i,k})'$.
By the definition of~$\widetilde Q_{b,n}$ and the independence
in its construction, its probability mass function is
\begin{equation*}
 f_{b,n}(w)
 :=\widetilde Q_{b,n}(\{w\})
 =2^{-n(k+k^2)}
   \prod_{i=1}^n\prod_{\ell=1}^k
   \del[1]{\pi_{\ell,n}^{(b)}}^{z_{i,\ell}}
   \del[1]{1-\pi_{\ell,n}^{(b)}}^{1-z_{i,\ell}},
 \qquad w\in\mathcal W_n.
\end{equation*}
Note that~$f_{b,n}(w)>0$ for every~$w\in\mathcal W_n$ and every~$b\in\cbr[0]{0,\hdots,k}$ implying, in particular, that~$\widetilde Q_{b,n}$ is absolutely continuous with respect to~$\widetilde Q_{0,n}$. Since~$\pi_{\ell,n}^{(b)}=\pi_{\ell,n}^{(0)}=p_n/3$
for~$\ell\neq b$, all factors indexed by~$\ell\neq b$ cancel in~$f_{b,n}(w)/f_{0,n}(w)$. The common factor~$2^{-n(k+k^2)}$ cancels as well.
Therefore, writing~$ N_{\ell,n}(w):=\sum_{i=1}^nz_{i,\ell}$,~$\ell\in[k],\ w\in\mathcal W_n$, the likelihood ratio of~$\widetilde Q_{b,n}$ with respect to~$\widetilde Q_{0,n}$ is
\begin{equation*}
 L_{b,n}(w)
 :=\frac{f_{b,n}(w)}{f_{0,n}(w)}
 =\prod_{i=1}^n
   \frac{p_n^{z_{i,b}}(1-p_n)^{1-z_{i,b}}}
        {(p_n/3)^{z_{i,b}}(1-p_n/3)^{1-z_{i,b}}}
 =
 3^{N_{b,n}(w)}
   \left(\frac{1-p_n}{1-p_n/3}\right)^{n-N_{b,n}(w)},
\end{equation*}
for~$w\in\mathcal W_n$ and~$b\in[k]$.
The factorization of~$f_{0,n}$ shows that under
$\widetilde Q_{0,n}$ the functions
$N_{1,n},\hdots,N_{k,n}$ are independent
$\operatorname{Binomial}(n,p_n/3)$ variables:
each~$N_{\ell,n}(w)$ sums the entries in column~$\ell$ of
$(z_{i,\ell})_{i,\ell}$, and different columns are independent. Thus, since~$L_{b,n}$ is a function of~$N_{b,n}$ only,
the likelihood ratios~$L_{1,n},\hdots,L_{k,n}$ are also independent
under~$\widetilde Q_{0,n}$.

Writing~$\E_0$ for integration with respect to
$\widetilde Q_{0,n}$, it holds by construction that~$ \E_0L_{b,n}=1$ for~$b\in[k]$. Furthermore, using that~$N_{b,n}\sim \operatorname{Binomial}(n,p_n/3)$ under~$\widetilde Q_{0,n}$ as well as the binomial theorem,
\begin{align*}
 \E_0L_{b,n}^2
 &=
 \sum_{t=0}^n
   \binom nt
   \left(\frac{p_n}{3}\right)^t
   \left(1-\frac{p_n}{3}\right)^{n-t}
   3^{2t}
   \left(\frac{1-p_n}{1-p_n/3}\right)^{2(n-t)}\\
 &=
 \sum_{t=0}^n
   \binom nt
   (3p_n)^t
   \left(\frac{(1-p_n)^2}{1-p_n/3}\right)^{n-t}\\
 &=
 \left(3p_n+\frac{(1-p_n)^2}{1-p_n/3}\right)^n=c_n^n.
\end{align*}
where~$c_n:=1+4p_n/(3-p_n)\leq1+2p_n$ as~$0<p_n<1$.
Recall that the likelihood ratios~$L_{1,n},\hdots,L_{k,n}$ are independent
under~$\widetilde Q_{0,n}$, and that~$\overline L_n:=k^{-1}\sum_{b=1}^kL_{b,n}$ is the likelihood
ratio of~$\frac1k\sum_{b=1}^k\widetilde Q_{b,n}$ with respect to~$\widetilde Q_{0,n}$. Hence,
\begin{equation*}
 \E_0(\overline L_n-1)^2
 =
\E_0\del[4]{\frac1{k}\sum_{b=1}^k(L_{b,n}-1)}^2
 =
 \frac1{k^2}\sum_{b=1}^k\E_0(L_{b,n}-1)^2
 =
 \frac{c_n^n-1}{k}.
\end{equation*}
Therefore, since~$ \mathsf{TV}(Q_{0,n},\overline Q_n)$ is upper bounded by~$ \mathsf{TV}(\widetilde Q_{0,n},k^{-1}\sum_{b=1}^k\widetilde Q_{b,n})$, it follows by the Cauchy-Schwarz inequality,~$1+x\leq e^x$, and recalling that~$p_n=\log(k)/(4n)$, that
\begin{equation}\label{sharp:TV}
 \begin{aligned}
 \mathsf{TV}(Q_{0,n},\overline Q_n)
 \leq 
 \E_0|\overline L_n-1|
 \leq 
 \sqrt{\frac{c_n^n-1}{k}}
 \leq
 \sqrt{\frac{(1+2p_n)^n-1}{k}}
 \leq
 \sqrt{\frac{e^{2np_n}}{k}}
 =k^{-1/4}
 \to0.
 \end{aligned}
\end{equation}

\medskip

\noindent\textbf{Step 4: Gaussian separation.}
For each~$\Gamma\in\cbr[0]{\Sigma,\Sigma_0}$ we now exhibit~$A_n\in\mc{R}_d$ such that~$\min_{1\leq b\leq k}[\gamma_{\Gamma(P_{0,n})}(A_n)-\gamma_{\Gamma(P_{b,n})}(A_n)]$ is bounded away from zero along the chosen subsequence. To this end, it is convenient to construct normal distributions with the same covariance matrices as the~$\widetilde X^{(b)}$ in~\eqref{sharp:observations} (cf.~\eqref{sharp:covariance}) via suitable normal random variables mimicking the construction of~$\widetilde X^{(b)}$ there: On a further auxiliary probability space, let~$Z_1,\hdots,Z_k$
be independent standard normals and let~$\zeta_1,\hdots,\zeta_k$
be independent~$\mathsf N_k(0_k,I_k)$ vectors, independent of the~$Z_u$. Set
\begin{equation*}
 \eta_\ell^{(0)}:=\frac{s_n}{\sqrt3}\zeta_\ell\quad\text{for }\ell\in[k]
 \qquad\text{and}\qquad
 \eta_\ell^{(b)}:=
 \begin{cases}
  s_n\zeta_b \quad&\text{if }\ell=b,\\
  (s_n/\sqrt3)\zeta_\ell\quad &\text{if }\ell\neq b,
 \end{cases}
 \qquad\text{for } b\in[k].
\end{equation*}
Now let~$\eta^{(b)}:=\del[1]{(\eta^{(b)}_1)',\hdots,(\eta^{(b)}_k)'}'\in\R^{k^2}$.
Define the centered Gaussian vector
\begin{equation*}
 G_{b,n}:=
 \begin{pmatrix}
  b_1Z_1\mathbf1_{k^2}+\eta^{(b)}\\
  \vdots\\
  b_1Z_k\mathbf1_{k^2}+\eta^{(b)}
 \end{pmatrix}\in\R^{k^3},
\end{equation*}
By~\eqref{sharp:covariance},~$G_{b,n}$ and
$\widetilde X^{(b)}$ have the same covariance matrix. Since~$\iota_n$ (defined in the end of Step 1) is linear, the centered Gaussian vector~$\iota_n(G_{b,n})\in\R^d$ therefore has law~$\gamma_{\Sigma(P_{b,n})}$. Fix~$\Gamma\in\cbr[0]{\Sigma,\Sigma_0}$ and let (recall the definition of~$\sigma_{0,n}$ and~$\sigma_{1,n}$ in~\eqref{sharp:variance})
\begin{equation*}
 (\rho_{0,n},\rho_{1,n})
 :=
 \begin{cases}
  (\sigma_{0,n},\sigma_{1,n})\quad &\text{if}\quad\Gamma=\Sigma_0,\\
  (b_1,b_1)\quad &\text{if}\quad\Gamma=\Sigma,
 \end{cases}
 \qquad
 A_n^\Gamma(t):=
 \begin{cases}
  (-\infty,t]^d \quad&\text{if}\quad\Gamma=\Sigma_0,\\
  (-\infty,b_1t]^d \quad &\text{if}\quad\Gamma=\Sigma.
 \end{cases}
\end{equation*}
Write~$M_j$ for the maximum of~$j\in\N$ independent standard normals and~$M_j'$ for an independent copy of~$M_j$. Define
\begin{equation}\label{sharp:comparisonmaxima}
 B_{0,n}:=\frac{b_1M_k+(s_n/\sqrt3)M_{k^2}}{\rho_{0,n}}
 \qquad\text{and}\qquad
 B_{1,n}:=\frac{b_1M_k+s_nM_k'}{\rho_{1,n}},
\end{equation}
where~$M_k$ is independent of~$M_{k^2}$ (and~$M_k'$). Write~$\eta_{\ell,j}^{(b)}$ for the~$j$th coordinate of~$\eta_\ell^{(b)}$ and note that ($\stackrel{d}{=}$ denoting equality in distribution)
 \begin{align*}
  \max_{1\leq j\leq k^3}\frac{G_{0,n,j}}{\rho_{0,n}}
  &=
  \max_{u,\ell,j\in[k]}
  \frac{b_1Z_u+\eta_{\ell,j}^{(0)}}{\rho_{0,n}}
  =
  \frac{b_1\max_{u\in[k]}Z_u+
        \max_{\ell,j\in[k]}\eta_{\ell,j}^{(0)}}{\rho_{0,n}}
  \ \stackrel{d}{=}\ B_{0,n},\\
  \max_{1\leq j\leq k^3}\del[2]{K_{b,n}^{-1}G_{b,n}}_j
  &\geq 
  \max_{u,j\in[k]}
  \frac{b_1Z_u+\eta_{b,j}^{(b)}}{\rho_{1,n}}
  =
  \frac{b_1\max_{u\in[k]}Z_u+
        \max_{j\in[k]}\eta_{b,j}^{(b)}}{\rho_{1,n}}
  \ \stackrel{d}{=}\ B_{1,n},
  \qquad b\in[k],
 \end{align*}
where~$K_{b,n}=I_k \otimes \text{diag}[(\rho_{0,n},\hdots\rho_{1,n},\hdots,\rho_{0,n})\otimes \mathbf{1}_k']$ with~$\rho_{1,n}$ being in the~$b$th position of $(\rho_{0,n},\hdots\rho_{1,n},\hdots,\rho_{0,n})'$ and $\text{diag}(v)$ returning a diagonal matrix with the elements of~$v$ on the diagonal. The inequality in the previous display follows from restricting the maximum to those entries of~$G_{b,n}$ that have variance~$\sigma_{1,n}^2$. The augmentation in~$\iota_n(G_{b,n})$ in case~$k^3<d$ repeats the first coordinate~$d-k^3$ times and therefore preserves the maximum value. Consequently, for every~$t\in\R$,
\begin{equation}\label{sharp:Gaussianbounds}
 \begin{aligned}
  \gamma_{\Gamma(P_{0,n})}(A_n^\Gamma(t))
  &=\P(B_{0,n}\leq t),\\
  \gamma_{\Gamma(P_{b,n})}(A_n^\Gamma(t))
  &\leq\P(B_{1,n}\leq t),\qquad b\in[k].
 \end{aligned}
\end{equation}
The inequality follows by retaining only the constraints
on the coordinates in the special block~$\ell=b$ of each ``row''. Note that the distribution of~$B_{1,n}$ does not depend on~$b\in[k]$. We now show that~$\P(B_{0,n}\leq t)$ and~$\P(B_{1,n}\leq t)$ are asymptotically separated for a suitable choice of~$t=t_n$. To this end, for~$j\geq2$, define
\begin{equation*}
 v_j:=\sqrt{2\log(j)},\qquad\text{and}\qquad
 w_j:=\sqrt{2\log(j)}-\frac{\log(\log(j))+\log(4\pi)}{2\sqrt{2\log(j)}}.
\end{equation*}
and the corresponding centered and normalized maxima
\begin{equation*}
 W_k:=v_k(M_k-w_k),\qquad
 W_{k^2}:=v_{k^2}(M_{k^2}-w_{k^2})\qquad\text{and}\qquad
 W_k':=v_k(M_k'-w_k).
\end{equation*}
The normal extreme-value theorem, as stated in Theorem~1.5.3 of~\cite{leadbetter}, gives that~$W_k$,~$W_{k^2}$, and~$W_k'$ converge in distribution to~$W$ as~$n\to\infty$, where~$\P(W\leq x)=\Lambda(x):=\exp(-e^{-x})$, $x\in\R$. Thus,~$W_k$,~$W_{k^2}$,
and~$W_k'$ are bounded in probability and
\begin{equation*}
 \frac{v_k}{v_{k^2}}=\frac1{\sqrt2},
 \qquad
 \frac{w_k}{v_k}\to1,
 \qquad\text{and}\qquad
 \frac{w_{k^2}}{v_k}\to\sqrt2.
\end{equation*}
To choose a threshold~$t_n\in\R$, note that for
\begin{equation*}
	 t_{0,n}:=\frac{b_1w_k+(s_n/\sqrt3)w_{k^2}}{\rho_{0,n}}\qquad \text{and}\qquad
 t_{1,n}:=\frac{(b_1+s_n)w_k}{\rho_{1,n}},
\end{equation*}
the identities~$M_j=w_j+W_j/v_j$ along with~$s_n\to 0$,~$\rho_{0,n}\to b_1$,~$\rho_{1,n}\to b_1$, and Slutsky's theorem give ($\stackrel{d}{\to}$ denoting convergence in distribution)
\begin{align*}
  v_k(B_{0,n}-t_{0,n})
  &=\frac{b_1}{\rho_{0,n}}W_k
    +\frac{s_n}{\sqrt6\,\rho_{0,n}}W_{k^2}
    \stackrel{d}{\to}W,\\
  v_k(B_{1,n}-t_{1,n})
  &=\frac{b_1}{\rho_{1,n}}W_k
    +\frac{s_n}{\rho_{1,n}}W_k'
    \stackrel{d}{\to}W,
\end{align*}
the corresponding distribution functions converging uniformly to~$\Lambda$ by P\'olya's theorem. Therefore, for any real sequence $t_n$ 
\begin{align*}
	\P\del[1]{B_{0,n}\leq t_n}
	&=
	\P\del[1]{v_k(B_{0,n}-t_{0,n})\leq v_k(t_n-t_{0,n})}
	=
	\Lambda(v_k(t_n-t_{0,n}))+o(1)\\
		\P\del[1]{B_{1,n}\leq t_n}
	&=
	\P\del[1]{v_k(B_{1,n}-t_{1,n})\leq v_k(t_n-t_{1,n})}
	=
	\Lambda(v_k(t_n-t_{1,n}))+o(1)
\end{align*}
We choose 
\begin{equation}\label{eq:tnDn}
 t_n:=\frac{t_{0,n}+t_{1,n}}2
 \qquad\text{and}\qquad
 D_n:=\frac{v_k}{2}(t_{1,n}-t_{0,n})
\end{equation}
so that~$v_k(t_n-t_{0,n})=D_n$,
$v_k(t_n-t_{1,n})=-D_n$, and define
\begin{equation*}
 F_{0,n}:=\P(B_{0,n}\leq t_n)=\Lambda(D_n)+o(1)
 \qquad\text{and}\qquad
 F_{1,n}:=\P(B_{1,n}\leq t_n)=\Lambda(-D_n)+o(1).
\end{equation*}
Since~$\Lambda$ is strictly increasing, it therefore
suffices to prove that~$\liminf_{n\to\infty}D_n>0$.
To this end, first note that~$\rho_{1,n}^{-1}-\rho_{0,n}^{-1}=0$ when~$\Gamma=\Sigma$,
whereas for~$\Gamma=\Sigma_0$,
\begin{equation*}
 \frac1{\rho_{1,n}}-\frac1{\rho_{0,n}}
 =
 \frac{\sigma_{0,n}-\sigma_{1,n}}{\sigma_{0,n}\sigma_{1,n}}
 =
  \frac{\sigma_{0,n}^2-\sigma_{1,n}^2}{\sigma_{0,n}\sigma_{1,n}(\sigma_{0,n}+\sigma_{1,n})}
 =
 -\frac{2s_n^2/3}
         {\sigma_{0,n}\sigma_{1,n}(\sigma_{0,n}+\sigma_{1,n})}
 =
 O(s_n^2).
\end{equation*}
It follows that, for either value of~$\Gamma\in\cbr[0]{\Sigma,\Sigma_0}$, using also that~$s_n\to 0$,
\begin{equation*}
 \frac{t_{1,n}-t_{0,n}}{s_nv_k}
 =
 b_1\frac{w_k}{v_k}
       \frac{\rho_{1,n}^{-1}-\rho_{0,n}^{-1}}{s_n}
 +\frac{w_k}{v_k}\frac1{\rho_{1,n}}
 -\frac{w_{k^2}}{v_k}\frac1{\sqrt3\rho_{0,n}}
 \to
 \frac1{b_1}\left(1-\sqrt{\frac23}\right)
 :=\frac{c_*}{b_1}.
\end{equation*}
Because~$v_k^2=2\log(k)=2q_n$, we obtain
\begin{equation}\label{sharp:Dasymptotic}
 \frac{D_n}{\frac{c_*}{b_1}s_nq_n}
 =
  \frac{t_{1,n}-t_{0,n}}{s_nv_k}\Big/ (c_*/b_1)
  \to 1\quad\text{and}\quad
\frac{c_*}{b_1}s_nq_n
 =
 \frac{ac_*}{b_1\,4^{r/2}}
       \frac{q_n^{1+r/2}}{n^{r/2}}
 =
 \frac{ac_*}{b_1\,4^{r/2}}
       \left(\frac{q_n}{n^\alpha}\right)^{1+r/2},
\end{equation}
where the last equality used~$\frac{r/2}{1+r/2}=\frac{r}{2+r}=\frac{m-2}{3m-2}=\alpha$ as~$r=(m-2)/m$, cf.~its definition prior to~\eqref{sharp:constants}.

Suppose first that~$\kappa<\infty$. By~\eqref{sharp:Dasymptotic} and~\eqref{eq:kappafiniteconv}
\begin{equation*}
 D_n
 \to D:=
 \frac{ac_*}{b_1\,4^{r/2}}
       \left(\frac{\kappa}{3}\right)^{1+r/2}>0.
\end{equation*}
Thus, continuity and strict monotonicity of~$\Lambda$ imply
\begin{equation}\label{sharp:Ffinite}
 F_{0,n}\to\Lambda(D),\qquad
 F_{1,n}\to\Lambda(-D)\qquad\text{and}\qquad
 F_{0,n}-F_{1,n}\to\Lambda(D)-\Lambda(-D)>0.
\end{equation}

Suppose next that~$\kappa=\infty$.
Then~\eqref{sharp:Dasymptotic} and~\eqref{eq:kappaINfiniteconv} give~$D_n\to\infty$. It follows that
\begin{equation}\label{sharp:Finfinite}
 F_{0,n}\to1\qquad\text{and}\qquad F_{1,n}\to 0.
\end{equation}

\medskip
\noindent\textbf{Step 5: Combining total variation closeness with Gaussian separation.} We now conclude by combining the total variation closeness of~$Q_{0,n}=P_{0,n}^{\otimes n}$ and~$\overline{Q}_n=\frac1k\sum_{b=1}^kP_{b,n}^{\otimes n}$ in~\eqref{sharp:TV} with the just established Gaussian separation.  
For every~$\mathsf{K}_n\in\mathfrak K_{n,d}$, contraction of total variation gives
	(cf., e.g., Theorem 1.70 and (1.94) in~\cite{liese2008statistical}), 
\begin{align*}
 \mathsf{TV}\del[4]{
 \mathsf{K}_nP_{0,n}^{\otimes n},
\frac1k\sum_{b=1}^k \mathsf{K}_nP_{b,n}^{\otimes n}}
 &=
 \mathsf{TV}\del[4]{
 \mathsf{K}_nP_{0,n}^{\otimes n},
 \mathsf{K}_n\del[3]{\frac1k\sum_{b=1}^kP_{b,n}^{\otimes n}}}\\
 &\leq
 \mathsf{TV}\del[3]{
 P_{0,n}^{\otimes n},
 \frac1k\sum_{b=1}^kP_{b,n}^{\otimes n}
 }
=:
 \delta_n\to 0,
\end{align*}
the convergence following from~\eqref{sharp:TV}. Thus, for~$A_n:=A_n^\Gamma(t_n)\in\mc R_d$ with~$t_n$ as in~\eqref{eq:tnDn} and any~$\mathsf{K}_n\in\mathfrak{K}_{n,d}$
\begin{align*}
	\envert[2]{\gamma_{\Gamma(P_{0,n})}(A_n)-\frac{1}{k}\sum_{b=1}^k \gamma_{\Gamma(P_{b,n})}(A_n)}
	&\leq 
	\envert[1]{\gamma_{\Gamma(P_{0,n})}(A_n)- (\mathsf{K}_nP_{0,n}^{\otimes n})(A_n)}\\
	&\ +
	\envert[2]{(\mathsf{K}_nP_{0,n}^{\otimes n})(A_n)-\frac1k\sum_{b=1}^k (\mathsf{K}_nP_{b,n}^{\otimes n})(A_n)}\\
	&\ +\envert[2]{\frac1k\sum_{b=1}^k (\mathsf{K}_nP_{b,n}^{\otimes n})(A_n)-\frac{1}{k}\sum_{b=1}^k\gamma_{\Gamma(P_{b,n})}(A_n)}\\
	&\leq
	2\max_{0\leq b\leq k}\envert[1]{(\mathsf{K}_nP_{b,n}^{\otimes n})(A_n)-\gamma_{\Gamma(P_{b,n})}(A_n)}+\delta_n,
\end{align*}	
so that by~\eqref{sharp:TV}, using also that~$\gamma_{\Gamma(P_{0,n})}(A_n)=F_{0,n}$ and~$\gamma_{\Gamma(P_{b,n})}(A_n)\leq F_{1,n}$ for all~$b\in[k]$, (cf.~\eqref{sharp:Gaussianbounds})
\begin{align*}
\max_{0\leq b\leq k}\envert[1]{(\mathsf{K}_nP_{b,n}^{\otimes n})(A_n)-\gamma_{\Gamma(P_{b,n})}(A_n)}
&\geq 
\frac{\gamma_{\Gamma(P_{0,n})}(A_n)-\frac{1}{k}\sum_{b=1}^k \gamma_{\Gamma(P_{b,n})}(A_n)-\delta_n}{2}\\
&\geq
\frac{F_{0,n}-F_{1,n}-\delta_n}{2}.	
\end{align*}
Thus, since~$P_{b,n}\in\mc{P}_{d}
^0(b_1,b_2,m)$ and the right-hand side does not depend on~$
\mathsf{K}_n\in\mathfrak{K}_{n,d}$, it also holds that
\begin{equation}\label{eq:LBbeforelimit}
	\inf_{\mathsf{K}_n\in\mathfrak K_{n,d}}
 \sup_{P\in\mc P_d^0(b_1,b_2,m)}
 \envert[1]{(\mathsf{K}_nP^{\otimes n})(A_n)-\gamma_{\Gamma(P)}(A_n)}
	\geq
	\frac{F_{0,n}-F_{1,n}-\delta_n}{2}
\end{equation}

If~$\kappa<\infty$, then~\eqref{eq:LBbeforelimit} together with~\eqref{sharp:Ffinite} and~$\delta_n\to 0$ yields, along the chosen subsequence,
\begin{equation}
	\liminf_{n\to\infty}	\inf_{\mathsf{K}_n\in\mathfrak K_{n,d}}
 \sup_{P\in\mc P_d^0(b_1,b_2,m)} \envert[1]{(\mathsf{K}_nP^{\otimes n})(A_n)-\gamma_{\Gamma(P)}(A_n)}
 \geq 
 \frac{\Lambda(D)-\Lambda(-D)}{2}
 >0.
\end{equation}

If~$\kappa=\infty$, then~\eqref{eq:LBbeforelimit} together with~\eqref{sharp:Finfinite} and~$\delta_n\to 0$ yields, along the chosen subsequence,
\begin{equation}
		\liminf_{n\to\infty}	\inf_{\mathsf{K}_n\in\mathfrak K_{n,d}}
 \sup_{P\in\mc P_d^0(b_1,b_2,m)} \envert[1]{(\mathsf{K}_nP^{\otimes n})(A_n)-\gamma_{\Gamma(P)}(A_n)}
 \geq 
 \frac{1}{2}.
\end{equation}
  
Finally, note that along the selected subsequence it holds for each~$\Gamma\in\cbr[0]{\Sigma,\Sigma_0}$ that $t_n/w_k\to 1$ (recall the definition of~$t_n=t_n(\Gamma)$ in~\eqref{eq:tnDn},~$w_k$ prior to this and use that~$s_n\to 0$,~$\rho_{l,n}(\Gamma)\to b_1$ for~$l\in\cbr[0]{0,1}$,~$w_{k^2}/w_k\to \sqrt{2}$). Therefore, because~$k^3\leq d$, it eventually holds that
\begin{equation*}
	0
	<
	t_n
	\leq 2w_k
	\leq 
	2\sqrt{2\log(k)}
	\leq
	2\sqrt{\log(d)} 
	\end{equation*}
For all remaining indices~$n$, including the finitely many
exceptional indices of the selected subsequence, set
~$t_n=2\sqrt{\log(d)}$.
Then~$0\leq t_n\leq 2\sqrt{\log(d)}$ for every~$n$.
  
If~$\log(d_n)/n^\alpha\to\kappa\in(0,\infty]$
on the full sequence, carry out the preceding construction
on the full sequence. For~$\kappa<\infty$, the argument
gives the strictly positive lower bound
~$(\Lambda(D)-\Lambda(-D))/2$ for the limit inferior;
for~$\kappa=\infty$, it gives the lower bound~$1/2$.
This proves the stated limit-inferior assertions.
\end{proof}

\section{Proof of Theorem~\ref{thm:SufficientWins}}\label{sec:WinsProof}

\subsection{A generalization}\label{sec:ContamSetting}
Throughout this appendix let~$m>2$,~$d \geq 2$ and~$n > 3$ (such that, e.g.,~$\sqrt{\log(d)}>0$) and let~$P\in\mc{P}_d(b_1,b_2,m)$. We prove a more general version of Theorem \ref{thm:SufficientWins} that also allows for some adversarial contamination. To be precise, an adversary inspects the (clean) observations~$X_1,\hdots, X_n$ and returns contaminated~$\tilde{X}_1,\hdots,\tilde{X}_n$ random vectors in~$\R^d$ satisfying 

\begin{equation}\label{eq:contamfrac}
	\envert[1]{\cbr[1]{i\in\cbr[0]{1,\hdots,n}:\tilde{X}_{i,j}\neq X_{i,j}}}
	\leq
	\overline{\eta}_n n \qquad \text{for all}\qquad j\in\cbr[0]{1,\hdots,d}
\end{equation}
to the statistician, where~$\overline{\eta}_n\in[0,1/2)$ is a non-random and known upper bound on the fraction of contaminated observations in each coordinate~$j$. Which of~$\tilde{X}_{i,j}$ differ from~$X_{i,j}$ as well as their values can depend on the uncontaminated sample~$X_1,\hdots, X_n$ and external randomization, but at most~$\overline{\eta}_n n$ observations can be altered for each coordinate~$j\in\cbr[0]{1,\hdots,d}$. The coordinatewise contamination restriction in~\eqref{eq:contamfrac} is satisfied whenever~$\envert[0]{\cbr[0]{i\in\cbr[0]{1,\hdots,n}:\tilde{X}_i\neq X_i}}\leq \overline{\eta}_n n$ imposed in~\cite{kock2025high} is. 

Write~$\mc{X}_{n,d}=(\R^d)^n$ and equip~$\mc{X}_{n,d}^2$ with its Borel sigma-algebra. Let~$(X_1,\hdots,X_n)$ denote the projection onto the first~$\mc{X}_{n,d}$ and~$(\tilde{X}_1,\hdots,\tilde{X}_n)$ the projection onto the second~$\mc{X}_{n,d}$. Let~$\P$ be a probability measure on~$\mc X_{n,d}^2$ with first marginal~$P^{\otimes n}$ on~$\mc X_{n,d}$, under which~\eqref{eq:contamfrac} holds. Thus,~$\P$ specifies the distribution of~$\del[0]{X_1,\hdots,X_n,\tilde{X}_1,\hdots,\tilde{X}_n}$. We refer to~\cite{kock2025high} for further discussion and references on adversarial contamination. 

The winsorized estimator considered here is a function of the observed contaminated sample~$\tilde{X}_1,\hdots,\tilde{X}_n$. In this more general setting, let~$\eps_n\in(0,0.5)$ be of the form
\begin{equation}\label{eq:epsfam}
	\eps_n=\lambda_1\cdot \overline{\eta}_n +\lambda_2\cdot \frac{\log(dn)}{n},\qquad \text{with }\lambda_1\in(1,\infty)\text{ and }\lambda_2\in (0,\infty)
\end{equation}  
and use~$\hat{\alpha}_j=\tilde X_{\lceil \eps_n n \rceil,j}^*$  and $\hat{\beta}_j=\tilde X_{\lfloor(1-\eps_n )n\rfloor+1,j}^*$, the order statistics of the~$j$th coordinate of the contaminated data, as the winsorization locations. Therefore, abbreviating~$\mu=\mu(P)$ and~$\Sigma=\Sigma(P)$, $S_{n,W}=S_{n,W}(P)\in\R^d$ has~$j$th coordinate
\begin{equation}\label{eqn:winsmeandef}
	S_{n,W,j}
	=
	n^{-1/2}\sum_{i=1}^n\del[1]{\phi_{\hat\alpha_j,\hat\beta_j}(\tilde{X}_{i,j})-\mu_j},\qquad j=1,\hdots,d,
\end{equation}

No contamination corresponds to~$\overline{\eta}_n=0$ in which case~$\eps_n$ in~\eqref{eq:epsfam} as well as~$S_{n,W}$ in~\eqref{eqn:winsmeandef} reduce to the clean-sample quantities defined in~\eqref{eq:epsfamMain} and~\eqref{eq:WinsMain}, respectively, for~$\lambda_2=\lambda$ (note that~$\lambda_1$ plays no role when~$\overline{\eta}_n=0$ and one can choose~$\lambda_1=2$, say). We shall now present a fixed-$n$ upper bound on
\begin{equation*}
	\rho_{n,W}
	=
	\rho_{n,W}(P)
	:=
	\sup_{H\in\mc H_d}
 \envert[1]{\P(S_{n,W}\in H)-\gamma_{\Sigma}(H)}
\end{equation*}
that is valid uniformly over~$P\in\mc{P}_d(b_1,b_2,m)$, all contamination schemes satisfying~\eqref{eq:contamfrac}, and from which Theorem~\ref{thm:SufficientWins} will follow. To this end, let
 \begin{equation*}
           R_{n,m}^{(1)}:=
            \begin{cases}
                \del[3]{\dfrac{\log^{3-\frac{2}{m}}(dn)}{n^{1-\frac{2}{m}}}}^{1/2}\quad &\text{if}\quad2<m\leq4,\\[3mm]
                \del[3]{\dfrac{\log^5(dn)}{n}}^{1/4}
                +\del[3]{\dfrac{\log^{3-\frac{2}{m}}(dn)}{n^{1-\frac{2}{m}}}}^{1/2}\quad &\text{if}\quad m>4,
            \end{cases}
        \end{equation*}
        and
        \begin{equation*}
        	R_{n,m}^{(2)}:=\del[3]{\frac{\log^{3-\frac{2}{m}}(dn)}{n^{1-\frac{2}{m}}}}^{1/2}\log^{\frac{1}{2}}(n),
        \end{equation*}
        and set
        \begin{equation}
        	R_{n,m}:=\min\del[1]{R_{n,m}^{(1)}, R_{n,m}^{(2)}}\label{eq:Rnm}
        \end{equation}

\begin{theorem}\label{thm:HDGaussSharp}
	Let~$P\in\mc{P}_d(b_1,b_2,m)$ and~\eqref{eq:contamfrac} be satisfied. If~$\eps_n\in(0,1/2)$, with~$\eps_n$ as in~\eqref{eq:epsfam}, then, for a constant~$C$ depending only on~$b_1,b_2,\lambda_1,\lambda_2$, and~$m$, we have
	\begin{equation}\label{eq:Gaussapprox}
		\begin{aligned}
			\rho_{n,W}
						&\leq 	
			C\del[4]{R_{n,m}				+
				\sbr[3]{\overline{\eta}_n^{1-\frac{1}{m}}+\sbr[2]{\frac{\log(dn)}{n}}^{1-\frac{1}{m}}}\sqrt{n\log(d)}}\\
			&\quad +
			C\del[4]{\log^2(d)\sbr[3]{\overline{\eta}_n^{1-\frac{2}{m}}+\sbr[2]{\frac{\log(dn)}{n}}^{1-\frac{2}{m}}}}^{1/2}.
		\end{aligned}
	\end{equation}
	Suppose additionally that~$\sqrt{n\log(d)}\overline{\eta}_n^{1-\frac{1}{m}}\to 0$. Then,
	\begin{equation*}
		\rho_{n,W}\to 0\ 
		\begin{cases}
			\text{for}\quad m\in(2,4]\quad &\text{if}\quad \log(d)/n^{\frac{m-2}{3m-2}}\to 0;\\  
			\text{for}\quad m\in(4,\infty)\quad &\text{if}\quad \log(d)/\del[1]{n^{\frac{m-2}{3m-2}}/\log^{\frac{m}{3m-2}}(n)}\to 0.
		\end{cases}
	\end{equation*}

\end{theorem}

\begin{proof}[Proof of Theorem~\ref{thm:SufficientWins}]
Use~$\lambda_1=2$,~$\lambda_2=\lambda$, and~$\overline{\eta}_n=0$ for all~$n\in\N$ in Theorem~\ref{thm:HDGaussSharp}, noting that~$C$ there does not depend on~$P\in\mc{P}_d(b_1,b_2,m)$.
\end{proof}

\subsection{Further notation and a useful decomposition}\label{sec:appnotation}
We begin by outlining a decomposition of~$S_{n,W}$ used in the proof of Theorem~\ref{thm:HDGaussSharp} along with some further notation used in the sequel. This decomposition is identical to the one in Appendix A in~\cite{kock2025high}. For~$p\in(0,1)$ and a random variable~$X$, denote by~$Q_p(X)$ the $p$-quantile of the distribution of~$X$, that is
	\begin{equation*}
		Q_p(X)=\inf\cbr[1]{z\in \R:\P(X\leq z)\geq p}.
	\end{equation*}
	For~$A_+=1-\lambda_1^{-1}\mathds{1}\cbr[0]{\overline{\eta}_n>0}$ and~$A_-=1+\lambda_1^{-1}\mathds{1}\cbr[0]{\overline{\eta}_n>0}$, define (recall~$\varepsilon_n$ from~\eqref{eq:epsfam})
	\begin{equation}\label{eq:c1}
		c_{1,n}:=-A_+W_0\del[1]{-e^{-(\frac{\log(dn)}{\eps_n n}+A_+)/A_+}}\in (0,A_+) 
	\end{equation}
	as well as
	\begin{equation}\label{eq:c2}
		c_{2,n}:=-A_-W_{-1}\del[1]{-e^{-(\frac{\log(dn)}{\eps_n n}+A_-)/A_-}}\in(A_-,\infty);
	\end{equation}
	here~$W_0$ and~$W_{-1}$ denote the principal and lower branches, respectively, of Lambert's~$W$ function (cf., e.g.,~\cite{Corless1996}). Because~$A_+\leq A_-$, note that~$0<c_{1,n}<c_{2,n}<\infty$. In the \emph{appendix only} we shall sometimes impose
	\begin{align}\label{eq:epscond}
		\mathfrak{T}_n^{(1)}:=2\eps_n +\frac{\log(dn)}{n}+\sqrt{\del[2]{\frac{\log(dn)}{n}}^2+4\frac{\log(dn)}{n}\eps_n}<\frac{1}{2},
	\end{align}
	which, however, is \emph{not} imposed in Theorem~\ref{thm:HDGaussSharp} as this will be argued to hold trivially when~\eqref{eq:epscond} is not satisfied.\footnote{Observe also that~\eqref{eq:epscond} is eventually satisfied under the assumptions of Theorem~\ref{thm:SufficientWins}.} By Lemma~G.2 in~\cite{kock2025high}, one has that~\eqref{eq:epscond} ensures that
	\begin{equation*}
		0
		<
		\eps_n\min(c_{1,n},c_{2,n})
		\leq 
		\eps_n(c_{1,n}+c_{2,n})<1/2.
	\end{equation*}
so that for~$j=1,\hdots,d$
	\begin{equation}\label{eq:alphas}
		\underline{\alpha}_{j}:=Q_{c_{1,n}\eps_n}(X_{1,j})\qquad\text{and}\qquad 	\overline{\alpha}_j:= Q_{c_{2,n}\eps_n}(X_{1,j})
	\end{equation}
	as well as
	\begin{equation}\label{eq:betas}
		\underline{\beta}_j:=Q_{1-c_{2,n}\eps_n
		}(X_{1,j})\qquad \text{and}
		\qquad\overline{\beta}_j:= Q_{1-c_{1,n}\eps_n}(X_{1,j}), 
	\end{equation}
	are well-defined and satisfy~$\underline{\alpha}_{j}\leq\underline{\beta}_j$ and~$\overline{\alpha}_j\leq \overline{\beta}_j$, implying in particular that the decompositions in~\eqref{eq:decomp} and~\eqref{eq:decomp2} below are well-defined. 
	
	By definition,~$\hat{\alpha}_j=\tilde X_{\lceil \eps_n n \rceil,j}^*$ and $\hat{\beta}_j=\tilde X_{\lfloor(1-\eps_n )n\rfloor+1,j}^*$. Lemma~G.3 in~\cite{kock2025high}, which is also valid under the contamination condition in~\eqref{eq:contamfrac}, the union bound, together with obvious monotonicity properties of~$(a,b) \mapsto \phi_{a,b}$, show that with probability at least~$1-\frac{4}{n}$ we have simultaneously for~$j=1,\hdots,d$ (the inequalities $\phi_{\underline{\alpha}_j, \underline{\beta}_j} \leq \phi_{\hat{\alpha}_j, \hat{\beta}_j} \leq \phi_{\overline{\alpha}_j, \overline{\beta}_j}$ and hence)
	\begin{equation}\label{eq:lbub}
		\frac{1}{n}\sum_{i=1}^n\sbr[1]{\phi_{\underline{\alpha}_j,\underline{\beta}_{j}}(\tilde{X}_{i,j})-\mu_j}
		\leq
		\frac{1}{n}\sum_{i=1}^n\sbr[1]{\phi_{\hat\alpha_j,\hat\beta_j}(\tilde{X}_{i,j})-\mu_j}
		\leq
		\frac{1}{n}\sum_{i=1}^n\sbr[1]{\phi_{\overline{\alpha}_j,\overline{\beta}_j}(\tilde{X}_{i,j})-\mu_j}.
	\end{equation}
	The far right-hand side of~\eqref{eq:lbub} can be decomposed as 
	\begin{align}
		\frac{1}{n}\sum_{i=1}^n\sbr[1]{\phi_{\overline{\alpha}_j,\overline{\beta}_j}(\tilde{X}_{i,j})-\mu_j}
		&=
		\underbrace{\frac{1}{n}\sum_{i=1}^n\sbr[1]{\phi_{\overline{\alpha}_j,\overline{\beta}_j}(\tilde{X}_{i,j})-\phi_{\overline{\alpha}_j,\overline{\beta}_j}(X_{i,j})}}_{\overline I_{n,j,1}}\notag\\
		&+
		\underbrace{\frac{1}{n}\sum_{i=1}^n\sbr[1]{\phi_{\overline{\alpha}_j,\overline{\beta}_j}(X_{i,j})-E_P\phi_{\overline{\alpha}_j,\overline{\beta}_j}(X_{i,j})}}_{\overline I_{n,j,2}}\notag\\
		&+\underbrace{\frac{1}{n}\sum_{i=1}^n\sbr[1]{E_P\phi_{\overline{\alpha}_j,\overline{\beta}_j}(X_{i,j})-\mu_j}}_{\overline I_{n,j,3}}\label{eq:decomp},
	\end{align}
	where for a probability measure~$P$ on~$\R^d$, we let~$E_P$ denote the corresponding expectation.
	Similarly, the left-hand side of~\eqref{eq:lbub} can be decomposed as
	\begin{equation}\label{eq:decomp2}
		\frac{1}{n}\sum_{i=1}^n\sbr[1]{\phi_{\underline{\alpha}_j,\underline{\beta}_j}(\tilde{X}_{i,j})-\mu_j}
		=
		\underline I_{n,j,1}+\underline I_{n,j,2}+\underline I_{n,j,3},
	\end{equation}
	with~$\underline I_{n,j,k}$, $j=1,\hdots, d,\ k=1,2,3$, defined analogously to the~$\overline I_{n,j,k}$. Define
	\begin{equation}\label{eq:I_13_max}
		\overline I_{n,k}=\max_{j=1,\hdots, d}\envert[1]{\overline I_{n,j,k}}\qquad\text{and}\qquad \underline I_{n,k}=\max_{j=1,\hdots, d}\envert[1]{\underline I_{n,j,k}},\qquad k=1,3,
	\end{equation}
	as well as
	\begin{equation}\label{eq:I_n2}
		\overline I_{n,2}=\del[1]{\overline I_{n,1,2},\hdots,\overline I_{n,d,2}}'\qquad\text{and}\qquad \underline I_{n,2}=\del[1]{\underline I_{n,1,2},\hdots,\underline I_{n,d,2}}'.	
	\end{equation}
	Throughout, $\overline{\R}=\R\cup\cbr[0]{-\infty,\infty}$. Thus, writing~$Y_{n,j}=n^{-1/2}\sum_{i=1}^n\sbr[1]{\phi_{\hat\alpha_j,\hat\beta_j}(\tilde{X}_{i,j})-\mu_j}$ and~$Y_n=\del[1]{Y_{n,1},\hdots,Y_{n,d}}'$, one has with probability at least~$1-\frac{4}{n}$ that
	\begin{equation}\label{eq:LBUPHD}
		\sqrt{n}(\underline{I}_{n,2}-\underline{I}_{n,1}-\underline{I}_{n,3})
		\leq 
		Y_n
		\leq 
		\sum_{k=1}^3\sqrt{n}\overline I_{n,k},
	\end{equation}
	where \emph{here and in the sequel} (i) all inequalities between vectors are understood elementwise, and (ii) with some abuse of notation we define, for a vector~$x \in \overline{\R}^m$, say, and a real number~$a$, the sum~$x + a$ coordinate-wise as~$(x_1 + a, \hdots, x_m + a)'$ , with the usual convention that~$\infty + a = \infty$ and~$-\infty + a = -\infty$. 
	
	Let~$\varepsilon_n$ be as in~\eqref{eq:epsfam} and recall the quantities defined in~\eqref{eq:alphas} and~\eqref{eq:betas}. Let~$\overline{\Sigma}_{\eps_n}=\overline{\Sigma}_{\eps_n}(P)$ and~$\underline{\Sigma}_{\eps_n}=\underline{\Sigma}_{\eps_n}(P)$ be the matrices with~$(j,k)$th entry
	\begin{equation*}
		\overline{\Sigma}_{\eps_n,j,k}
		=
		\overline{\Sigma}_{\eps_n,j,k}(P)
		=
		E_P\sbr[2]{\del[1]{\phi_{\overline{\alpha}_j,\overline{\beta}_j}(X_{1,j})-E_P \phi_{\overline{\alpha}_j,\overline{\beta}_j}(X_{1,j})}\del[1]{\phi_{\overline{\alpha}_k,\overline{\beta}_k}(X_{1,k})-E_P \phi_{\overline{\alpha}_k,\overline{\beta}_k}(X_{1,k})}} 
	\end{equation*}
	and
	\begin{equation*}
		\underline{\Sigma}_{\eps_n,j,k}
		=
		\underline{\Sigma}_{\eps_n,j,k}(P)
		=
		E_P\sbr[2]{\del[1]{\phi_{\underline{\alpha}_j,\underline{\beta}_j}(X_{1,j})-E_P \phi_{\underline{\alpha}_j,\underline{\beta}_j}(X_{1,j})}\del[1]{\phi_{\underline{\alpha}_k,\underline{\beta}_k}(X_{1,k})-E_P \phi_{\underline{\alpha}_k,\underline{\beta}_k}(X_{1,k})}}, 
	\end{equation*}
	respectively. By Lemma B.1 in~\cite{kock2025high} there exists a positive constant~$C$ depending only on~$b_2$,~$\lambda_2$, and~$m$ such that under~\eqref{eq:epscond}
	
	\begin{equation}\label{eq:truncvar}
			\max_{j,k\in\cbr[0]{1,\hdots,d}}\del[2]{\envert[1]{\overline{\Sigma}_{\eps_n,j,k}-\Sigma_{j,k}}\vee\envert[1]{\underline{\Sigma}_{\eps_n,j,k}-\Sigma_{j,k}}}
			\leq
			C\eps_n^{1-\frac{2}{m}}=:\mathfrak{T}_n^{(2)}.
		\end{equation}

	\begin{lemma}\label{lem:HDGaussBD}
		Let~$P\in\mc{P}_d(b_1,b_2,m)$ and~$\eps_n$ be chosen as in~\eqref{eq:epsfam} satisfying~\eqref{eq:epscond} and suppose that there exists a strictly positive constant~$\mathfrak{b}_1$ such that for all~$j=1,\hdots,d$
		\begin{equation*}
			E_P\sbr[1]{\phi_{\underline{\alpha}_j,\underline{\beta}_j}(X_{1,j})-E_P\phi_{\underline{\alpha}_j,\underline{\beta}_j}(X_{1,j})}^2\geq\mathfrak{b}_1^2\quad
			\text{and}
			\quad E_P\sbr[1]{\phi_{\overline{\alpha}_j,\overline{\beta}_j}(X_{1,j})-E_P\phi_{\overline{\alpha}_j,\overline{\beta}_j}(X_{1,j})}^2\geq \mathfrak{b}_1^2.
		\end{equation*}
		Then, recalling the definitions of~$\overline{I}_{n,2}$ and~$\underline{I}_{n,2}$ in~\eqref{eq:I_n2} and~$R_{n,m}$ in~\eqref{eq:Rnm},
		\begin{equation*}
			\sup_{H\in \mc{H}_d}\envert[2]{P^{\otimes n}\del[1]{\sqrt{n}\overline{I}_{n,2}\in H}-\gamma_{\overline{\Sigma}_{\eps_n}}(H)}
			\ \bigvee\ \sup_{H\in\mc{H}_d}\envert[2]{P^{\otimes n}\del[1]{\sqrt{n}\underline{I}_{n,2}\in H}-\gamma_{\underline{\Sigma}_{\eps_n}}(H)}
			\leq
			CR_{n,m}
		\end{equation*}
		where~$C$ is a constant depending only on~$\mathfrak{b}_1$,~$b_2$,~$\lambda_1$,~$\lambda_2$, and~$m$.
	\end{lemma}
	\begin{proof}
		We establish~$A_n:=\sup_{H\in \mc{H}_d}\envert[0]{P^{\otimes n}\del[0]{\sqrt{n}\overline{I}_{n,2}\in H}-\gamma_{\overline{\Sigma}_{\eps_n}}(H)}\leq CR_{n,m}$ by showing that the left-hand side is bounded from above by~$C R_{n,m}^{(1)}$ and~$C R_{n,m}^{(2)}$ (for potentially different constants~$C$ depending only on $\mathfrak{b}_1$,~$b_2$,~$\lambda_1$,~$\lambda_2$, and~$m$). By an analogous argument, it also holds that~$\sup_{H\in \mc{H}_d}\envert[0]{P^{\otimes n}\del[0]{\sqrt{n}\underline{I}_{n,2}\in H}-\gamma_{\underline{\Sigma}_{\eps_n}}(H)}\leq CR_{n,m}$.
		
	 For the purpose of proving~$A_n\leq CR_{n,m}^{(k)}$ for each~$k\in\cbr[0]{1,2}$, we note some properties of
	 \begin{equation*}
	 	U_{i,j}:=\phi_{\overline{\alpha}_j,\overline{\beta}_j}(X_{i,j})-E_P\phi_{\overline{\alpha}_j,\overline{\beta}_j}(X_{i,j}).
	 \end{equation*}
	 
	 First, for~$i=1,\hdots,n$ and~$j=1,\hdots,d$
		\begin{align*}
			U_{i,j}=\phi_{\overline{\alpha}_j,\overline{\beta}_j}(X_{i,j})-E_P\phi_{\overline{\alpha}_j,\overline{\beta}_j}(X_{i,j})
			\leq
			\overline{\beta}_j-\overline{\alpha}_j=Q_{1-c_{1,n}\eps_n}(X_{1,j})-Q_{c_{2,n}\eps_n}(X_{1,j}),	
		\end{align*}
		Furthermore, by Lemma~G.1 in~\cite{kock2025high} and~$\eps_n\geq \lambda_2\frac{\log(dn)}{n}$, the far right-hand side of the previous display is upper bounded by
		\begin{align*}
			\frac{b_2}{(c_{1,n}\eps_n)^{1/m}}+\frac{b_2}{(c_{2,n}\eps_n)^{1/m}}
			\leq
			C_1\del[3]{\frac{n}{\log(dn)}}^{1/m}=:B_n,
		\end{align*}
		where, by Lemma~G.2 in~\cite{kock2025high}, $C_1$ is a constant depending only on~$b_2$,~$\lambda_1$,~$\lambda_2$ and~$m$. Similarly, 
		\begin{align*}
			U_{i,j}=\phi_{\overline{\alpha}_j,\overline{\beta}_j}(X_{i,j})-E_P\phi_{\overline{\alpha}_j,\overline{\beta}_j}(X_{i,j})
			\geq
			\overline{\alpha}_j-\overline{\beta}_j
			=
			Q_{c_{2,n}\eps_n}(X_{1,j})-Q_{1-c_{1,n}\eps_n}(X_{1,j})
			\geq
			-B_n.
		\end{align*}
		Thus, by construction,~$\max_{1\leq j\leq d}\max_{1\leq i\leq n}|U_{i,j}|\leq B_n$ and by assumption,~$E_PU_{1,j}^2\geq \mathfrak{b}_1^2$.        
        Furthermore, for~$2<m\leq 4$,
        \begin{equation*}
        	E_P\sbr[1]{\phi_{\overline{\alpha}_j,\overline{\beta}_j}(X_{1,j})-E_P\phi_{\overline{\alpha}_j,\overline{\beta}_j}(X_{1,j})}^4
        	\leq
        	E_P\envert[1]{\phi_{\overline{\alpha}_j,\overline{\beta}_j}(X_{1,j})-E_P\phi_{\overline{\alpha}_j,\overline{\beta}_j}(X_{1,j})}^mB_n^{4-m}
        	\leq
        	b_2^mB_n^{4-m},
        \end{equation*} 
        the last inequality following from, e.g., Corollary 3 in~\cite{cstudd}. In case~$m>4$, it holds by Corollary 3 in~\cite{cstudd} and Lyapunov's inequality yield that $E_P\sbr[1]{\phi_{\overline{\alpha}_j,\overline{\beta}_j}(X_{1,j})-E_P\phi_{\overline{\alpha}_j,\overline{\beta}_j}(X_{1,j})}^4\leq b_2^4$. Thus, writing~$r=m\wedge 4$,
        \begin{equation*}
        		EU_{1,j}^4
    			 =
        		E_P\sbr[1]{\phi_{\overline{\alpha}_j,\overline{\beta}_j}(X_{1,j})-E_P\phi_{\overline{\alpha}_j,\overline{\beta}_j}(X_{1,j})}^4
        		\leq
        		b_2^rB_n^{4-r}.
        \end{equation*}

		We can now apply (A.2) of supplementary Lemma A.1 in~\cite{chernozhuokov2022improved} to show that $\sup_{H\in \mc{H}_d}\envert[0]{P^{\otimes n}\del[0]{\sqrt{n}\overline{I}_{n,2}\in H}-\gamma_{\overline{\Sigma}_{\eps_n}}(H)}\leq CR_{n,m}^{(1)}$ for the desired constant~$C$: Writing~$U_i=(U_{i,1},\hdots,U_{i,d})'$, it trivially holds that
		\begin{equation*}
			|-U_{i,j}|=|U_{i,j}|\leq B_n,\qquad~E_P(-U_{i,j})^2=E_PU_{i,j}^2\qquad\text{and}\qquad E_P(-U_{i,j})^4=E_PU_{i,j}^4.
		\end{equation*}
  Hence, by (A.2) of Lemma A.1 in~\cite{chernozhuokov2022improved} with~$X_i$ there being~$(U_i',-U_i')'\in\R^{2d}$,~$\Delta_0=0$,~$\Delta_1\leq n^{-1}b_2^rB_n^{4-r}$,~$\kappa_n=B_n$,~$c=\mathfrak{b}_1^2$,~$p=2d$, and~$\Sigma$ there being the covariance matrix of $(U_1',-U_1')'$ (having upper-left block~$\overline{\Sigma}_{\eps_n}$) yields the existence of a constant~$C_2$ depending only on~$\mathfrak{b}_1$,~$b_2$ and~$m$ such that
  \begin{align*}
            \sup_{H\in\mc{H}_d^f}\envert[2]{P^{\otimes n}\del[1]{\sqrt{n}\overline{I}_{n,2}\in H}-\gamma_{\overline{\Sigma}_{\eps_n}}(H)}
            &\leq 
            C_2\del[4]{\del[3]{\frac{B_n^{4-r}\log^5(dn)}{n}}^{1/4}
            +\del[3]{\frac{B_n^2\log^3(dn)}{n}}^{1/2}}\\
            &\leq 
            C_3R_{n,m}^{(1)},
        \end{align*}	
        where~$C_3$ depends only on~$\mathfrak{b}_1$,~$b_2$,~$\lambda_1$,~$\lambda_2$, and~$m$ and~$\mc{H}_d^f$ is the collection of hyperrectangles with finite endpoints. The last inequality used the definition of~$B_n$ and~$R_{n,m}$ in~\eqref{eq:Rnm}, separately for~$2<m\leq4$ and~$m>4$. For unrestricted~$H\in\mc{H}_d$, non-empty, and~$M>0$, write~$H_{M}:=H\cap [-M,M]^d$, which for~$M$ sufficiently large is a non-empty element of~$\mc{H}_d^f$. It follows from upward continuity of probability measures that
\begin{equation*}
	\envert[2]{P^{\otimes n}\del[1]{\sqrt{n}\overline{I}_{n,2}\in H}-\gamma_{\overline{\Sigma}_{\eps_n}}(H)}
	=
	\lim_{M\to\infty}\envert[2]{P^{\otimes n}\del[1]{\sqrt{n}\overline{I}_{n,2}\in H_M}-\gamma_{\overline{\Sigma}_{\eps_n}}(H_M)}
	\leq
	 C_3R_{n,m}^{(1)}. 
\end{equation*}
Since~$H\in\mc{H}_d$ (non-empty) was arbitrary, and~$C_3R_{n,m}^{(1)}$ does not depend on~$H$, as well as~$H=\emptyset$ being trivial, this establishes~$A_n\leq C_3R_{n,m}^{(1)}$.
        
      We next invoke Corollary 1 of~\cite{CCKK26} to show that~$A_n\leq C R_{n,m}^{(2)}$. To this end, let~$\sigma_{n,j}^2:=E_PU_{1,j}^2\geq \mathfrak{b}_1^2$ and~$D_n=\text{diag}(\sigma_{n,1},\hdots, \sigma_{n,d})$. For any~$H\in \mc{H}_d$ and~$x\in\R^d$ one has that~$x\in H$ if and only if~$D_n^{-1}x\in H_{D_n}:=D_n^{-1}H\in\mc{H}_d$ so that (recall the definition of~$\overline{I}_{n,2}$ in~\eqref{eq:I_n2}) 
         \begin{equation*}
         \envert[2]{P^{\otimes n}\del[1]{\sqrt{n}\overline{I}_{n,2}\in H}-\gamma_{\overline{\Sigma}_{\eps_n}}(H)}
         =	
         \envert[3]{P^{\otimes n}\del[2]{n^{-1/2}\sum_{i=1}^n[D_n^{-1}U_i]\in H_{D_n}}-\gamma_{\overline{\Sigma}_{\eps_n,0}}(H_{D_n})},
         \end{equation*}
         where~$\overline{\Sigma}_{\eps_n,0}:=D_n^{-1}\overline{\Sigma}_{\eps_n}D_n^{-1}$.
         Thus,~$D_n^{-1}U_1$ has correlation matrix~$\overline{\Sigma}_{\eps_n,0}$; in particular its variances are one. Applying Corollary 1 of~\cite{CCKK26} with~$X_i$ there equal to~$D_n^{-1}U_i$ and~$B=\mathfrak{b}_1^{-1}B_n$, and using~$g_n\leq C'\sqrt{\log(d)}$ in the notation of that corollary for a universal constant~$C'>0$ yields that the right-hand side of the previous display is bounded from above by 
         \begin{equation*}
         	C_4\del[3]{\frac{B_n\log^{3/2}(d)}{n^{1/2}}}\sqrt{1+\log_+\del[2]{\frac{n^{1/2}}{B_n\log^{1/2}(d)}}}
         	\leq
         	C_4\del[3]{\frac{B_n^2\log^3(d)}{n}}^\frac{1}{2}\log^{\frac{1}{2}}(n) 
         	\leq
         	C_4R_{n,m}^{(2)}   
         	 \end{equation*}
         	 for a constant~$C_4$ depending at most on $\mathfrak{b}_1$,~$b_2$,~$\lambda_1$,~$\lambda_2$, and~$m$ (potentially changing values across the inequalities); here~$\log_+(z)=0\vee \log(z)$ for~$z\in(0,\infty)$, and the first inequality in the previous display follows from ($C_5$ being a constant depending only on~$\mathfrak{b}_1$, $b_2$,~$\lambda_1$,~$\lambda_2$ and~$m$ but changing values across inequalities)
         	 \begin{equation*}
         	 \frac{n^{1/2}}{B_n\log^{1/2}(d)}
         	 \leq
         	 C_5n^{\frac{1}{2}-\frac{1}{m}}\frac{\log^{1/m}(dn)}{\log^{1/2}(d)}
         	 \leq
         	 C_5	n^{\frac{1}{2}-\frac{1}{m}}\del[3]{\frac{\log(dn)}{\log(d)}}^{\frac{1}{2}}
         	 \leq
         	 C_5n^{\frac{1}{2}-\frac{1}{m}}\log^{1/2}(n);
         	 \end{equation*}
         	 the last inequality considering separately the case of~$d\leq n$ and~$d>n$.
         	 
         	 Therefore,~$A_n\leq C_4 R_{n,m}^{(2)}$ since~$H\in\mc{H}_d$ was arbitrary. That~$A_n\leq CR_{n,m}$ follows with~$C=\max(C_3,C_4)$.
                	\end{proof}

\begin{lemma}\label{lem:covcompare}
		Let~$m>2$ and assume that~$\sqrt{n\log(d)}\overline{\eta}_n^{1-\frac{1}{m}}\to 0$ as well as~$\log(d)/n^{\frac{m-2}{3m-2}}\to 0$. Then, 		%
		\begin{equation}\label{eq:conv1}
			\log^2(d)\del[3]{\overline{\eta}_n^{1-\frac{2}{m}}+\sbr[2]{\frac{\log(dn)}{n}}^{1-\frac{2}{m}}}\to 0,
		\end{equation}	
		and
		\begin{equation}\label{eq:conv2}
			\del[3]{\overline{\eta}_n^{1-\frac{1}{m}}+\sbr[2]{\frac{\log(dn)}{n}}^{1-\frac{1}{m}}}\sqrt{n\log(d)}\to 0.
		\end{equation}
		Finally, for~$R_{n,m}$ as defined in~\eqref{eq:Rnm}, it holds that~$R_{n,m}\to 0$ for~$m\in(2,4]$ whereas~$R_{n,m}\to 0$ for~$m\in(4,\infty)$ if, in addition,~$\log(d)/\del[1]{n^{\frac{m-2}{3m-2}}/\log^{\frac{m}{3m-2}}(n)}\to 0$.
	\end{lemma}
	\begin{proof}
	We begin by proving~\eqref{eq:conv1}, to which end we first establish that $\log^2(d)\overline{\eta}_n^{1-\frac{2}{m}}\to 0$. Because~$\sqrt{n\log(d)}\overline{\eta}_n^{1-\frac{1}{m}}\to 0$, it holds that~$\overline{\eta}_n=o\del[1]{1/[n\log(d)]^{\frac{m}{2(m-1)}}}$. Thus,
		\begin{equation*}
		\log^2(d)\overline{\eta}_n^{1-\frac{2}{m}}
		=
		o\del[4]{\frac{\log^2(d)}{[n\log(d)]^{(m-2)/(2(m-1))}}}
		=
		o\del[4]{\sbr[3]{\frac{\log(d)}{n^{\frac{m-2}{3m-2}}}}^{\frac{3m-2}{2(m-1)}}}
		=
		o(1).	
		\end{equation*} 
		Next,
		\begin{equation*}
			\log^2(d)\sbr[2]{\frac{\log(d)}{n}}^{1-\frac{2}{m}}
			=
			\frac{\log(d)^{3-\frac{2}{m}}}{n^{1-\frac{2}{m}}}
			\to 0
			\iff 
			\frac{\log(d)}{n^{\frac{m-2}{3m-2}}}\to 0
			,
		\end{equation*}
		the latter convergence being by assumption. Finally, 
		\begin{equation*}
			\log^2(d)\sbr[2]{\frac{\log(dn)}{n}}^{1-\frac{2}{m}}\to 0
		\end{equation*}
		by a standard subsequence argument, considering separately the cases of subsequences along which i)~$n\leq d$ for which the convergence follows from the penultimate display and ii)~$d< n$ for which the convergence follows from~$m>2$. This establishes~\eqref{eq:conv1}.
		
		To prove~\eqref{eq:conv2}, note that~$\sqrt{n\log(d)}\overline{\eta}_n^{1-\frac{1}{m}}\to 0$ by assumption and it remains to prove 
		\begin{equation*}
			\sqrt{n\log(d)}\del[2]{\frac{\log(dn)}{n}}^{1-\frac{1}{m}}\to 0.
		\end{equation*}
		To this end, note that
		\begin{equation*}
			\sqrt{n\log(d)}\del[2]{\frac{\log(d)}{n}}^{1-\frac{1}{m}}
			=
			\frac{\sbr[1]{\log(d)}^{\frac{3}{2}-\frac{1}{m}}}{n^{\frac{1}{2}-\frac{1}{m}}}\to 0
			\iff
			\frac{\log(d)}{n^{\frac{m-2}{3m-2}}}\to 0,	
		\end{equation*}
		where the latter convergence is by assumption. The convergence in the penultimate display now follows by considering separately subsequences along which~$n\leq d$ and~$d<n$, respectively (as in the concluding argument of the proof of~\eqref{eq:conv1}). 
			
		Finally, we turn to the statements concerning~$R_{n,m}\to 0$ without repeating subsequence arguments. For~$m\in(2,4]$, use that~$R_{n,m}\leq R_{n,m}^{(1)}\to 0$, the convergence being by~$\log(d)/n^{\frac{m-2}{3m-2}}\to 0$. For~$m\in(4,\infty)$, use that~$R_{n,m}\leq R_{n,m}^{(2)}$,	write	
		\begin{equation*}
			R_{n,m}^{(2)}
			=
			\del[3]{\frac{\log^{3-\frac{2}{m}}(dn)}{n^{1-\frac{2}{m}}/\log(n)}}^{1/2}		
		\end{equation*}
		and note that
		\begin{equation*}
			\del[3]{\frac{\log^{3-\frac{2}{m}}(d)}{n^{1-\frac{2}{m}}/\log(n)}}^{1/2}	
			=
			\del[3]{\frac{\log(d)}{n^{\frac{m-2}{3m-2}}/\log^{\frac{m}{3m-2}}(n)}}^{\frac{1}{2}(3-\frac{2}{m})}	
			\to 0,
		\end{equation*}
		the convergence being by assumption.
	\end{proof}

Let~$\mc{O}_d$ be the collection of lower orthant sets of the form
\begin{equation*}
	O=\cbr[1]{x\in\R^d:x_j\leq t_j\text{ for all }j=1,\hdots,d}
\end{equation*}
 for some~$t\in\overline{\R}^d$. The following lemma shows that it suffices to establish Theorem~\ref{thm:HDGaussSharp} over~$\mc{O}_d$. 

\begin{lemma}\label{lem:onesidedtorectangle}
		Let~$P\in\mc{P}_d(b_1,b_2,m)$ and suppose Theorem~\ref{thm:HDGaussSharp} holds for~$\rho_{n, W}$ throughout replaced by
		\begin{equation}
			\rho_{n, W}^{os} := \sup_{O\in\mc{O}_d}\envert[1]{\P\del[1]{S_{n,W}\in O}-\gamma_{\Sigma}\del[0]{O}},
		\end{equation}
		and with~$\lambda_2$ in the definition of~$\eps_n$ in~\eqref{eq:epsfam} replaced by~$\lambda_{2,n}\in[\lambda_2/2,\lambda_2]$. Then, Theorem~\ref{thm:HDGaussSharp} holds (with a different constant~$C$ depending only on~$b_1,b_2,\lambda_1,\lambda_2$, and~$m$).
	\end{lemma}
	
	\begin{proof}
		Fix~$a$ and~$b$ in~$\overline{\R}^d$ such that~$-\infty \leq a_j\leq b_j\leq \infty$ and let~$H=\bigtimes_{j=1}^d [a_j,b_j]\in\mc{H}_d$. Clearly,
		\begin{equation}\label{eq:twodim1}
			\P\del[1]{S_{n,W}\in H}
			=
			\P\del[1]{(S_{n,W}',-S_{n,W}')'\leq (b',-a')'}.
		\end{equation}
		Furthermore, by~\eqref{eqn:winsmeandef}
		\begin{equation*}
			-S_{n,W,j}
			=
			n^{-1/2}\sum_{i=1}^n\del[1]{\phi_{-\hat\beta_j,-\hat\alpha_j}(-\tilde{X}_{i,j})-(-\mu_j)}
			,\qquad j=1,\hdots,d.
		\end{equation*}
		Thus, the vector~$(X_i',-X_i')'\in\R^{2d}$ has distribution~$Q=P\circ \del[1]{x\mapsto (x,-x)}^{-1}\in\mc{P}_{2d}(b_1,b_2,m)$ and~$(\tilde{X}_i',-\tilde{X}_i')'\in\R^{2d}$ satisfies~\eqref{eq:contamfrac} with~$d$ there replaced by~$2d$ (but with the same~$m$, $b_1$,~$b_2$ and~$\overline{\eta}_n$). Note that the covariance matrix of~$(X_1',-X_1')'$ is
		\begin{equation*}
			\Xi= 
			\begin{pmatrix}
				\Sigma & -\Sigma\\
				-\Sigma & \Sigma
			\end{pmatrix}.
		\end{equation*}
	 Note also that for fixed~$\lambda_1\in(1,\infty)$ and~$\lambda_2\in(0,\infty)$,~$\eps_n$ defined in~\eqref{eq:epsfam} can be written as
		\begin{equation*}
			\eps_n
			=			
			\lambda_1\cdot \overline{\eta}_n +\lambda_{2,n}\cdot \frac{\log(2dn)}{n},\qquad\text{where}\qquad\lambda_{2,n}=\lambda_2\cdot\frac{\log(dn)}{\log(2dn)}\in[\lambda_2/2,\lambda_2],
		\end{equation*}
		the inclusion following from~$d\geq 2$ and~$n>3$. Hence, the version of Theorem~\ref{thm:HDGaussSharp} that replaces~$\rho_{n, W}$ by~$\rho_{n, W}^{os}$ and~$\lambda_2$ by~$\lambda_{2,n}\in[\lambda_2/2,\lambda_2]$ (and which is assumed to hold true) applies to~$(S_{n,W}',-S_{n,W}')'$ and yields, for~$C^*$ a constant depending only on~$b_1,b_2,\lambda_1,\lambda_2$, and~$m$, the upper bound (recall~$R_{n,m}$ is defined in~\eqref{eq:Rnm}) 
		\begin{align*}
			C^*\del[4]{R_{n,m}
				+
				\sbr[3]{\overline{\eta}_n^{1-\frac{1}{m}}+\sbr[2]{\frac{\log(dn)}{n}}^{1-\frac{1}{m}}}\sqrt{n\log(d)}} 
			+
			C^*\del[4]{\log^2(d)\sbr[3]{\overline{\eta}_n^{1-\frac{2}{m}}+\sbr[2]{\frac{\log(dn)}{n}}^{1-\frac{2}{m}}}}^{1/2}
		\end{align*}
		on (noting also that~$[-\infty,b]\times [-\infty,-a]\in\mc{O}_{2d}$)
		\begin{equation}\label{eq:twodim2}
			\envert[1]{\P\del[1]{(S_{n,W},-S_{n,W})\leq (b',-a')'}-\gamma_{\Xi}\del[1]{[-\infty,b]\times [-\infty,-a]}} ,
		\end{equation}
		Furthermore,~$\gamma_{\Xi}$ is the image measure of~$\gamma_{\Sigma}$ under the mapping~$\R^d\ni z\mapsto (z,-z)\in\R^{2d}$, that is~$\gamma_{\Xi}=\gamma_{\Sigma}\circ \del[1]{z\mapsto (z,-z)}^{-1}$. Thus,
		\begin{align}\label{eq:twodim3}
			\gamma_{\Xi}\del[1]{[-\infty,b]\times [-\infty,-a]}
			=
			\gamma_\Sigma\del[1]{z\in\R^d: (z,-z)\in [-\infty,b]\times [-\infty,-a]}\
			=
		\gamma_{\Sigma}(H)
		\end{align}
		Combining~\eqref{eq:twodim1} and~\eqref{eq:twodim3} with the upper bound in~\eqref{eq:twodim2} obtained above (which does not depend on~$a$ or~$b$) delivers the claim. 
	\end{proof}
	
	\subsection{Proof of Theorem~\ref{thm:HDGaussSharp}}
	By Lemma~\ref{lem:onesidedtorectangle} it suffices to prove~\eqref{eq:Gaussapprox} for hyperrectangles~$H$ of the form~$H=\bigtimes_{j=1}^d[-\infty,t_j]$ with~$ t_j\in \overline{\R}$ and~$\eps_n$ as defined in~\eqref{eq:epsfam}, but with~$\lambda_2$ replaced by~$\lambda_{2,n}\in[\lambda_2/2,\lambda_2]$. Throughout the proof, we consider the case where~$\lambda_{2,n}=\lambda_2$ for notational convenience. Inspection of the argument shows that it remains valid for~$\lambda_{2,n}\in[\lambda_2/2,\lambda_2]$. Furthermore, recalling the definitions of~$\mathfrak{T}_n^{(1)}$ and~$\mathfrak{T}_n^{(2)}$ in~\eqref{eq:epscond} and \eqref{eq:truncvar}, respectively, and using that~$\log(dn)/n\leq \lambda_2^{-1}\eps_n$, one has with~$C$ being the constant from~\eqref{eq:truncvar} that
	\begin{equation}\label{eq:auxgood}
		\mathfrak{T}_n^{(1)}+\mathfrak{T}_n^{(2)}
		\leq
		2\eps_n+\lambda_2^{-1}\eps_n+\sqrt{\lambda_2^{-2}\eps_n^2+4\lambda_2^{-1}\eps_n^2}+C\eps_n^{1-\frac{2}{m}}
		\leq
		C_1\eps_n^{1-\frac{2}{m}},
	\end{equation}
	for~$C_1$ a constant depending only on~$b_2$,~$\lambda_2$~and~$m$, and where we used that~$\eps_n\leq \eps_n^{1-\frac{2}{m}}$ as~$\eps_n\in(0,0.5)$. In the following we first suppose that
	\begin{equation}\label{eq:CaseGood}
		C_1\eps_n^{1-\frac{2}{m}}
		<
		\frac{(b_1^2\wedge 1)}{2},
	\end{equation}
	which, by the penultimate display, implies that~$\mathfrak{T}_n^{(1)}<1/2$ such that~$\eqref{eq:epscond}$ is satisfied. By Lemma~G.4 of~\cite{kock2025high}, which remains valid under the coordinate-wise contamination condition in~\eqref{eq:contamfrac}, there exists a constant~$C_2$ depending only on~$b_2$, $\lambda_1$,~$\lambda_2$, and~$m$ such that
	\begin{equation}\label{eq:I1}
		\underline{I}_{n,1}\vee \overline{I}_{n,1}
		\leq
		C_2\overline{\eta}_n^{1-\frac{1}{m}}
		:=I_{n,1}
	\end{equation}
	with~$I_{n,1}$ being non-random. Thus, from~\eqref{eq:LBUPHD} it holds with probability at least~$1-\frac{4}{n}$ that
	\begin{equation*}
		\underline{Y}_{n}:=\sqrt{n}\underline{I}_{n,2}-\sqrt{n}(I_{n,1}+\underline{I}_{n,3})
		\leq 
		Y_n
		\leq 
		\sqrt{n}\overline{I}_{n,2}+\sqrt{n}(I_{n,1}+\overline{I}_{n,3})=:\overline{Y}_n.
	\end{equation*}
	Next, it holds for all~$t=(t_1,\hdots,t_d)$ with~$t_j\in \overline{\R}$ that
	\begin{equation*}
		\P\del[1]{Y_n\leq t}
		\geq
		\P\del[1]{\overline{Y}_n\leq t, Y_n\leq \overline{Y}_n}
		\geq
		\P\del[1]{\overline{Y}_n\leq t}+\P\del[1]{Y_n\leq \overline{Y}_n}-1
		\geq
		\P\del[1]{\overline{Y}_n\leq t}-\frac{4}{n}.
	\end{equation*}
	Similarly,
	\begin{equation*}
		\P\del[1]{Y_n\leq t}
		\leq
		\P\del[1]{\underline{Y}_n\leq t}+\P\del[1]{\underline{Y}_n\not\leq Y_n}
		\leq	
		\P\del[1]{\underline{Y}_n\leq t}+\frac{4}{n}.
	\end{equation*}
	Thus, recalling that~$H=\bigtimes_{j=1}^d[-\infty,t_j]$,
	\begin{align}\label{eq:1overn}
		\sup_{t\in\overline{\R}^d}\envert[2]{\P\del[1]{Y_n\leq t}-\gamma_{\Sigma}(H)}	
		&\leq
		\sup_{t\in\overline{\R}^d}\envert[2]{\P\del[1]{\underline{Y}_n\leq t}-\gamma_{\Sigma}(H)}\notag\\
		&+\sup_{t\in\overline{\R}^d}\envert[2]{\P\del[1]{\overline{Y}_n\leq t}-\gamma_{\Sigma}(H)}
		+\frac{4}{n}.
	\end{align}
	We proceed by bounding the second summand on the right-hand side of the previous display (the argument for the first summand is analogous and hence skipped). To this end, note that~\eqref{eq:auxgood} and~\eqref{eq:CaseGood} imply that~$\mathfrak{T}_n^{(2)}\leq b_1^2/2$ which then together with \eqref{eq:truncvar} implies that~$\min_{j=1,\hdots, d}\overline{\Sigma}_{\eps_n,j,j}\geq b_1^2/2$. Next, by the definition of~$\overline{Y}_n$
	\begin{equation*}
		\P\del[1]{\overline{Y}_n\leq t}
		=
		P^{\otimes n}\del[1]{\sqrt{n} \overline{I}_{n,2}\leq t-\sqrt{n}\del[0]{I_{n,1}+\overline{I}_{n,3}}},\qquad t\in \overline{\R}^d,
	\end{equation*}
	whence, noting also that~$I_{n,1}$ and~$\overline{I}_{n,3}$ are non-random, Lemma~\ref{lem:HDGaussBD} applied with~$\mathfrak{b}_1^2=b_1^2/2$ implies the existence of a constant~$C_3$ depending only on~$b_1,b_2,\lambda_1,\lambda_2$, and~$m$ such that
	\begin{equation*}
		\sup_{t\in\overline{\R}^d}\envert[3]{\P\del[1]{\overline{Y}_n\leq t}-\P^*\del[3]{\overline{Z}\leq t-\sqrt{n}\del[0]{I_{n,1}+\overline{I}_{n,3}}}}
		\leq 
		C_3R_{n,m},	
	\end{equation*}
	where~$\overline{Z}\sim \mathsf{N}_d(0,\overline{\Sigma}_{\eps_n})$ is defined on an auxiliary probability space with probability measure~$\P^*$.
	Furthermore, by the Gaussian anti-concentration inequality as stated in Theorem 1 of~\cite{chernozhukov2017detailed}, (cf.~also Lemma A.1 in~\cite{chernozhukov2017central}),
	\begin{align*}
		0
		&\geq
		\P^*\del[1]{\overline{Z}\leq t-\sqrt{n}\del[0]{I_{n,1}+\overline{I}_{n,3}}}-\P^*\del[1]{\overline{Z}\leq t}\\
		&\geq
		-\frac{\sqrt{n}\del[0]{I_{n,1}+\overline{I}_{n,3}}}{\sqrt{\min_{j=1,\hdots, d}\overline{\Sigma}_{\eps_n,j,j}}}\del[1]{\sqrt{2\log(d)}+2}\\
		&\geq
		-C_4\del[3]{\overline{\eta}_n^{1-\frac{1}{m}}+\sbr[2]{\frac{\log(dn)}{n}}^{1-\frac{1}{m}}}\sqrt{n\log(d)},
	\end{align*}
	the final inequality following from Lemma~B.3 in~\cite{kock2025high}, \eqref{eq:I1}, $\min_{j=1,\hdots, d}\overline{\Sigma}_{\eps_n,j,j}\geq b_1^2/2$, and~$C_4$ being a constant depending on~$b_1$, $b_2$,~$\lambda_1,\lambda_2$, and~$m$ only. Thus, combining the previous two displays, there exists a constant~$C_5$ depending on~$b_1,b_2,\lambda_1,\lambda_2$, and~$m$ only such that when~\eqref{eq:CaseGood} is satisfied one has
	\begin{align*}
		\sup_{t\in\overline{\R}^d}\envert[2]{\P\del[1]{\overline{Y}_n\leq t}-\P^*\del[1]{\overline{Z}\leq t}}
		\leq
		C_5\del[4]{R_{n,m}
			+
			\sbr[3]{\overline{\eta}_n^{1-\frac{1}{m}}+\sbr[2]{\frac{\log(dn)}{n}}^{1-\frac{1}{m}}}\sqrt{n\log(d)}}.
	\end{align*}
	Finally, by the Gaussian-to-Gaussian comparison inequality as stated in Proposition 2.1 of~\cite{chernozhuokov2022improved} and~\eqref{eq:truncvar}, there exists a constant~$C_6$, depending on~$b_1$,~$b_2$,~$\lambda_1$, $\lambda_2$, and~$m$ only ($C_6$ changes value in the second inequality below), such that for~$Z\sim \mathsf{N}_d(0,\Sigma)$ defined on the same probability space as~$\overline{Z}$
	\begin{align*}
		\sup_{t\in\overline{\R}^d}\envert[2]{\P^*\del[1]{\overline{Z}\leq t}-\P^*\del[1]{Z\leq t}}	
		&\leq
		C_6\del[2]{\log^2(d)\eps_n^{1-\frac{2}{m}}}^{1/2}\\
		&\leq
		C_6\del[4]{\log^2(d)\sbr[3]{\overline{\eta}_n^{1-\frac{2}{m}}+\sbr[2]{\frac{\log(dn)}{n}}^{1-\frac{2}{m}}}}^{1/2}.
	\end{align*}
	The previous two displays yield that there exists a constant~$C_7$ depending only on~$b_1$,~$b_2$,~$\lambda_1,\lambda_2$, and~$m$ such that when~\eqref{eq:CaseGood} is satisfied it holds that
	\begin{align*}
		\sup_{t\in\overline{\R}^d}\envert[2]{\P\del[1]{\overline{Y}_n\leq t}-\P^*\del[1]{Z\leq t}}
		\leq
		C_7\del[4]{\log^2(d)\sbr[3]{\overline{\eta}_n^{1-\frac{2}{m}}+\sbr[2]{\frac{\log(dn)}{n}}^{1-\frac{2}{m}}}}^{1/2}\\
		+\ C_7\del[4]{R_{n,m}
			+
			\sbr[3]{\overline{\eta}_n^{1-\frac{1}{m}}+\sbr[2]{\frac{\log(dn)}{n}}^{1-\frac{1}{m}}}\sqrt{n\log(d)}}.
	\end{align*}
	This implies the claimed bound on~$\rho_{n,W}$ in~\eqref{eq:Gaussapprox} by using that~$4/n$ in~\eqref{eq:1overn} is dominated by, e.g., $\sbr[1]{\frac{\log(dn)}{n}}^{1-\frac{1}{m}}\sqrt{n\log(d)}$ (if necessary adjust~$C_7$) as well as~$\P^*\del[1]{Z\leq t}=\gamma_{\Sigma}(H)$. 
	
	If, on the other hand,~\eqref{eq:CaseGood} is not satisfied, then 
	\begin{align*}
		C_1\eps_n^{1-\frac{2}{m}}
		\geq
		\frac{(b_1^2\wedge 1)}{2}	\quad\iff \quad 
		\eps_n
		\geq
		\sbr[3]{\frac{(b_1^2\wedge 1)}{2C_1}}^{\frac{m}{m-2}}=:K(b_1,b_2,\lambda_2,m),	
	\end{align*}
	where~$K=K(b_1,b_2,\lambda_2,m)$ does not depend on~$n$. Thus, by the definition of~$\eps_n$, 
	\begin{align*}
		\overline{\eta}_n\geq\frac{K}{2\lambda_1}\qquad\text{or}\qquad \frac{\log(dn)}{n}\geq\frac{K}{2\lambda_2}, 
	\end{align*}
	In either case, since~$\rho_{n,W}\leq 1$ the bound in~\eqref{eq:Gaussapprox} remains valid by adjusting the constant~$C$ there, if necessary.
	
	Finally, suppose that~$\sqrt{n\log(d)}\overline{\eta}_n^{1-\frac{1}{m}}\to 0$ and~$\log(d)/n^{\frac{m-2}{3m-2}}\to 0$. If~$m\in(2,4]$, it follows by Lemma~\ref{lem:covcompare} that~$\rho_{n,W}\to 0$. If~$m\in(4,\infty)$ and additionally $\log(d)/\del[1]{n^{\frac{m-2}{3m-2}}/\log^{\frac{m}{3m-2}}(n)}\to 0$ (implying that~$\log(d)/n^{\frac{m-2}{3m-2}}\to 0$), then Lemma~\ref{lem:covcompare} again yields that~$\rho_{n,W}\to 0$. 
	
\section{Proof of Theorem \ref{thm:CorrTargetSharp}}
As was the case for Theorem~\ref{thm:SufficientWins}, we establish a more general version of Theorem~\ref{thm:CorrTargetSharp} in the setting allowing coordinatewise adversarial contamination detailed in Section~\ref{sec:ContamSetting}. We maintain~$m>2$,~$d \geq 2$ and~$n > 3$. Although the coordinatewise contamination scheme in~\eqref{eq:contamfrac} is more general than the one in~\cite{kock2025high}, their estimator of~$\Sigma$, whose diagonal elements are used to standardize the corresponding entries of~$S_{n,W}$, retains the exact same precision guarantees in this more general setting. To prove this, we briefly recall the definition of their estimator.  Set
\begin{equation}\label{eq:epsprime}
	\eps_n'= \lambda_1'\cdot\overline{\eta}_n+\lambda_2'\cdot\frac{\log(d^2n)}{n},\qquad \lambda_1'\in(1,\infty)\text{ and }\lambda_2'\in (0,\infty).
\end{equation}
\emph{Like for~$\lambda_1$ and~$\lambda_2$, we take~$\lambda_1'$ and~$\lambda_2'$ as fixed throughout}.

Writing~$\hat{a}_j=\tilde{X}_{\lceil \eps_n'n\rceil,j}^*$,~$\hat{b}_j=\tilde{X}_{\lfloor (1-\eps_n')n\rfloor+1,j}^*$, and~$\tilde{\mu}_{n,j}=n^{-1}\sum_{i=1}^n\phi_{\hat{a}_j,\hat{b}_j}(\tilde{X}_{i,j})$ for~$j=1,\hdots,d$, define~$\tilde{\Sigma}_n$ as the matrix with entries 
\begin{align}\label{eq:tildeSigma}
	\tilde{\Sigma}_{n,j,k}
	=
	\frac{1}{n}\sum_{i=1}^n\sbr[1]{\phi_{\hat{a}_j,\hat{b}_j}(\tilde{X}_{i,j})-\tilde{\mu}_{n,j}}\sbr[1]{\phi_{\hat{a}_k,\hat{b}_k}(\tilde{X}_{i,k})-\tilde{\mu}_{n,k}},\qquad 1\leq j,k\leq d.
\end{align}
The estimator~$\tilde{\Sigma}_n$, which is a winsorized analog of the sample covariance matrix, is positive semi-definite and symmetric by virtue of being a Gram matrix. We note that the proof of Theorem~\ref{thm:covestimGRAM} below dispenses with the finite-sample condition on~$\eps_n'$ imposed in (9) in~\cite{kock2025high} and remains valid under the more demanding contamination in~\eqref{eq:contamfrac}. For~$A_+'=1-\lambda_1'^{-1}\mathds{1}\cbr[0]{\overline{\eta}_n>0}$, let
		\begin{equation}\label{eq:c1'}
		c_{1,n}':=-A_+'W_0\del[1]{-e^{-(\frac{\log(d^2n)}{\eps_n' n}+A_+')/A_+'}}\in(0,A_+'),
		\end{equation}
with~$W_0$ denoting the principal branch of Lambert's~$W$ function.
\begin{theorem}\label{thm:covestimGRAM}
	Let~$P\in\mc{P}_d(b_1,b_2,m)$ and~\eqref{eq:contamfrac} be satisfied. If~$\eps_n'\in(0,1/2)$, with~$\eps_n'$ as in~\eqref{eq:epsprime}, then, for a constant~$C$ depending only on~$b_2,\lambda_1',\lambda_2'$ and~$m$, we have
	\begin{align}\label{eq:covestimGram}
		\P\del[4]{\max_{1\leq j,k\leq d}\envert[1]{\tilde{\Sigma}_{n,j,k}-\Sigma_{j,k}}> C\sbr[3]{\overline{\eta}_n^{1-\frac{2}{m}}+\del[2]{\frac{\log(dn)}{n}}^{1-\frac{1}{(m/2)\wedge 2}}}}
		\leq \frac{24}{n}.
	\end{align}	
\end{theorem}

\begin{proof}[Proof of Theorem~\ref{thm:covestimGRAM}]
		Fix~$1\leq j,k\leq d$ and note that
		\begin{align}
			\tilde{\Sigma}_{n,j,k}
			&=
			\frac{1}{n}\sum_{i=1}^n\sbr[1]{\phi_{\hat{a}_j,\hat{b}_j}(\tilde{X}_{i,j})-\mu_j-(\tilde{\mu}_{n,j}-\mu_j)}\sbr[1]{\phi_{\hat{a}_k,\hat{b}_k}(\tilde{X}_{i,k})-\mu_k-(\tilde{\mu}_{n,k}-\mu_k)}\notag\\
			&=
			\frac{1}{n}\sum_{i=1}^n\sbr[1]{\phi_{\hat{a}_j,\hat{b}_j}(\tilde{X}_{i,j})-\mu_j}\sbr[1]{\phi_{\hat{a}_k,\hat{b}_k}(\tilde{X}_{i,k})-\mu_k}-(\tilde{\mu}_{n,j}-\mu_j)(\tilde{\mu}_{n,k}-\mu_k)\notag\\
			&=
			\frac{1}{n}\sum_{i=1}^n\phi_{\check{a}_j,\check{b}_j}(\tilde{Y}_{i,j})\phi_{\check{a}_k,\check{b}_k}(\tilde{Y}_{i,k})-(\tilde{\mu}_{n,j}-\mu_j)(\tilde{\mu}_{n,k}-\mu_k)\label{eq:covdecomp},
		\end{align}
		where~$Y_{i,j}=X_{i,j}-\mu_j$,~$\tilde{Y}_{i,j}=\tilde{X}_{i,j}-\mu_j$,~$\check{a}_j=\tilde{Y}_{\lceil \eps_n'n\rceil,j}^*=\tilde{X}_{\lceil \eps_n'n\rceil,j}^*-\mu_j=\hat{a}_j-\mu_j$, and~$\check{b}_j=\tilde{Y}_{\lfloor (1-\eps_n')n\rfloor+1,j}^*=\tilde{X}_{\lfloor (1-\eps_n')n\rfloor+1,j}^*-\mu_j=\hat{b}_j-\mu_j$. Next, introduce the auxiliary
		\begin{equation*}
			\eps_n''= \lambda_1'\cdot(2\overline{\eta}_n)+\lambda_2'\cdot\frac{\log(d^2n)}{n}\in[\eps_n',2\eps_n'],\qquad \lambda_1'\in(1,\infty)\text{ and }\lambda_2'\in (0,\infty).	
		\end{equation*}
		We now separately consider the case of
		\begin{equation}\label{eq:epscond''}
	2\eps_n'' +\frac{\log(d^2n)}{n}+\sqrt{\del[2]{\frac{\log(d^2n)}{n}}^2+4\frac{\log(d^2n)}{n}\eps_n''}<1.	
\end{equation}
 being satisfied or not.
 
 \textbf{Case:~\eqref{eq:epscond''} holds:} Since~$\eps_n'\leq \eps_n''$, \eqref{eq:epscond''} is also satisfied with~$\eps_n''$ replaced by~$\eps_n'$. Thus, by Theorem 2.1 in~\cite{Wins1} (the assumptions of which are satisfied for~$P\in\mc{P}_d(b_1,b_2,m)$ and under~\eqref{eq:contamfrac}) applied with~$\delta=6/(d^2n)$,~$\eta=\overline{\eta}_n$,~$\lambda_i=\lambda_i'$ for~$i=1,2$ (so that~$\eps$ there is~$\eps_n'$ and (5) there is satisfied), there exists a constant~$C_1$ depending only on~$b_2,\lambda_1',\lambda_2'$, and~$m$, such that with probability at least~$1-12/(d^2n)$ 
		\begin{align}
			|\tilde{\mu}_{n,j}-\mu_j||\tilde{\mu}_{n,k}-\mu_k|
			&\leq 
			C_1\del[3]{\overline{\eta}_n^{1-\frac{1}{m}}+\sbr[2]{\frac{\log(dn)}{n}}^{1-\frac{1}{m\wedge 2}}}^2\notag\\
			&\leq
			C_1\del[3]{\overline{\eta}_n^{1-\frac{2}{m}}+\sbr[2]{\frac{\log(dn)}{n}}^{1-\frac{1}{(m/2)\wedge 2}}},\label{eq:meanestimerror}
		\end{align}
		the second estimate using that~$\overline{\eta}_n$ and~$\frac{\log(dn)}{n}$ are both strictly smaller than 1 by~\eqref{eq:epscond''} (and potentially adjusting~$C_1$). Next, write~$U_{i,j,k}=Y_{i,j}Y_{i,k}$,~$\tilde{U}_{i,j,k}=\tilde{Y}_{i,j}\tilde{Y}_{i,k}$,~$\check{\alpha}_{j,k}=\tilde{U}_{\lceil \eps_n''n\rceil,j,k}^*$, and~$\check{\beta}_{j,k}=\tilde{U}_{\lfloor (1-\eps_n'')n\rfloor+1,j,k}^*$. Note that~$\E U_{1,j,k}=\Sigma_{j,k}$,~$\E\envert[0]{U_{1,j,k}-\E U_{1,j,k}}^{m/2}
		\leq 2^{m/2}\E|U_{1,j,k}|^{m/2}\leq 
		2^{m/2}b_2^m$, and by~\eqref{eq:contamfrac}
		\begin{equation*}
		\envert[1]{\cbr[1]{i\in\cbr[0]{1,\hdots,n}: \tilde{U}_{i,j,k}\neq U_{i,j,k}}}\leq 2\overline{\eta}_nn,\qquad\text{for all }j,k\in\cbr[0]{1,\hdots,d}.	
		\end{equation*}
		Therefore, by Theorem 2.1 in~\cite{Wins1} applied with~$\delta=6/(d^2n)$,~$\eta=2\overline{\eta}_n$,~$m$ \emph{there} being~$m/2>1$,~$\lambda_i=\lambda_i'$ for~$i=1,2$ (so that~$\eps$ there is~$\eps_n''$ and (5) there is satisfied), there exists a constant~$C_2$ depending only on~$b_2,\lambda_1',\lambda_2'$, and~$m$ such that with probability at least~$1-6/(d^2n)$ 
		\begin{align*}
			\envert[2]{\frac{1}{n}\sum_{i=1}^n\phi_{\check{\alpha}_{j,k},\check{\beta}_{j,k}}(\tilde{U}_{i,j,k})-\Sigma_{j,k}}
			\leq C_2\del[3]{\overline{\eta}_n^{1-\frac{2}{m}}+\sbr[2]{\frac{\log(dn)}{n}}^{1-\frac{1}{(m/2)\wedge 2}}}.
		\end{align*}
		Thus, defining
		\begin{align*}
			\mathfrak{S}_{n,j,k}
			:=
			\frac{1}{n}\sum_{i=1}^n\phi_{\check{a}_j,\check{b}_j}(\tilde{Y}_{i,j})\phi_{\check{a}_k,\check{b}_k}(\tilde{Y}_{i,k})-\frac{1}{n}\sum_{i=1}^n\phi_{\check{\alpha}_{j,k},\check{\beta}_{j,k}}(\tilde{U}_{i,j,k}),
		\end{align*}
		in order to prove the theorem in the present case, it suffices (by the union bound over~$d^2$) to show that there exists a constant~$C_3$ depending only on~$b_2,\lambda_1',\lambda_2'$, and~$m$ such that with probability at least~$1-6/(d^2n)$ 
		\begin{equation}\label{eq:Toshow}
			\envert[1]{\mathfrak{S}_{n,j,k}}
			\leq C_3\del[3]{\overline{\eta}_n^{1-\frac{2}{m}}+\sbr[2]{\frac{\log(dn)}{n}}^{1-\frac{1}{(m/2)\wedge 2}}}.	
		\end{equation}
		To establish~\eqref{eq:Toshow}, note that with
		\begin{align*}
			A_{n,j,k}:=\cbr[2]{i\in\cbr[1]{1,\hdots,n}: \check{\alpha}_{j,k}\leq \tilde{U}_{i,j,k}\leq \check{\beta}_{j,k} \text{ and }\check{a}_j\leq \tilde{Y}_{i,j}\leq\check{b}_j\text{ and }\check{a}_k\leq \tilde{Y}_{i,k}\leq\check{b}_k},
		\end{align*}
		one has that
		\begin{align*}
			\phi_{\check{\alpha}_{j,k},\check{\beta}_{j,k}}\del[0]{\tilde{U}_{i,j,k}}-\phi_{\check{a}_j,\check{b}_j}(\tilde{Y}_{i,j})\phi_{\check{a}_k,\check{b}_k}(\tilde{Y}_{i,k})
			=
			\tilde{U}_{i,j,k}-\tilde{Y}_{i,j}\tilde{Y}_{i,k}
			=0\qquad\text{for }i\in A_{n,j,k}. 
		\end{align*}
		Hence,
		\begin{align*}
			\mathfrak{S}_{n,j,k}
			&=
			\frac{1}{n}\sum_{i\in A_{n,j,k}^c}\sbr[1]{\phi_{\check{a}_j,\check{b}_j}(\tilde{Y}_{i,j})\phi_{\check{a}_k,\check{b}_k}(\tilde{Y}_{i,k})-\phi_{\check{\alpha}_{j,k},\check{\beta}_{j,k}}\del[0]{\tilde{U}_{i,j,k}}},
		\end{align*}	
		and for every~$i\in A_{n,j,k}^c$
		\begin{align}\label{eq:tailterms}
			\envert[2]{\phi_{\check{\alpha}_{j,k},\check{\beta}_{j,k}}\del[0]{\tilde{U}_{i,j,k}}-\phi_{\check{a}_j,\check{b}_j}(\tilde{Y}_{i,j})\phi_{\check{a}_k,\check{b}_k}(\tilde{Y}_{i,k})}
			\leq \del[1]{|\check{\alpha}_{j,k}| \vee |\check{\beta}_{j,k}|}	+\del[1]{|\check{a}_j| \vee |\check{b}_j|}\del[1]{|\check{a}_k| \vee |\check{b}_k|}.
		\end{align}
		To bound the right-hand side of~\eqref{eq:tailterms}, note that by Lemma~B.5 and Remark B.1 in~\cite{Wins1} applied with~$\delta=6/(d^2n)$,~$\eta=\overline{\eta}_n$,~$\lambda_i=\lambda_i'$ for~$i=1,2$ so that~$\eps$ there is~$\eps_n'\geq \lambda_1'\overline{\eta}_n$, it holds that
			\begin{align}\label{eq:boundorderstats}
			\del[1]{|\check{a}_j| \vee |\check{b}_j|}\del[1]{|\check{a}_k| \vee |\check{b}_k|}
			\leq &
			\del[1]{|Q_{c_{1,n}'\eps_n'}(Y_{1,j})|\vee  |Q_{1-c_{1,n}'\eps_n'}(Y_{1,j})|}\notag\\
			\cdot & \del[1]{|Q_{c_{1,n}'\eps_n'}(Y_{1,k})|\vee  |Q_{1-c_{1,n}'\eps_n'}(Y_{1,k})|}.
		\end{align}
		with probability at least~$1-\frac{4}{d^2n}$, for~$c_{1,n}'$ as defined in~\eqref{eq:c1'}. Here we also used that $Q_{c_{1,n}'\eps_n'}(Y_{1,j})\leq \check{a}_j\leq \check{b}_j\leq Q_{1-c_{1,n}'\eps_n'}(Y_{1,j})$ and~$Q_{c_{1,n}'\eps_n'}(Y_{1,k})\leq \check{a}_k\leq \check{b}_k\leq Q_{1-c_{1,n}'\eps_n'}(Y_{1,k})$ because~$\eps_n'\in(0,1/2)$ by assumption. 
		
	Use Lemma~B.5 and Remark B.1 in~\cite{Wins1} once more, but with~$\delta=6/(d^2n)$,~$\eta=2\overline{\eta}_n$,~$\lambda_i=\lambda_i'$ for~$i=1,2$ so that~$\eps$ there is~$\eps_n''\geq \lambda_1'\cdot(2\overline{\eta}_n)$. It holds with probability at least~$1-\frac{2}{d^2n}$ that
		\begin{equation*}
			|\check{\alpha}_{j,k}| \vee |\check{\beta}_{j,k}|
			\leq
			|Q_{c_{1,n}''\eps_n''}(U_{1,j,k})|\vee |Q_{1-c_{1,n}''\eps_n''}(U_{1,j,k})|	.
		\end{equation*}
	for
	\begin{equation*}
		c_{1,n}'':=-A_+'W_0\del[1]{-e^{-(\frac{\log(d^2n)}{\eps_n'' n}+A_+')/A_+'}}\in(0,A_+').
	\end{equation*}
		Thus, with probability at least~$1-\frac{6}{d^{2}n}$ the right-hand side of~\eqref{eq:tailterms} is no greater than
		\begin{align*}
			|Q_{c_{1,n}''\eps_n''}(U_{1,j,k})|\vee |Q_{1-c_{1,n}''\eps_n''}(U_{1,j,k})|
			+&\del[1]{|Q_{c_{1,n}'\eps_n'}(Y_{1,j})|\vee  |Q_{1-c_{1,n}'\eps_n'}(Y_{1,j})|}\\
			\cdot & \del[1]{|Q_{c_{1,n}'\eps_n'}(Y_{1,k})|\vee  |Q_{1-c_{1,n}'\eps_n'}(Y_{1,k})|}.
		\end{align*}
		Next, recall that~$\E U_{1,j,k}=\Sigma_{j,k}$,~$\E\envert[0]{U_{1,j,k}-\E U_{1,j,k}}^{m/2}
		\leq 2^{m/2}\E|U_{1,j,k}|^{m/2}\leq 
		2^{m/2}b_2^m$, $\E Y_{1,j}=0$, and~$\E |Y_{1,j}|^m\leq b_2^m$ for all~$1\leq j, k\leq d$. Hence, Lemma~C.1 in~\cite{Wins1} implies that the previous display is bounded from above by 
		\begin{align*}
			\del[3]{|\Sigma_{j,k}|+\frac{2b_2^2}{(c_{1,n}''\eps_n'')^{2/m}}}+\frac{b_2^2}{(c_{1,n}'\eps_n')^{2/m}}
			\leq
			\frac{4b_2^2}{(c_{1,n}'\eps_n')^{2/m}},
		\end{align*}
		the inequality following from~$|\Sigma_{j,k}|\leq b_2^2$,~$\eps_n''\geq \eps_n'$ and~$c_{1,n}''\geq c_{1,n}'$ as~$W_0$ is strictly increasing. Observe that~$|A_{n,j,k}^c|\leq (2\eps_n'' +4\eps_n')n\leq 8\eps_n'n$ because~$\eps_n''\leq 2\eps_n'$. Thus, with probability at least~$1-\frac{6}{d^2n}$ it holds that
		\begin{align*}
			\envert[1]{\mathfrak{S}_{n,j,k}}
			\leq
			8\eps_n'\cdot \frac{4b_2^2}{(c_{1,n}'\eps_n')^{2/m}}
			\leq 
			\frac{32b_2^2}{(c_{1,n}')^{2/m}}\cdot{\eps_n'}^{1-\frac{2}{m}}
			\leq 
			C_3\del[2]{\overline{\eta}_n^{1-\frac{2}{m}}+\sbr[2]{\frac{\log(d^2n)}{n}}^{1-\frac{2}{m}}},
		\end{align*}
		where we inserted~$\eps_n'$ from~\eqref{eq:epsprime}, used subadditivity of~$z\mapsto z^{1-\frac{2}{m}}$, the lower bound on~$c_{1,n}'$ from Lemma~B.3 in~\cite{Wins1}, and~$C_3$ is the desired constant depending only on~$b_2,\lambda_1',\lambda_2'$, and~$m$. Hence,~\eqref{eq:Toshow} (and thus also~\eqref{eq:covestimGram}) follows in case~\eqref{eq:epscond''} holds because~$1-\frac{2}{m}\geq 1-\frac{1}{(m/2)\wedge 2}$ and~$\frac{\log(d^2n)}{n}<1$ by~\eqref{eq:epscond''}.		
		
		\textbf{Case~\eqref{eq:epscond''} does not hold:} 
		Because~$\log(d^2n)/n\leq (\lambda_2')^{-1}\eps_n''$ it follows from (the negation of) \eqref{eq:epscond''} that~$\eps_n''\geq C_4$, for~$C_4$ a constant depending on~$\lambda_2'$ only. Thus, also~$\eps_n'\geq \eps_n''/2\geq C_4/2$. Next, note that, cf.~the argument following~\eqref{eq:covdecomp},
		\begin{equation*}
			\tilde{\mu}_{n,j}-\mu_j
			=
			\frac{1}{n}\sum_{i=1}^n(\phi_{\hat{a}_j,\hat{b}_j}(\tilde{X}_{i,j})-\mu_j)
			=
			\frac{1}{n}\sum_{i=1}^n\phi_{\check{a}_j,\check{b}_j}(\tilde{Y}_{i,j})\in\sbr[1]{\check{a}_j,\check{b}_j},\qquad j\in\cbr[0]{1,\hdots,d}.
		\end{equation*}
		Thus, in combination with~\eqref{eq:covdecomp} and~\eqref{eq:boundorderstats} (the latter only using that~$\eps_n'\in(0,1/2)$), it holds for every pair~$(j,k)$ with probability at least~$1-\frac{4}{d^2n}$ that
		\begin{align*}
		|\tilde{\Sigma}_{n,j,k}|	
		&\leq
		\frac{1}{n}\sum_{i=1}^n|\phi_{\check{a}_j,\check{b}_j}(\tilde{Y}_{i,j})\phi_{\check{a}_k,\check{b}_k}(\tilde{Y}_{i,k})|+|\tilde{\mu}_{n,j}-\mu_j||\tilde{\mu}_{n,k}-\mu_k|\\
		&\leq
		2\cdot\del[1]{|Q_{c_{1,n}'\eps_n'}(Y_{1,j})|\vee  |Q_{1-c_{1,n}'\eps_n'}(Y_{1,j})|}
			\cdot \del[1]{|Q_{c_{1,n}'\eps_n'}(Y_{1,k})|\vee  |Q_{1-c_{1,n}'\eps_n'}(Y_{1,k})|}.
		\end{align*}
By Lemma~C.1 in~\cite{Wins1}
\begin{equation*}
	\del[1]{|Q_{c_{1,n}'\eps_n'}(Y_{1,j})|\vee  |Q_{1-c_{1,n}'\eps_n'}(Y_{1,j})|}
			\cdot \del[1]{|Q_{c_{1,n}'\eps_n'}(Y_{1,k})|\vee  |Q_{1-c_{1,n}'\eps_n'}(Y_{1,k})|}
			\leq
			\frac{b_2^2}{(c_{1,n}'\eps_n')^{2/m}}
			\leq C_5,
\end{equation*}
the last inequality using the lower bound on~$c_{1,n}'$ from Lemma~B.3 in~\cite{Wins1} and that~$\eps_n'\geq C_4/2$; here~$C_5$ is a constant depending only on~$b_2,\lambda_1',\lambda_2'$, and~$m$.  Thus, using also that~$\max_{1\leq j,k\leq d}|\Sigma_{j,k}|\leq b_2^2$, it holds with probability at least~$1-\frac{4}{d^2n}$
\begin{equation*}
|\tilde{\Sigma}_{n,j,k}-\Sigma_{j,k}|
\leq
|\tilde{\Sigma}_{n,j,k}|+|\Sigma_{j,k}|
\leq 
2C_5+b_2^2.	
\end{equation*} 
 Finally, because~$\eps_n'\geq C_4/2$, it must be the case that 
		\begin{equation*}
			\overline{\eta}_n\geq \frac{C_4}{4\lambda_1'}\qquad\text{or}\qquad \frac{2\log(dn)}{n}
			=
			\frac{\log(d^2n^2)}{n}
			\geq
			\frac{\log(d^2n)}{n}
			\geq 
			\frac{C_4}{4\lambda_2'},
		\end{equation*}
Therefore, the far right-hand side of the penultimate display is not greater than
 \begin{equation*}
 	C_6\del[3]{\overline{\eta}_n^{1-\frac{2}{m}}+\sbr[2]{\frac{\log(dn)}{n}}^{1-\frac{1}{(m/2)\wedge 2}}}
 \end{equation*}
 for~$C_6$ a constant depending only on~$b_2,\lambda_1',\lambda_2'$, and~$m$. The conclusion now follows from the union bound over~$1\leq j,k\leq d$.
\end{proof}

Theorem~\ref{thm:CorrTargetSharp} follows from Theorem~\ref{thm:HDGauss_studentized} below, which is valid under a modest amount of adversarial contamination. Let~$\tilde{\sigma}_{n,j}=(\tilde{\Sigma}_{n,j,j})^{1/2}$ be the square root of the diagonal elements of the estimator~$\tilde{\Sigma}_n$ of~$\Sigma$, cf.~\eqref{eq:tildeSigma}. Write~$D=D(P)=\text{diag}(\sigma_{2,1}(P),\hdots,\sigma_{2,d}(P))$ for the square root of the diagonal elements of~$\Sigma(P)$ so that~$\Sigma_0=D^{-1}\Sigma D^{-1}$. Define~$S_{n,W,S}=S_{n,W,S}(P)\in\R^d$ via

\begin{equation}\label{eq:S_nWS}
	S_{n,W,S,j}=\frac{1}{\sqrt{n}\tilde{\sigma}_{n,j}}\sum_{i=1}^n\sbr[1]{\phi_{\hat\alpha_j,\hat\beta_j}(\tilde{X}_{i,j})-\mu_j},\qquad j=1,\hdots,d,
\end{equation}
with~$S_{n,W,S,j}=0$ if~$\tilde{\sigma}_{n,j}=0$.
\begin{theorem}\label{thm:HDGauss_studentized}
	Let~$P\in\mc{P}_d(b_1,b_2,m)$ and~\eqref{eq:contamfrac} be satisfied.  If~$\eps_n,\eps_n'\in(0,1/2)$, with~$\eps_n$ as in~\eqref{eq:epsfam} and~$\eps_n'$ as in~\eqref{eq:epsprime}, then 	
	\begin{equation}\label{eq:HDGaussstudentized}
		\rho_{n,W,S}
		:=
		\sup_{H\in\mc{H}_d}\envert[2]{\P\del[1]{S_{n,W,S}\in H}-\gamma_{\Sigma_0}(H)}
		\leq
		C\mathfrak{A}_n,
	\end{equation}
	where~$\mathfrak{A}_n$ is the upper bound on~$\rho_{n,W}$ in~\eqref{eq:Gaussapprox} of Theorem~\ref{thm:HDGaussSharp} and~$C$ is a constant depending only on~$b_1,b_2,\lambda_1,\lambda_2,\lambda_1',\lambda_2'$ and~$m$. Suppose additionally that~$\sqrt{n\log(d)}\overline{\eta}_n^{1-\frac{1}{m}}\to 0$. Then,
	\begin{equation*}
		\rho_{n,W,S}\to 0\ 
		\begin{cases}
			\text{for}\quad m\in(2,4]\quad &\text{if}\quad \log(d)/n^{\frac{m-2}{3m-2}}\to 0;\\  
			\text{for}\quad m\in(4,\infty)\quad &\text{if}\quad \log(d)/\del[1]{n^{\frac{m-2}{3m-2}}/\log^{\frac{m}{3m-2}}(n)}\to 0.
		\end{cases}
	\end{equation*}
\end{theorem}

\begin{proof}[Proof of Theorem~\ref{thm:CorrTargetSharp}]
	Use~$\lambda_1=\lambda_1'=2$,~$\lambda_2=\lambda$,~$\lambda_2'=\lambda'$ and~$\overline{\eta}_n=0$ for all~$n\in\N$ in Theorem~\ref{thm:HDGauss_studentized}, noting that~$C$ there does not depend on~$P\in\mc{P}_d(b_1,b_2,m)$.
\end{proof}

We will use the following notation in the proof of Theorem~\ref{thm:HDGauss_studentized}.  For~$x\in\R^d$, let~$\enVert[0]{x}_\infty=\max_{j=1,\hdots,d}|x_j|$. For any $A\subseteq \R^d$ and~$\zeta>0$, let $$A^{\zeta,\infty}=\cbr[1]{x\in\R^d:\inf_{y\in A}\|x-y\|_\infty\leq \zeta}.$$ Furthermore, $$A^{-\zeta,\infty}=\cbr[1]{x\in\R^d:\mc{B}_\infty(x,\zeta) \subseteq A} \quad \text{ where } \quad \mc{B}_{\infty}(x,\zeta)=\cbr[1]{y\in\R^d:\|y-x\|_{\infty} \leq \zeta}.$$
	
	\begin{proof}[Proof of Theorem~\ref{thm:HDGauss_studentized}]
		The proof follows the same strategy as that of Theorem 4.1 in~\cite{kock2025high}, the main difference being that it uses Theorem~\ref{thm:HDGaussSharp} of Section~\ref{sec:ContamSetting} to conclude a Gaussian approximation to the distribution of the auxiliary~$V_n\in\R^d$ defined in~\eqref{eq:V_n} below instead of the weaker Theorem 2.1 of that reference. We provide the details for completeness. We first establish~\eqref{eq:HDGaussstudentized}. Assume that
		\begin{align}\label{eq:CaseNormalize}
		C_1\sbr[3]{\overline{\eta}_n^{1-\frac{2}{m}}+\del[2]{\frac{\log(dn)}{n}}^{1-\frac{1}{(m/2)\wedge 2}}}\leq \frac{b_1^2}{2},
		\end{align}
		where~$C_1$ is the constant from Theorem~\ref{thm:covestimGRAM} depending only on~$b_2,\lambda_1',\lambda_2'$, and~$m$. If~\eqref{eq:CaseNormalize} is violated, we conclude~\eqref{eq:HDGaussstudentized} with similar arguments as those employed at the end of the proof of Theorem~\ref{thm:HDGaussSharp}. Thus, suppose in what follows that~\eqref{eq:CaseNormalize} is satisfied.

		Let~$V_n\in\R^d$ have entries
		\begin{equation}\label{eq:V_n}
			V_{n,j}=\frac{1}{\sqrt{n}\sigma_{2,j}}\sum_{i=1}^n\sbr[1]{\phi_{\hat\alpha_j,\hat\beta_j}(\tilde{X}_{i,j})-\mu_j},\qquad j=1,\hdots,d.	
		\end{equation}
		By Theorem~\ref{thm:covestimGRAM},~$\min_{j=1,\hdots,d}\sigma_{2,j}\geq b_1$, and the mean-value theorem, there exists an event~$E_n$ of probability at least~$1-24/n$ such that on~$E_n$ it holds that~$\tilde{\sigma}_{n,j} \geq b_1/\sqrt{2} > 0$ for every~$j = 1, \hdots, d$ (we used~\eqref{eq:CaseNormalize}) and
		\begin{equation}
			A_n
			:=
			\enVert[1]{S_{n,W,S}-V_n}_\infty
			\leq 
			\max_{j=1,\hdots,d}\envert[3]{\frac{1}{\tilde{\sigma}_{n,j}}-\frac{1}{\sigma_{2,j}}}
			\cdot\max_{j=1,\hdots,d}\envert[3]{\frac{1}{\sqrt{n}}\sum_{i=1}^n\sbr[1]{\phi_{\hat\alpha_j,\hat\beta_j}(\tilde{X}_{i,j})-\mu_j}}\label{eq:An1}.
		\end{equation}
		In addition, there exists a constant~$C_2$ depending only on $b_1,b_2,\lambda_1',\lambda_2'$, and~$m$ so that on~$E_n$
		\begin{equation}
			\max_{j=1,\hdots,d}\envert[3]{\frac{1}{\tilde{\sigma}_{n,j}}-\frac{1}{\sigma_{2,j}}}
			=
			\max_{j=1,\hdots,d}\frac{|\tilde{\sigma}_{n,j}-\sigma_{2,j}|}{\tilde{\sigma}_{n,j}\sigma_{2,j}}
			\leq 
			C_2\del[3]{\overline{\eta}_n^{1-\frac{2}{m}}+\del[2]{\frac{\log(dn)}{n}}^{1-\frac{1}{(m/2)\wedge 2}}}\label{eq:An2}.
		\end{equation}
		Furthermore, as~$\sigma_{2,j}\leq b_2$, the union bound and~$\P^*\del[0]{|\bm{z}|> t}\leq 2\exp(-t^2/2)$ for~$t\geq 0$ and~$\bm{z} \sim \mathsf{N}_1(0,1)$ yields for~$Z \sim \mathsf{N}_d(0, \Sigma)$ that [here~$\bm{z}$ and~$Z$ are defined on an auxiliary probability space with probability measure~$P^*$]
		\begin{equation*}
	 P^*\del{\max_{j=1,\hdots,d}|Z_j|> b_2\sqrt{2\log(2dn)}}
			\leq
			\frac{1}{n}.
		\end{equation*}
		Thus, by Theorem~\ref{thm:HDGaussSharp}	it follows that there exists a constant~$K_1$ depending only on~$b_1,b_2,\lambda_1,\lambda_2,$ and~$m$ such that
		\begin{align}
			\P\del[3]{\max_{j=1,\hdots,d}\envert[3]{\frac{1}{\sqrt{n}}\sum_{i=1}^n\sbr[1]{\phi_{\hat\alpha_j,\hat\beta_j}(\tilde{X}_{i,j})-\mu_j}}> b_2\sqrt{2\log(2dn)}}
			\leq 
			\mathfrak{A}_n+\frac{1}{n}
			\leq
			K_1\mathfrak{A}_n\label{eq:An3}.
		\end{align} 
Hence, by~\eqref{eq:An1}--\eqref{eq:An3} there exists a constant~$C_3$ depending only on~$b_1, b_2,\lambda_1',\lambda_2'$, and~$m$ such that
		\begin{equation}
			\P\del[3]{A_n\leq C_3\sqrt{\log(dn)}\sbr[2]{\overline{\eta}_n^{1-\frac{2}{m}}+\del[2]{\frac{\log(dn)}{n}}^{1-\frac{1}{(m/2)\wedge 2}}}}
			\geq 
			1-\frac{24}{n}-K_1\mathfrak{A}_n
			\geq
			1-K\mathfrak{A}_n\label{eq:An4},
		\end{equation}
		where~$K$ depends only on~$b_1,b_2,\lambda_1,\lambda_2$, and~$m$. Thus, writing
		\begin{equation*}
			\overline{A}_n:=C_3\sqrt{\log(dn)}\sbr[3]{\overline{\eta}_n^{1-\frac{2}{m}}+\del[2]{\frac{\log(dn)}{n}}^{1-\frac{1}{(m/2)\wedge 2}}},
		\end{equation*}
		it holds for all~$H\in\mc{H}_d$ that
		\begin{equation}\label{eq:setdecomp1}
			\cbr[1]{S_{n,W,S}\in H}
			\subseteq 
			\cbr[1]{V_n\in H^{\overline{A}_n,\infty}}
			\bigcup \cbr[1]{A_n>\overline{A}_n}.
		\end{equation} 
		Note that~$DV_n=S_{n,W}$ and~$DH\in\mc{H}_d$ for every~$H\in\mc{H}_d$. Therefore, since~$H^{\overline{A}_n,\infty}\in\mc{H}_d$ and recalling that~$\Sigma=D\Sigma_0D$, Theorem~\ref{thm:HDGaussSharp} yields
			\begin{equation*}
			\envert[2]{\P\del[1]{V_n\in H^{\overline{A}_n,\infty}}-\gamma_{\Sigma_0}\del[1]{H^{\overline{A}_n,\infty}}}
			=
				\envert[2]{\P\del[1]{S_{n,W}\in DH^{\overline{A}_n,\infty}}-\gamma_{\Sigma}\del[1]{DH^{\overline{A}_n,\infty}}}
			\leq
			\mathfrak{A}_n.
		\end{equation*}
		Furthermore, by Lemma~E.1 in~\cite{kock2025high} (applied with~$\delta_1=-\overline{A}_n$ and~$\delta_2=\overline{A}_n$), noting that~$\Sigma_{0,j,j}=1$ for all~$j=1,\hdots,d$, 
		\begin{equation*}
			0
			\leq 
			\gamma_{\Sigma_0}\del[1]{H^{\overline{A}_n,\infty}}
			-
			\gamma_{\Sigma_0}\del[0]{H}
			\leq 
			\overline{A}_n\del[1]{\sqrt{2\log(d)}+4}
			\leq
			7\overline{A}_n\sqrt{\log(d)},
		\end{equation*}
		the last inequality following from~$d\geq 2$. Using this in~\eqref{eq:setdecomp1} along with~\eqref{eq:An4} yields that
		\begin{equation*}
			\P\del[1]{S_{n,W,S}\in H}-\gamma_{\Sigma_0}\del[0]{H}
			\leq
			(1+K)\mathfrak{A}_n+7\overline{A}_n\sqrt{\log(d)}. 	
		\end{equation*}
		Since 
		\begin{equation*}
			\cbr[1]{V_n\in H^{-\overline{A}_n,\infty}}
			\subseteq 
			\cbr[1]{S_{n,W,S}\in H}	
			\bigcup \cbr[1]{A_n>\overline{A}_n},
		\end{equation*}
		the same arguments as those following~\eqref{eq:setdecomp1} lead to
		\begin{equation*}
			-(1+K)\mathfrak{A}_n-7\overline{A}_n\sqrt{\log(d)}
			\leq
			\P\del[1]{S_{n,W,S}\in H}-\gamma_{\Sigma_0}\del[0]{H},
		\end{equation*}
		so that, since the upper and lower bounds just obtained do not depend on~$H\in\mc{H}_d$,
		\begin{equation*}
			\sup_{H\in\mc{H}_d}\envert[1]{\P\del[1]{S_{n,W,S}\in H}-\gamma_{\Sigma_0}\del[0]{H}}
			\leq
			(1+K)\mathfrak{A}_n+7\overline{A}_n\sqrt{\log(d)}.
		\end{equation*}
		Finally, we show that if~$\log(dn)/n\leq 1$, then~$\overline{A}_n\sqrt{\log(d)}\leq C_4\mathfrak{A}_n$ for a suitable constant~$C_4$ depending only on $b_1,b_2,\lambda_1,\lambda_2,\lambda_1',\lambda_2'$, and~$m$, which yields~\eqref{eq:HDGaussstudentized}. To this end, note that up to the constant~$C_3$, $\overline{A}_n\sqrt{\log(d)}$ equals
		\begin{equation*}
		\sqrt{\log(d)\log(dn)}\sbr[3]{\overline{\eta}_n^{1-\frac{2}{m}}+\del[2]{\frac{\log(dn)}{n}}^{1-\frac{1}{(m/2)\wedge 2}}}
		=	
		\sqrt{\log(d)n}\sbr[3]{\overline{\eta}_n^{1-\frac{2}{m}}\del[2]{\frac{\log(dn)}{n}}^\frac{1}{2}+\del[2]{\frac{\log(dn)}{n}}^{\frac{3}{2}-\frac{1}{(m/2)\wedge 2}}}.
		\end{equation*}
		Because i)
		\begin{equation*}
			\overline{\eta}_n^{1-\frac{2}{m}}\del[2]{\frac{\log(dn)}{n}}^\frac{1}{2}
			\leq	
				\del[2]{\overline{\eta}_n\vee \frac{\log(dn)}{n}}^{\frac{3}{2}-\frac{2}{m}}
				\leq
				 \del[2]{\overline{\eta}_n\vee \frac{\log(dn)}{n}}^{1-\frac{1}{m}}
				 \leq
				 \overline{\eta}_n^{1-\frac{1}{m}}+\del[2]{\frac{\log(dn)}{n}}^{1-\frac{1}{m}},	
				\end{equation*}
		the second inequality using that~$\overline{\eta}_n\vee \log(dn)/n\leq 1$ and~$\frac{3}{2}-\frac{2}{m}\geq 1-\frac{1}{m}$, as well as ii)~$\frac{3}{2}-\frac{1}{(m/2)\wedge 2}\geq 1-\frac{1}{m}$, the right-hand side of the penultimate display is bounded from above by
		\begin{equation*}
			2\sqrt{\log(d)n}\sbr[3]{\overline{\eta}_n^{1-\frac{1}{m}}+\del[2]{\frac{\log(dn)}{n}}^{1-\frac{1}{m}}}
			\leq C_4\mathfrak{A}_n.
		\end{equation*}
		Here, the inequality used that the left-hand side is one of the constituent terms of~$\mathfrak{A}_n$ up to the factor 2. $\overline{A}_n\sqrt{\log(d)}\leq C_4\mathfrak{A}_n$ follows by adjusting~$C_4$. In case~$\log(dn)/n>1$,~\eqref{eq:HDGaussstudentized} follows from adjusting~$C$ there (if needed) because~$\rho_{n,W,S}\leq 1$.
		
		Finally, since~$\mathfrak{A}_n\to 0$ by~\eqref{eq:conv1},~\eqref{eq:conv2}, and~$R_{n,m}\to 0$ of Lemma~\ref{lem:covcompare}, it follows that~$\rho_{n,W,S}\to 0$. 
	\end{proof}

\section{Additional results}
\subsection{No universal barrier is possible without restrictions on~$\Gamma$}\label{sec:structure_needed} 

The following proposition shows that without any restrictions on the covariance functional~$\Gamma$, there exist (somewhat trivial) statistics~$T_n\in\mathfrak{T}_{n,d}$ for which the Gaussian approximation error in~\eqref{eq:error} converges to zero as~$n\to\infty$ irrespectively of how fast~$d=d_n$ grows with~$n$. Let~$\mathbf 1_d=(1,\ldots,1)'\in\R^d$ and~$\sigma_{m,1}^m=\int |x_1|^m P(dx)$ for~$m\in[2,\infty)$. 
\begin{proposition}\label{prop:arbitrarytarget}
Let~$\Gamma_*(P):=\sigma_{2,1}^2(P)\mathbf 1_d\mathbf 1_d'$ and for $(x_1,\ldots,x_n)\in(\R^d)^n$ set
\begin{equation*}
 T_n^*(x_1,\ldots,x_n)
 :=
 \del[3]{\frac{1}{\sqrt n}\sum_{i=1}^nx_{i,1}} \mathbf 1_d
\end{equation*}
Then, for every sequence $d=d_n$ of natural numbers, irrespective of its growth rate%
\begin{equation}\label{eq:convTstar}
		 \lim_{n\to\infty}
 \sup_{P\in\mc P_d^0(b_1,b_2,m)}\sup_{H\in\mc H_d}
 \envert[1]{P^{\otimes n}(T_n^*\in H)-\gamma_{\Gamma_*(P)}(H)}
 =
 0.
\end{equation}
\end{proposition}
Although~$\Gamma_*$ in Proposition~\ref{prop:arbitrarytarget} is a valid covariance functional (it is symmetric and positive semi-definite) it ``ignores'' all but the first coordinate of the~$X_i$. 
\begin{proof}[Proof of Proposition~\ref{prop:arbitrarytarget}]
Let $\Phi$ denote the cdf of the standard normal distribution and let
$\delta=(m-2)\wedge1\in(0,1]$. By the Berry--Esseen inequality as stated in (3) of~\cite{shevtsova2011asymptotically}
\begin{align}
 \sup_{t\in\R}
 \envert[3]{P^{\otimes n}\del[3]{
 \frac{\sum_{i=1}^nX_{i,1}}
 {\sqrt n}\leq t}-\Phi\del[2]{\frac{t}{\sigma_{2,1}(P)}}}
 &=
 \sup_{t\in\R}
 \envert[3]{P^{\otimes n}\del[3]{
 \frac{\sum_{i=1}^nX_{i,1}}
 {\sqrt n\sigma_{2,1}(P)}\leq t}-\Phi(t)}\notag\\
 &\leq
 \frac{C_\delta\sigma_{2+\delta,1}^{2+\delta}(P)}
 {\sigma_{2,1}^{2+\delta}(P)n^{\delta/2}}
 \leq 
 C_{m,b_1,b_2}n^{-\delta/2},                \label{eq:BE}
\end{align}
uniformly over $P\in\mc P_d^0(b_1,b_2,m)$ for a universal constant~$C_\delta>0$ depending only on~$\delta\in(0,1]$ and the final inequality following
from $2+\delta\leq m$, Lyapunov's inequality, and the moment bounds defining
$\mc P_d^0(b_1,b_2,m)$, cf.~\eqref{eq:dist_class0}. Next, for every~$H=\prod_{j=1}^d[a_j,b_j]\in\mc H_d$ note that
\begin{equation*}
	\cbr[1]{x\in\R: x\mathbf 1_d\in H}
	= 
	\cbr[2]{x\in\R:\max_{1\leq j\leq d}a_j\leq x\leq\min_{1\leq j\leq d}b_j}.
\end{equation*}
In case~$l:=\max_{1\leq j\leq d}a_j> \min_{1\leq j\leq d}b_j=:u$, the sets in the previous display are empty and it trivially holds that
\begin{equation*}
	P^{\otimes n}(T_n^*\in H)
	=
	P^{\otimes n}\del[3]{\frac{\sum_{i=1}^nX_{i,1}}{n^{1/2}}\in \emptyset}
	=
	0
	=
	\mathsf{N}_1\del[1]{0,\sigma_{2,1}^2(P)}(\emptyset)
	=
	\gamma_{\Gamma_*(P)}(H).
\end{equation*}
Thus, suppose that~$l\leq u$. By a standard argument dealing with the potential discontinuity of~$t\mapsto P^{\otimes n}(n^{-1/2}\sum_{i=1}^nX_{i,1}\leq t)$ at~$l$ and~\eqref{eq:BE}
\begin{align*}
\envert[2]{P^{\otimes n}(T_n^*\in H)-\gamma_{\Gamma_*(P)}(H)}
&=
\envert[3]{P^{\otimes n}\del[3]{\frac{\sum_{i=1}^nX_{i,1}}{n^{1/2}}\in [l,u]}
-
\sbr[2]{\Phi\del[2]{\frac{u}{\sigma_{2,1}(P)}}-\Phi\del[2]{\frac{l}{\sigma_{2,1}(P)}}}}\\
&\leq
2\cdot C_{m,b_1,b_2}n^{-\delta/2}
\end{align*}
In either case, the upper bound on~$\envert[0]{P^{\otimes n}(T_n^*\in H)-\gamma_{\Gamma_*(P)}(H)}$ does not depend on~$P\in \mc{P}_d^0(b_1,b_2,m)$,~$H\in\mc{H}_d$, or~$d$. Thus,~\eqref{eq:convTstar} follows.
\end{proof}

\subsection{Details of Remark~\ref{rem:conjectureselfnormalized}}\label{sec:remarkdetails}
Let~$U_n\in\R^d$ be the self-normalized sums as defined in~\eqref{eq:selfnorm} with $U_{n,j}=0$ when~$\sum_{i=1}^nX_{i,j}^2=0$ for~$j\in\cbr[0]{1,\hdots,d}$. Set~$m=3$, corresponding to~$\tau=1$ in~\cite{qiu2025self}. Their condition (2.1) is satisfied uniformly over~$P\in \mc P_d^0(b_1,b_2,3)$ because $\max_{1\leq j\leq d}E_P|X_{1,j}|^3/[E_P(X_{1,j}^2)]^{3/2}\leq b_2^3/b_1^3<\infty $. Let~$d=d_n=\lceil\exp(n^{1/7})\rceil$ so that~$\lim_{n\to\infty}\log(d)/n^{1/7}=1$. For~$t_n$ as in the statement of Theorem~\ref{thm:sharpStronger} in Appendix~\ref{sec:proofSharp} with~$\Gamma=\Sigma_0$, let
\begin{equation*}
a_{n,d}=a_{n,d}(P):=P^{\otimes n}\del[2]{\max_{1\leq j\leq d}U_{n,j}\leq t_n}\qquad \text{and}\qquad b_{n,d}=b_{n,d}(P):=\gamma_{\Sigma_0(P)}\del[1]{(-\infty,t_n]^d}	
\end{equation*}
so that
\begin{equation*}
\envert[3]{\frac{P^{\otimes n}\del[1]{\max_{1\leq j\leq d}U_{n,j}> t_n}}{1-\gamma_{\Sigma_0(P)}\del[1]{(-\infty,t_n]^d}}-1}
=
\envert[2]{\frac{1-a_{n,d}}{1-b_{n,d}}-1}
=
\envert[2]{\frac{b_{n,d}-a_{n,d}}{1-b_{n,d}}}
\geq
|a_{n,d}-b_{n,d}|.
\end{equation*}
Therefore, Theorem~\ref{thm:sharpStronger} in Appendix~\ref{sec:proofSharp}  implies, in particular, that 
\begin{equation*}
 \liminf_{n\to\infty}
 \sup_{P\in\mc P_d^0(b_1,b_2,3)}
 \envert[3]{\frac{P^{\otimes n}\del[1]{\max_{1\leq j\leq d}U_{n,j}> t_n}}{1-\gamma_{\Sigma_0(P)}\del[1]{(-\infty,t_n]^d}}-1}>0.
\end{equation*}
Because~$0\leq t_n\leq 2\sqrt{\log(d)}=O(n^{1/14})$, so that~$t_n/n^{\tau/(4+2\tau)}=t_n/n^{1/6}\to 0$, the thresholds~$t_n$ lie within the conjectured range. Hence, the Gaussian approximation conjectured in (2.16) of Remark 2.4 in~\cite{qiu2025self} is not true \emph{uniformly} over~$\mc P_d^0(b_1,b_2,3)$ on which the correlation matrix~$\Sigma_0(P)$ is unrestricted.\footnote{Positive approximation results are established in \cite{qiu2025self} under suitable restrictions on~$\Sigma_0$.}

\end{appendix}

\bibliographystyle{ecta} 
\bibliography{ref}		

@article{su2026sharp,
  title={Sharp High-Dimensional Self-Normalized Central Limit Theorems},
  author={Su, Zi-Yao and Zhao, Yi-Kun},
  journal={arXiv preprint arXiv:2611.31060},
  year={2026}
}

@article{wu2026counterexample,
  title={A Counterexample to Gaussian Approximation for Maxima of Self-Normalized Sums},
  author={Wu, Wei Biao},
  year={2026},
  journal={preprint}
}

@article{gao2026cubic,
  title={Cubic-Root Gaussian Approximation under Unrestricted Covariance},
  author={Gao, Zijun and Zhang, Weihan},
  journal={arXiv preprint arXiv:2608.30221},
  year={2026}
}

@article{CCKK26,
  title={Sharp High-Dimensional Gaussian Approximations
over Rectangles without Covariance Restrictions},
  author={Chernozhukov, Victor and Chetverikov, Denis and Kato, Kengo and Koike, Yuta},
  journal={arXiv preprint arXiv:2611.21474},
  year={2026}
  }

@article{LopesMasini26,
  title={High-Dimensional Gaussian and Bootstrap Approximation: Near~$1/\sqrt{n}$ Rates with General Covariance},
  author={Lopes, Miles and Masini, Ricardo},
  journal={arXiv preprint arXiv:2611.17923},
  year={2026}
}

@article{chang2025berry,
  title={The Berry-Esseen Bound for High-dimensional Self-normalized Sums},
  author={Chang, Woonyoung and Takatsu, Kenta and Urban, Konrad and Kuchibhotla, Arun Kumar},
  journal={arXiv preprint arXiv:2501.08979},
  year={2025}
}

@article{shevtsova2011asymptotically,
  title={On the asymptotically exact constants in the Berry--Esseen--Katz inequality},
  author={Shevtsova, IG},
  journal={Theory of Probability \& Its Applications},
  volume={55},
  number={2},
  pages={225--252},
  year={2011},
  publisher={SIAM}
}

@book{liese2008statistical,
  author    = {Liese, Friedrich and Miescke, Klaus-J.},
  title     = {Statistical Decision Theory:
               Estimation, Testing, and Selection},
  series    = {Springer Series in Statistics},
  publisher = {Springer},
  year      = {2008}
}

@article{aronow2026optimal,
  title={Optimal central limit theorem for bounded random variables in high dimensions},
  author={Aronow, PM and Lopatto, Patrick},
  journal={arXiv preprint arXiv:2609.07067},
  year={2026}
}

@article{kock2025high,
  title={High-dimensional {Gaussian} and bootstrap approximations for robust means},
  author={Kock, Anders Bredahl and Preinerstorfer, David},
  journal={Annals of Applied Probability (forthcoming) and arXiv preprint arXiv:2504.08435},
  year={2025},
  note={Version of June 29, 2026},
  url={https://arxiv.org/abs/2504.08435}
}

@article{das2025necessary,
  title={Necessary and sufficient conditions for high dimensional Central Limit Theorem under moment conditions},
  author={Das, Debraj and Lahiri, Soumendra},
  journal={arXiv preprint arXiv:2512.22312},
  year={2025}
}

@article{das2024central,
  title={Central limit theorem and near classical Berry-Esseen rate for self normalized sums in high dimensions},
  author={Das, Debraj},
  journal={Bernoulli},
  volume={30},
  number={1},
  pages={278--303},
  year={2024},
  publisher={Bernoulli Society for Mathematical Statistics and Probability}
}

@article{qiu2025self,
  title={Self-normalized Cram{\'e}r type moderate deviation theorem for Gaussian approximation},
  author={Qiu, Jingkun and Chen, Song Xi and Shao, Qi-Man},
  journal={The Annals of Statistics},
  volume={53},
  number={3},
  pages={1319--1346},
  year={2025},
  publisher={Institute of Mathematical Statistics}
}

@article{Corless1996,
  author  = {Corless, R. M. and Gonnet, G. H. and Hare, D. E. G. and Jeffrey, D. J. and Knuth, D. E.},
  title   = {On the {L}ambert {$W$} function},
  journal = {Advances in Computational Mathematics},
  year    = {1996},
  volume  = {5},
  number  = {1},
  pages   = {329--359},
  doi     = {10.1007/BF02124750}
}

@article{cstudd,
 ISSN = {00034851, 21688990},
 URL = {http://www.jstor.org/stable/2239660},
 author = {Y. S. Chow and W. J. Studden},
 journal = {The Annals of Mathematical Statistics},
 number = {3},
 pages = {1106--1108},
 publisher = {Institute of Mathematical Statistics},
 title = {Monotonicity of the variance under truncation and variations of Jensen's inequality},
 urldate = {2025-04-11},
 volume = {40},
 year = {1969}
}

@article{Wins1,
  title={Winsorized mean estimation with heavy tails and adversarial contamination},
  author={Kock, Anders Bredahl and Preinerstorfer, David},
  journal={Electronic Journal of Statistics},
  volume={20},
  number={2},
  pages={3225--3254},
  year={2026},
  publisher={The Institute of Mathematical Statistics and the Bernoulli Society}
}

@article{liu2024robust,
  title={Robust max statistics for high-dimensional inference},
  author={Liu, Mingshuo and Lopes, Miles E},
  journal={arXiv preprint arXiv:2409.16683},
  year={2024}
}

@article{koike2024high,
  title={High-dimensional bootstrap and asymptotic expansion},
  author={Koike, Yuta},
  journal={Probability Theory and Related Fields},
  pages={1051--1125},
  volume={195},
  year={2026},
  publisher={Springer}
}

@article{fang2023high,
  title={High-dimensional Central Limit Theorems by Stein's method in the degenerate case},
  author={Fang, Xiao and Koike, Yuta and Liu, Song-Hao and Zhao, Yi-Kun},
  journal={arXiv preprint arXiv:2305.17365},
  year={2023}
}

@article{fang2021high,
  title={High-dimensional central limit theorems by Stein's method},
  author={Fang, Xiao and Koike, Yuta},
  journal={Annals of Applied Probability},
  volume={31},
  number={4},
  pages={1660--1686},
  year={2021}
}

@article{kock2024remark,
  title={A remark on moment-dependent phase transitions in high-dimensional Gaussian approximations},
  author={Kock, Anders Bredahl and Preinerstorfer, David},
  journal={Statistics \& Probability Letters},
  volume={211},
  pages={110149},
  year={2024},
  publisher={Elsevier}
}

@article{das2021central,
  title={Central Limit Theorem in high dimensions: The optimal bound on dimension growth rate},
  author={Das, Debraj and Lahiri, Soumendra},
  journal={Transactions of the American Mathematical Society},
  volume={374},
  number={10},
  pages={6991--7009},
  year={2021}
}

@article{zhangwu2017,
author = {Danna Zhang and Wei Biao Wu},
title = {{Gaussian approximation for high dimensional time series}},
volume = {45},
journal = {Annals of Statistics},
number = {5},
publisher = {Institute of Mathematical Statistics},
pages = {1895--1919},
year = {2017},
doi = {10.1214/16-AOS1512},
URL = {https://doi.org/10.1214/16-AOS1512}
}

@book{leadbetter,
  title={Extremes and Related Properties of Random Sequences and Processes},
  author={Leadbetter, Malcolm R and Lindgren, Georg and Rootz{\'e}n, Holger},
  year={1983},
  publisher={Springer: New York}
}

@article{lopes2020bootstrapping,
  title={Bootstrapping max statistics in high dimensions: Near-parametric rates under weak variance decay and application to functional and multinomial data},
  author={Lopes, Miles and Lin, Zhenhua and M{\"u}ller, Hans-Georg},
  volume={48},
  journal={Annals of Statistics},
  issue={2},
  pages={1214--1229},
  year={2020}
}

@article{kuchibhotla2021high,
  title={High-dimensional CLT: Improvements, non-uniform extensions and large deviations},
  author={Kuchibhotla, Arun Kumar and Mukherjee, Somabha and Banerjee, Debapratim},
  volume={27},
  journal={Bernoulli},
  issue={1},
  pages={192--217},
  year={2021}
}

@article{chernozhukov2023nearly,
  title={Nearly optimal central limit theorem and bootstrap approximations in high dimensions},
  author={Chernozhukov, Victor and Chetverikov, Denis and Koike, Yuta},
  journal={Annals of Applied Probability},
  volume={33},
  number={3},
  pages={2374--2425},
  year={2023},
  publisher={Institute of Mathematical Statistics}
}

@article{koike2021notes,
  title={Notes on the dimension dependence in high-dimensional central limit theorems for hyperrectangles},
  author={Koike, Yuta},
  journal={Japanese Journal of Statistics and Data Science},
  volume={4},
  pages={257--297},
  year={2021},
  publisher={Springer}
}

@article{kuchibhotla2020high,
  title={High-dimensional CLT for sums of non-degenerate random vectors: $n^{-1/2}$-rate},
  author={Kuchibhotla, Arun Kumar and Rinaldo, Alessandro},
  journal={arXiv preprint arXiv:2009.13673},
  year={2020}
}

@article{lopes2022central,
  title={Central limit theorem and bootstrap approximation in high dimensions: Near 1/n rates via implicit smoothing},
  author={Lopes, Miles},
  journal={Annals of Statistics},
  volume={50},
  number={5},
  pages={2492--2513},
  year={2022},
  publisher={Institute of Mathematical Statistics}
}

@article{deng2020beyond,
  title={Beyond Gaussian approximation: Bootstrap for maxima of sums of independent random vectors},
  journal={Annals of Statistics},
  pages={3643--3671},
  volume={48},
  issue={6},
  author={Deng, Hang and Zhang, Cun-Hui},
  year={2020}
}

@article{chernozhukov2013gaussian,
  title={Gaussian approximations and multiplier bootstrap for maxima of sums of high-dimensional random vectors},
  author={Chernozhukov, Victor and Chetverikov, Denis and Kato, Kengo},
  journal={Annals of Statistics},
  volume={41},
  issue={6},
  pages={2786--2819},
  year={2013}
}

@article{chernozhuokov2022improved,
  title={Improved central limit theorem and bootstrap approximations in high dimensions},
  author={Chernozhukov, Victor and Chetverikov, Denis and Kato, Kengo and Koike, Yuta},
  journal={Annals of Statistics},
  volume={50},
  number={5},
  pages={2562--2586},
  year={2022},
  publisher={Institute of Mathematical Statistics}
}

@article{chernozhukov2023high,
  title={High-dimensional data bootstrap},
  author={Chernozhukov, Victor and Chetverikov, Denis and Kato, Kengo and Koike, Yuta},
  journal={Annual Review of Statistics and Its Application},
  volume={10},
  pages={427--449},
  year={2023}
}

@article{chernozhukov2017central,
  title={Central limit theorems and bootstrap in high dimensions},
  author={Chernozhukov, Victor and Chetverikov, Denis and Kato, Kengo},
  journal={The Annals of Probability},
  pages={2309--2352},
  volume={45},
  year={2017}
}

@article{chernozhukov2017detailed,
  title={Detailed proof of Nazarov's inequality},
  author={Chernozhukov, Victor and Chetverikov, Denis and Kato, Kengo},
  journal={arXiv preprint arXiv:1711.10696},
  year={2017}
}

@article{resende2024robust,
  title={Robust high-dimensional Gaussian and bootstrap approximations for trimmed sample means},
  author={Resende, Lucas},
  journal={arXiv preprint arXiv:2410.22085v2 of 28 July 2026},
  year={2024}
}

\end{document}